\documentclass[english]{artjlt}

\usepackage{latexsym}
\usepackage{amsfonts}
\usepackage{amssymb}
\usepackage{amsmath}
\usepackage{enumerate}
\usepackage{cite}
\usepackage{hyperref}  
\title{Maximal Subgroups of Compact Lie Groups}

\author{Fernando Antoneli\thanks{Work partially 
supported by CNPq (Conselho Nacional de Desenvolvimento 
Cient\'{\i}fico e Tecno\-l\'ogico), Brazil.} ,
Michael Forger\footnotemark[1] \;
and 
Paola Gaviria}

\lastname{Antoneli, Forger and Gaviria}

\msc{22E15}

\keywords{Lie groups, Lie algebras, Compact groups, Maximal subgroups}

\address{Fernando Antoneli \\
Escola Paulista de Medicina \\
Universidade Federal de S\~ao Paulo \\
S\~ao Paulo, SP \\
BR 04039-062 -- Brazil \\
fernando.antoneli@unifesp.br
}

\address{Michael Forger \\
Instituto de Matem\'atica e Estat\'{\i}stica \\
Universidade de S\~ao Paulo \\
S\~ao Paulo, SP, PO Box 66281 \\
BR 05315-970 -- Brazil  \\
forger@ime.usp.br
}

\address{Paola Gaviria \\
Instituto de Matem\'atica e Estat\'{\i}stica \\
Universidade de S\~ao Paulo \\
S\~ao Paulo, SP, PO Box 66281 \\
BR 05315-970 -- Brazil \\
pgaviria@ime.usp.br
}

\renewcommand{\firstheadline}{\begin{minipage}[t]{9cm}
    Journal of Lie Theory\hfill

    Volume {\bf 22}\ (2012)\ 949--1031\hfill

    \copyright\  2012\ Heldermann Verlag\hfill
    \end{minipage}\hfill
      }

\newcommand{\Z}{\mathbb{Z}}

\newcommand{\R}{\mathbb{R}}
\newcommand{\C}{\mathbb{C}}

\DeclareMathOperator{\tr}{tr}

\newcommand{\SPG}{\mathit{Sp\>\!}}

\begin{document}

\maketitle

\begin{abstract}
\noindent
 This report aims at giving a general overview on the classification of
 the maximal subgroups of compact Lie groups (not necessarily connected).
 \linebreak
 In the first part, it is shown that these fall naturally into three types:
 (1)~those of trivial type, which are simply defined as inverse images
 of maximal subgroups of the corresponding component group under the
 canonical projection and whose classification constitutes a problem
 in finite group theory, (2)~those of normal type, whose connected
 one-component is a normal subgroup, and (3)~those of normalizer
 type, which are the normalizers of their own connected one-component.
 It is also shown how to reduce the classification of maximal subgroups
 of the last two types to: (2)~the classification of the finite maximal
 $\Sigma$-invariant subgroups of centerfree connected compact simple
 Lie groups and (3)~the classification of the $\Sigma$-primitive
 sub\-algebras of compact simple Lie algebras, where $\Sigma$ is
 a subgroup of the corresponding outer automorphism group.
 In the second part, we explicitly compute the normalizers of
 the primitive subalgebras of the compact classical Lie algebras
 (in the corresponding classical groups), thus arriving at the
 complete classification of all (non-discrete) maximal subgroups
 of the compact classical Lie groups.
\end{abstract}

\vspace{5mm}

\section{Introduction}

In this paper, we address an important problem from the theory of
Lie groups, namely that of classifying their maximal subgroups.%
\footnote{In this paper, the term ``maximal \ldots\ subgroup'' will
always mean ``maximal \ldots\ closed subgroup'', where \ldots\ stands
for a string of other possible adjectives. A more detailed discussion,
with precise definitions, can be found in Section~2, for the case where
this string is empty, and in Section~3, for the case where this string
represents the expression ``$\Gamma$-invariant'' or ``$\Sigma$-invariant''.}
This is closely related (though not completely equivalent) to
another classical problem from the \mbox{theory} of Lie groups,
namely that of classifying their primitive actions, which goes
all the way back to Sophus Lie and has over decades received
a great deal of attention in the literature: see Dynkin~%
\cite{Dy1,Dy2,Ti}, Golubitsky et al.~\cite{Go,GR},
Chekalov~\cite{Ch} and Komrakov~\cite{Ko1,Ko2,Ko3}.

Of course, the same two problems and, in particular, the notion of a
primi\-tive action can also be formulated in the context of abstract groups.
A transitive action of an abstract group $G$ on a set $X$ is said to
be \emph{primitive} if the only $G$-invariant equivalence relations
$\, R \subset X \times X \,$ are the trivial ones: $\, R = X \times X$
and $\, R = \{ (x,x) \,|\; x \in X \}$.\footnote{Often it is also 
assumed that the action is effective, but this is not essential.}%
\addtocounter{footnote}{-1}
This is equivalent to $X$ being a homogeneous space $G/H$ where the
stability group~$H$ is a maximal subgroup of~$G$, so in this context
the two problems are completely equivalent.
For Lie groups, the corresponding concept is defined somewhat differently.
A transitive action of a Lie group $G$ \linebreak on a manifold $M$ is said
to be \emph{(Lie) primitive} if $M$ admits no $G$-invariant \mbox{foliation}
with leaves of positive dimension smaller than $\dim M$.\footnotemark\,\,
This is equivalent \linebreak to $M$ being a homogeneous space $G/H$
where the stability group~$H$ satisfies the \mbox{following} weaker
maximality condition: for any Lie subgroup $\tilde{H}$ of~$G$ such
that $\, H \subset \tilde{H} \subset G$, $\tilde{H}_0 = H_0 \,$ or
$\, \tilde{H}_0 = G_0$, where the index~$0$ denotes taking the
connected one-component.
Clearly every action of a Lie group which is primitive as the action
of an abstract group is also Lie primitive, but not conversely.

\setcounter{page}{950}

The infinitesimal version of both problems leads to the quest for a
classification of the maximal subalgebras of complex semisimple Lie
algebras.
This was completely solved in the early 1950s by Morozov~\cite{Mo1,Mo2},
Karpelevich~\cite{Ka}, Borel and de~Siebenthal~\cite{BS} and Dynkin~%
\cite{Dy1,Dy2}; see also the Bourbaki review by Tits~\cite{Ti}.
However, the classification of maximal subalgebras of Lie algebras
(or at least of certain classes of Lie algebras such as compact or
semisimple or, more generally, reductive Lie algebras) only
provides a classification of maximal \emph{connected} subgroups
of \emph{connected} Lie groups (within the corresponding class):
dealing with the discrete parts is quite another story.
In fact, although maximal subalgebras of a Lie algebra do give rise
to maximal subgroups of any corresponding connected Lie group
(namely, as we shall see, by taking their normalizer under the
adjoint representation), it is \emph{not} true that maximal sub%
groups are necessarily associated with maximal subalgebras: an
extreme counterexample is provided by the trivial subalgebra
$\{0\}$, corresponding to discrete maximal subgroups.
Thus we must face the question as to what is the class
of subalgebras corresponding to the maximal subgroups.

The study of the global version including non-connected groups was
initiated in the 1970s by Golubitsky~\cite{Go}, who determined the
maximal rank primitive subgroups of the classical complex Lie groups.
Soon after, his work was extended to the exceptional complex
Lie groups by Golubitsky and Rothschild~\cite{GR}.
Finally, this classification was completed by including primitive
subgroups of any rank, first for the classical complex Lie groups
by Chekalov~\cite{Ch} and later for the exceptional complex Lie
groups by Komrakov~\cite{Ko3}.

It is clear from this brief summary that most of the efforts
were centered around the notion of primitive actions, leaving
that of maximal subgroups as an auxiliary concept.
It is not difficult to see that when the ambient group is a
simple Lie group, the non-discrete maximal subgroups can be
obtained from the primitive subgroups: the former are precisely
the normalizers of the connected one-components of the latter.
However, this correspondence is less clear when the ambient
group is not simple.
In~fact, it is shown in~\cite[p.~179]{Go} that if a semisimple
Lie group has more than three simple factors, then it does not
admit any primitive subgroups, while it is obvious that it
does admit non-discrete maximal subgroups.
Thus although the results mentioned above constitute important
steps towards a full classification of maximal subgroups, some
ingredients are still missing.

The main purpose of the present paper is to fill these gaps
and to provide a comprehensive treatment of the subject,
which does not seem to be available in the literature.
In particular, our approach includes the study of maximal
subgroups of Lie groups which are not necessarily connected:
a problem which apparently has never before been
investigated systematically and which cannot be neglected
if one wants to arrive at a conceptually consistent picture.
Indeed, such a picture can only emerge if one requires the
maximal subgroups considered to belong to the same category
as their ambient groups: either they should both be assumed
to be connected or else they should both be allowed to have
a non-trivial component group.
But apart from logical coherence and completeness, such an
approach is also of practical use; in fact, it is imperative
if one wants to iterate the process in order to construct
descending chains of subgroups in which each subgroup is
maximal in the previous one: this is an important
ingredient in studies of symmetry breaking~\cite{HHF}.
On the other hand, one cannot of course hope for a completely
general solution: some restriction on the type of Lie groups
involved will certainly have to be imposed.
Not surprisingly, an adequate category for which a complete
theory can be developed is that of compact Lie groups, but it
can be verified with little effort that practically all results
obtained in this context continue to hold within the larger
category of reductive Lie groups, which contains that of
semisimple Lie groups~-- compact as well as non-compact.

The paper is divided into two parts: the first (Sections 2-5) is
devoted to the general theory whereas the second (Sections 6-8)
deals with the classification problem for the classical groups.
In~Section~2, we begin by specifying the precise mathematical
setup and, in particular, by giving an exact definition of the
term ``maximal subgroup''.
We~also collect some elementary facts about the action of a
Lie group on its connected one-component and on its
Lie algebra by automorphisms and about normalizers.
In~Section~3, we analyze the relation between subgroups of a Lie
group~$G$ and subgroups of its connected one-component~$G_0$ (see
Proposition~\ref{prop:invsg}), which leads us to introduce the
important and useful concept of a maximal $\Gamma$-invariant
subgroup of~$G_0$, where $\, \Gamma = G/G_0 \,$ is the component
group of~$G$.
Next, we prove as our first main \mbox{theorem} that maximal
subgroups of Lie groups fall into three distinct types which
we shall refer to as the trivial type, the normal type and
the normalizer type (see Theorem~\ref{theo:maxsg1}).
Briefly, for any Lie group~$G$ with connected one-component~$G_0$
and component group $\, \Gamma = G/G_0$, a maximal subgroup of~$G$
of trivial type%
\footnote{This terminology, borrowed from the theory of finite groups
(see Remark~\ref{rmk:maxsg3}), is perhaps unfortunate, and the authors
are open to suggestions for a better one.}
is obtained as the union of the connected components of~$G$ labelled by
the elements of a maximal subgroup of~$\Gamma$ (so this type can only
exist if $G$ is not connected), whereas any other maximal subgroup $M$
of~$G$ must meet every connected component of~$G$ and is classified
according to what is the normalizer $N_G(M_0)$ in~$G$ of its connected
one-component $M_0$, since maximality of~$M$ leaves only two options:
\begin{itemize}
 \item $N_G(M_0) = G$: this is called the normal type since $M_0$ is
       a normal subgroup of~$G$,
 \item $N_G(M_0) = M$: this is called the normalizer type since $M$
       itself is the normalizer of its own connected one-component,
       or equivalently, of its own Lie algebra.
\end{itemize}
In~Section~4, we introduce the concept of a $\Sigma$-quasiprimitive
subalgebra of a Lie algebra~$\mathfrak{g}$, containing as a special
case that of a $\Sigma$-primitive subalgebra, where $\Sigma$ is any
given subgroup of the outer automorphism group $\mathrm{Out}%
(\mathfrak{g})$ of~$\mathfrak{g}$.
In~Section~5, we turn to reductive Lie groups and prove first that,
at least in this context, the $\Sigma$-quasiprimitive subalgebras
of~$\mathfrak{g}$ correspond precisely to the the maximal subgroups
of~$G$ of normal type or of normalizer type, where $G$ is any Lie
group with Lie algebra~$\mathfrak{g}$ whose component group $\Gamma$
projects to~$\Sigma$ under the natural homomorphism from~$\Gamma$
to~$\mathrm{Out}(\mathfrak{g})$ induced by the adjoint representation
of~$G$ (see Theorem~\ref{theo:maxsg2}).%
\footnote{In this situation, we shall use the terms ``$\Sigma$-(quasi)%
primitive'' and ``$\Gamma$-(quasi)primitive'' inter\-changeably. Of course,
the prefix will be omitted as soon as $\Sigma$ is trivial.}
Subsequently, we show how the classification of the maximal subgroups
of reductive Lie groups can, in a sequence of steps, be reduced to that
of the maximal $\Sigma$-invariant subgroups of centerfree connected
simple Lie groups, where $\Sigma$ runs through the subgroups of the
respective outer automorphism group.
Together with the results obtained in the preceding sections, this procedure
reduces the classification of the maximal subgroups of a general compact Lie
group~$G$ with connected one-component~$G_0$, component group $\, \Gamma
= G/G_0 \,$ and Lie algebra~$\mathfrak{g}$ to
\begin{enumerate}[(1)]
\item the classification of the maximal subgroups of $\Gamma$, for the
      trivial type,
\item the classification of the discrete (hence finite) maximal
      $\Sigma_s$-invariant subgroups of~$G_s$, for the normal type,
\item the classification of the $\Sigma_s$-primitive subalgebras
      of~$\mathfrak{g}_s$, whose normalizers provide the non-discrete
      maximal $\Sigma_s$-invariant subgroups of~$G_s$, for the
      normalizer type,
\end{enumerate}
where $G_s$ is any one of the simple factors of the quotient of~$G_0$ by
its center $Z(G_0)$, $\mathfrak{g}_s$ is its Lie algebra and $\Sigma_s$ is
any subgroup of the outer automorphism group $\mathrm{Out}(\mathfrak{g}_s)$
of~$\mathfrak{g}_s$.
Note that the first item is part of the problem of classifying the maximal
subgroups of finite groups, which constitutes an important issue in the
theory of finite groups that has been vigorously investigated in the last
decade; see, for example,~\cite{AS,KL,LPS,LS2,No,NW}.
On the other hand, the second item (even with the simplifying assumption
that $G$ is connected and hence $\, \Gamma = \{1\}$, or more generally,
that $G$ acts on~$G_0$ by inner automorphisms and hence $\, \Sigma =
\{1\}$) leads directly to the problem of classifying the finite maximal
subgroups of compact connected simple Lie groups.
For the classical Lie groups, it is known that this may be reduced to
a problem in the representation theory of finite groups~\cite{GR1},
whose solution in the simplest case $A_1$ is stated in Example~%
\ref{ex:SO3} below, whereas for the exceptional Lie groups, it
constitutes a highly non-trivial problem which has been the subject
of recent (and still ongoing) research by experts in the theory of
finite groups~\cite{Gr,GR1,GR2}.

In view of this situation, we shall in the second part of this paper
deal with the third of the above items: the classification of the
maximal subgroups of compact simple Lie groups of normalizer type,
or equivalently, of the $\Sigma$-primitive subalgebras of compact
simple Lie algebras $\mathfrak{g}$, where $\Sigma$ is a subgroup of
the (finite) outer automorphism group $\mathrm{Out}(\mathfrak{g})$
of~$\mathfrak{g}$.
Note that with only one exception, $\mathrm{Out}(\mathfrak{g})$ (and
hence $\Sigma$) is either trivial or equal to $\mathbb{Z}_2$; in fact,
\[
 \mathrm{Out}(\mathfrak{g})~
 =~\left\{ \begin{array}{ccc}
               \{1\}     & \mbox{for} & \mbox{$A_1$,
                                              $B_n$ ($n \geqslant 2$),
                                              $C_n$ ($n \geqslant 3$),
                                              $E_7$, $E_8$, $F_4$, $G_2$}
            \\[2mm]
            \mathbb{Z}_2 & \mbox{for} & \mbox{$A_n$ ($n \geqslant 2$),
                                              $D_n$ ($n \geqslant 5$),
                                              $E_6$}
            \\[2mm]
                S_3      & \mbox{for} & \mbox{$D_4$}
           \end{array} \right\}~.
\]
Therefore, our strategy will be to first deal with the situation where
$\Sigma$ is trivial and leave the determination of $\Sigma$-primitive
subalgebras for non-trivial $\Sigma$ to a second step; this will be
greatly facilitated by the fact that any primitive subalgebra which
is also $\Sigma$-invariant is automatically $\Sigma$-primitive.
In~Section~6, we discuss how the classification problem for
primitive subalgebras (and analogously, for maxi\-mal subalgebras)
of compact simple Lie algebras can be translated to the complex
setting, where it becomes a classification problem for primitive
reductive sub\-algebras (and analogously, for maximal reductive
subalgebras) of complex simple Lie algebras.
In~Section~7, we use this method to translate the existing classification
of primitive reductive subalgebras of the complex classical Lie algebras~%
\cite{Ch,Dy1,Dy2,Go} to the compact setting and, in a second step, extend
it to a classification of $\Sigma$-primitive subalgebras of all compact
classical Lie algebras except~$\mathfrak{so}(8)$.
Finally, in Section~8, we address the task of computing the corresponding
normalizers in the respective classical groups~-- a task that, strangely
enough, is largely neglected in the existing literature.
We begin by assembling a few general tools that are needed in these
computations and then proceed to a case by case treatment of each
of the classical groups $SU(n)$, $SO(n)$ and $\SPG(n)$.
The results are collected in Tables \ref{tab:SUG}--\ref{tab:SPG}.
It should be pointed out that we restrict ourselves to listing maximal
subgroups whose connected one-component is not simple since the others
are explicitly given in terms of irreducible representations (see
Theorem~\ref{theo:DYNEXC}), with a ``list of exceptions'' for which
we refer the reader to the original paper by Dynkin~\cite{Dy2}.
This result will finally fill in a gap left by Dynkin~\cite[p.~247]{Dy2},
promising that the computation of the normalizers of the maximal connected
subgroups would appear in a separate article, which has apparently never
been published.

\section{General setup and basic definitions}

In order to specify the precise mathematical setup for the problem
to be treated in this paper, we must first of all give an exact
definition of what we mean by a ``maximal subgroup'' of a Lie
group and, analogously, a ``maximal subalgebra'' of a Lie algebra.
For Lie algebras, this definition is standard.
\begin{Definition} \label{def:maxsa}
 Let $\mathfrak{g}$ be a Lie algebra.
 A \emph{maximal subalgebra} of~$\mathfrak{g}$ is a proper sub\-algebra
 $\mathfrak{m}$ of~$\mathfrak{g}$ such that if $\tilde{\mathfrak{m}}$
 is any subalgebra of~$\mathfrak{g}$ with $\, \mathfrak{m}
 \subset \tilde{\mathfrak{m}} \subset \mathfrak{g}$, then
 $\, \tilde{\mathfrak{m}} = \mathfrak{m} \,$ or
 $\, \tilde{\mathfrak{m}} = \mathfrak{g}$.
 The same terminology applies when the expression ``subalgebra''
 is everywhere replaced by the expression ``ideal''.
\end{Definition}
\noindent
In the case of Lie groups, however, the corresponding definition
is ambiguous to the extent that it depends on the specific lattice
$\mathcal{L}$ of subgroups in which one chooses to search for
maximal elements.
The general convention on terminology here is that a maximal ... sub%
group of a Lie group~$G$ is a maximal element in the lattice \linebreak
of ... subgroups of~$G$, where ... is a string of adjectives that defines
$\mathcal{L}$~-- adjectives such as ``closed'' or ``Lie'' or ``connected''
or ``discrete'' or ``finite'', etc.%
\footnote{Abstract subgroups of Lie groups which are not Lie
subgroups will be disregarded.}
Explicitly, we have, for example:
\begin{Definition} \label{def:maxsg}
 Let $G$ be a Lie group.
 A \emph{maximal closed subgroup} of~$G$ is a proper closed subgroup
 $M$ of~$G$ such that if $\tilde{M}$ is any closed subgroup of~$G$ with
 $\, M \subset \tilde{M} \subset G$, then $\, \tilde{M} = M \,$ or
 $\, \tilde{M} = G$.
 The same terminology applies when the expression ``closed subgroup''
 is everywhere replaced by the expression ``closed normal subgroup''.
\end{Definition}
\noindent
Similarly, we may define the concept of a \emph{maximal Lie subgroup},
as well as those of a \emph{maximal connected closed/Lie subgroup} and
of a \emph{maximal discrete closed/Lie subgroup}, of a Lie group~$G$,
repeating Definition~\ref{def:maxsg} with the word ``closed''
every\-where replaced by the word ``Lie'' or by the expressions
``connected closed/Lie'' and ``discrete closed/Lie'', respectively.%
\footnote{Connected Lie subgroups are also called analytic subgroups.
Discrete closed/Lie subgroups are simply closed/Lie subgroups with
trivial Lie algebra. In the closed case, these are discrete in the
topology of~$G$, whereas in the Lie case, they are only discrete in
their own topology but not in that of~$G$; quite to the contrary,
they may very well be dense in~$G$.}

For connected Lie groups, it is natural to consider only connected
(but not necessarily closed) Lie subgroups, since the Lie correspondence~%
\cite[p.~47]{Kn} establishes a bijection between the lattice of connected
Lie subgroups of a connected Lie group~$G_0$ and the lattice of subalgebras
of its Lie algebra~$\mathfrak{g}$: hence it also establishes a bijection
between maximal connected Lie subgroups of~$G_0$ and maximal subalgebras
of~$\mathfrak{g}$.
This implies that maximal connected Lie subgroups always exist and
that any connected Lie subgroup is contained in a maximal one; in
fact, it can always be realized as the last member in a finite
descending chain of connected Lie subgroups where each is maximal
in the preceding one.%
\footnote{The same favorable situation prevails for finite groups:
they always admit maximal subgroups and any subgroup is contained
in a maximal one: in fact, it can always be realized as the last
member in a finite descending chain of subgroups where each is
maximal in the preceding one.}


On the other hand, for compact Lie groups, it is natural to consider
only closed (but not necessarily connected) subgroups, since these are
again compact Lie groups.
Unfortunately, the Lie correspondence is almost entirely lost when we
exclude Lie subgroups which are not closed but allow for Lie subgroups
which are not connected.
As a result, we cannot even guarantee that maximal closed subgroups
exist at all, nor that a given closed subgroup is contained in a
maximal one. \linebreak
In fact, there is an elementary general obstruction:
\begin{proposition} \label{prop:nomsg}
 A connected abelian Lie group contains no maximal closed subgroups.
\end{proposition}
\begin{proof}
 To begin with, note that it suffices to show that a connected abelian
 Lie group contains no maximal discrete closed subgroups.
 (Indeed, if $G$ is a connected abelian Lie group and $H$ is a closed
 subgroup of~$G$ with connected one-component~$H_0$, then $G/H_0$ is
 a connected abelian Lie group and $H/H_0$ is a discrete closed sub%
 group of~$G/H_0$, so if $H$ were maximal among the closed subgroups
 of~$G$, $H/H_0$ would have to be maximal among the discrete closed
 subgroups of~$G/H_0$.)
 To prove this, suppose that $G$ is a connected abelian Lie group,
 so
 $\, G = \mathbb{R}^n/\Lambda \cong \mathbb{T}^p \times \mathbb{R}^q \,$
 where $\Lambda$ is some lattice in~$\mathbb{R}^n$, $\mathbb{T}$ denotes
 the unit circle (one-dimensional torus) and $p$ is the dimension of the
 subspace of~$\mathbb{R}^n$ generated by $\Lambda$ \cite[Theorem~6.20,
 p.~155]{SW}.
 Now it is well known that the discrete closed subgroups of~%
 $\mathbb{R}^n$ are precisely the lattices in~$\mathbb{R}^n$~%
 \cite[Lemma~6.18, p.~155]{SW}, so it is clear that discrete
 closed subgroups of~$G$ are of the form $L/\Lambda$ where $L$ is
 a lattice  in~$\mathbb{R}^n$ containing $\Lambda$ as a sublattice.
 But for any of these lattices~$L$, there is always a bigger one
 in which it is contained.
\end{proof}
\noindent
Fortunately, the situation is quite different for nonabelian Lie
groups~-- otherwise, the present paper would be pointless.
To see this, let us look at a simple example:
\begin{example} \label{ex:SO3}
 Consider the group $\, G = SO(3) \,$ of all rotations in~$\mathbb{R}^3$.
 Then up to conjugacy, the non-discrete proper closed subgroups (and in
 fact, the non-discrete proper Lie subgroups) of~$G$ are
 \[
  \begin{array}{c}
   SO(2)~=~\mbox{the two-dimensional rotation group}~, \\[1mm]
   O(2)~=~\mbox{the two-dimensional orthogonal group}~,
  \end{array}
 \]
 whereas the discrete closed (and hence finite) subgroups of~$G$ are
 \[
  \begin{array}{c}
   \mathbb{Z}_n~=~\mbox{the cyclic group of order $n$}~, \\[1mm]
   \mathrm{D}_n~=~\mbox{the dihedral group of order $2n$}~, \\[1mm]
   \mathrm{T}~=~\mbox{the tetrahedral group}~\cong~A_4~, \\[1mm]
   \mathrm{O}~=~\mbox{the cubic/octahedral group}~\cong~S_4~, \\[1mm]
   \mathrm{I}~=~\mbox{the dodecahedral/icosahedral group}~\cong~A_5~.
  \end{array}
 \]
 (The first statement follows by noting that if $H$ is a non-discrete
 proper Lie subgroup of~$\, G = SO(3) \,$ with connected one-component~%
 $H_0$, then $H_0$ must be equal to its maximal torus $\, T = SO(2) \,$
 and $H$ itself must be contained in its normalizer $\, N_G(T) = O(2)$.%
 \footnote{Here, we use the fact that any Lie subgroup must normalize
 its own connected one-component; this general feature will be used
 extensively in what follows.}
 But $\, N_G(T)/T = \mathbb{Z}_2$, so there is no other alternative:
 $H$ must be either one or the other.
 For the second statement, see \cite[p.~192]{DFN1}
 and~\cite[p.~103]{GSS2}.) \linebreak
 Concerning maximality, note first that $SO(2)$ is a maximal connected
 Lie subgroup of~$SO(3)$, since its Lie algebra $\mathfrak{so}(2)$ is
 a maximal subalgebra of $\mathfrak{so}(3)$, but fails to be a maximal
 closed subgroup of~$SO(3)$, since $\, SO(2) \subset O(2)$.
 However, $O(2)$ is a maximal closed subgroup of~$SO(3)$.
 Geometrically, $SO(2)$ is a circle and $\, O(2)~\cong~\mathbb{Z}_2
 \ltimes SO(2)$ is the disjoint union of two circles, generated by
 $SO(2)$ and any orthogonal matrix of determinant $-1$, such as
 \[
  \left( \begin{array}{ccc} 
          0 & 1  \\
          1 & 0 \\
         \end{array} \right)~,
 \]    
 which represents reflection along the main diagonal.
 Similarly, looking at the list of finite subgroups above, we see
 that the first three fail to be maximal closed subgroups of~$SO(3)$,
 since $\, \mathbb{Z}_n \subset D_n \subset O(2) \,$ and
 $\, \mathrm{T} \subset \mathrm{O}$, but the last two are.%
 \footnote{Although it is true that $S_4$ can be embedded into~$A_5$
 as a subgroup, $\mathrm{O}$ cannot be embedded into~$\mathrm{I}$
 within $SO(3)$, even up to conjugacy: the symmetry group of the
 octahedron is not a subgroup of the symmetry group of the
 icosahedron.}
\end{example}

Of course, we may wonder what would happen if we were to replace the
concept of a maximal closed subgroup by that of a maximal Lie subgroup.
To get an idea, let us look again at the subgroups of the rotation group.
\begin{example}
 Consider again the group $\, G = SO(3) \,$ of all rotations in~$\mathbb{R}^3$.
 Then it follows from the arguments given in Example~\ref{ex:SO3} above
 that the orthogonal subgroup $O(2)$ is a maximal Lie subgroup of~$SO(3)$.
 On the other hand, none of the finite subgroups listed in Example~%
 \ref{ex:SO3} above is a maximal Lie subgroup of~$SO(3)$.
 For example, the octahedral group $\mathrm{O}$ is contained in a bigger
 subgroup $SO(3,\mathbb{Q})$, which is clearly a discrete Lie subgroup.
 The same is true for the icosahedral group $\mathrm{I}$, with
 $\mathbb{Q}$ replaced by $\mathbb{Q}(\sqrt 5)$.
 Moroever, none of the discrete Lie subgroups $SO(3,\mathbb{Q}(\alpha_1,
 \ldots,\alpha_k))$ of~$SO(3)$, where $\{\alpha_1,\ldots,\alpha_k\}$ is
 a finite set of irrational numbers linearly independent over $\mathbb{Q}$,
 is maximal, since one can always extend any one of them to a bigger one
 by passing to a bigger field extension of~$\mathbb{Q}$, with more generators.
\end{example}
\noindent
In general, it is clear that maximal Lie subgroups fall into two distinct
classes: closed maximal Lie subgroups and dense maximal Lie subgroups.
(Indeed, if $M$ is a maximal Lie subgroup of~$G$, then its closure
$\bar{M}$ is a Lie subgroup of~$G$ containing~$M$ and hence
$\, \bar{M} = M \,$ or $\, \bar{M} = G$.)
Now closed maximal Lie subgroups are maximal closed subgroups,
but the two notions do \emph{not} coincide since the former are
maximal among all Lie subgroups while the latter are only maximal
among all closed subgroups.
(Taking into account that the expression ``closed Lie subgroup''
is a tautology and is always abbreviated to ``closed subgroup'',
we may express this fact by the statement that the adjectives
``closed'' and ``maximal'' do not commute.)
This discrepancy is not specific to the above example but rather,
at least for discrete subgroups, a general and systematic feature,
since the same argument as in the example shows that if $G$ is
any compact connected matrix Lie group ($G \subset SO(n)$), then
a finite subgroup of~$G$ can be a maximal closed subgroup of~$G$
but cannot be a closed maximal Lie subgroup of~$G$, since it is
contained in a dense discrete Lie subgroup of~$G$ of the form
$\, G \cap SO(n,\mathbb{Q}(\alpha_1,...,\alpha_k))$, where the
$\alpha$'s are chosen such that all matrix elements of all
elements of~$H$ are rational linear combinations of the
$\alpha$'s.
Moreover, none of these dense discrete Lie subgroups is a maximal
Lie subgroup.
In fact, we know of no example of a dense maximal Lie subgroup,
and the arguments given here support the conjecture that such
subgroups do not exist.

These considerations show that studying maximal Lie subgroups is largely
uninteresting: the closed ones can be found among the maximal closed
subgroups, and on the dense ones, practically nothing is known, not
even regarding the question whether they exist at all.
And even if they do, they are uninteresting at least from the point
of view of representation theory, since by continuity, the restriction
of an irreducible representation of a Lie group to a dense Lie subgroup
remains irreducible: there is no branching.
Thus in order to simplify the terminology, we adhere to the general
convention already mentioned in the footnote on the title page:
\vspace{2ex}
\begin{center}
 \emph{Convention: ``maximal subgroup'' always means
       ``maximal closed subgroup''.}
\end{center}
\vspace{2ex}

With this convention in mind, let us collect a few initial observations
on the problem of determining the maximal subgroups of Lie groups, and in
particular, of compact Lie groups.
As in the case of maximal connected Lie subgroups, the Lie correspondence
(i.e., the one-to one correspondence between connected Lie subgroups of a
connected Lie group and subalgebras of its Lie algebra~\cite[p.~47]{Kn})
plays an important role, but in a different way, because it is clear
that if we allow for Lie subgroups (more specifically, closed subgroups)
which are not connected, several different ones may have the same
connected one-component and hence correspond to the same Lie
subalgebra, so we must face the question which of them, if~any,
is a maximal one.
That this is not completely straightforward can already be seen from
Example~\ref{ex:SO3} above: it shows that, even when the ambient Lie
group~$G$ is connected, the classification of its maximal subgroups
is by no means equivalent to the classification of the maximal sub%
algebras of its Lie algebra~$\mathfrak{g}$, simply because the
connected one-component $M_0$ of a maximal subgroup $M$ of~$G$
and the corresponding subalgebra $\mathfrak{m}$ are not
necessarily maximal.
Quite to the contrary, $M_0$ and $\mathfrak{m}$ may even be trivial,
which means that $M$ must be discrete, but that does of course not
mean that $M$ must be trivial.
On the other hand, $M_0$ may be maximal among the connected Lie subgroups
of~$G$ and yet fail to be maximal among all the Lie subgroups of~$G$ or
even the closed subgroups of~$G$.
Additional complications arise when the ambient Lie group~$G$ is~not
connected, but these can be handled by introducing a modified concept
of maximality, giving rise to the notions of a ``maximal invariant''
subalgebra or subgroup and of a ``(quasi)primitive'' subalgebra;
these will be specified in a mathematically precise manner in
Definitions~\ref{def:maxisa}, \ref{def:maxisg} and~\ref{def:primsa}
below, with the conventions stated in Remarks~\ref{rmk:gaminv}
and~\ref{rmk:siginv}.

Before embarking on this program, let us fix some notation that will
be used constantly throughout this paper, often without further mention.
In all that follows, $G$ will always denote a Lie group and $\mathfrak{g}$
will denote its Lie algebra.%
\footnote{All Lie groups and Lie algebras considered in this paper
are assumed to be finite-dimensional, without further mention.}
The~case of main interest to us is when $G$ is compact, but since many
of the results in the first part of the paper hold more generally, we
have refrained from assuming this right from the start.
Similarly, a case of particular interest and importance is when $G$
is connected, but the fact that maximal subgroups of connected Lie
groups need not be connected, together with the desire to repeat
the procedure of passing to a maximal subgroup in order to construct
chains of subgroups where each subgroup is maximal in the previous
one, leads us to refrain from assuming this right from the start.
Thus in general, $G_0$ will denote the connected one-component
of~$G$, $\Gamma$ will denote the component group of~$G$, that
is, the discrete quotient group~$G/G_0$, and $\pi$ will denote
the canonical projection from~$G$ to~$\Gamma$, so we have the
following short exact sequence:
\begin{equation} \label{eq:GREXT1}
 \{1\}~\longrightarrow~G_0~\longrightarrow~G~\stackrel{\pi}{\longrightarrow}~
                       \Gamma~\longrightarrow~\{1\}~.
\end{equation}
Hence $G$ is an upwards extension of~$G_0$ by~$\Gamma$ (or downwards extension
of~$\Gamma$ by~$G_0$):
\begin{equation} \label{eq:GREXT2}
 G~=~G_0 \,.\, \Gamma
\end{equation}
Further, we shall write $\mathrm{Aut}(G_0)$ for the group of automorphisms
of~$G_0$, $\mathrm{Inn}(G_0)$ for the normal subgroup of inner automorphisms
of~$G_0$ (i.e., automorphisms of~$G_0$ given by conjugation with elements
of~$G_0$) and $\mathrm{Out}(G_0)$ for the quotient group
\[
 \mathrm{Out}(G_0) = \mathrm{Aut}(G_0)/\mathrm{Inn}(G_0)
\]
which, by abuse of language, is called the group of outer automorphisms
of~$G_0$ even though its elements are not automorphisms but rather
equivalence classes of auto\-morphisms.
Similarly, we shall write $\mathrm{Aut}(\mathfrak{g})$ for the group of
automorphisms of~$\mathfrak{g}$, $\mathrm{Inn}(\mathfrak{g})$ for the normal
subgroup of inner automorphisms of~$\mathfrak{g}$ (i.e., automorphisms
of~$\mathfrak{g}$ of the form $\mathrm{Ad}(g_0)$ with $\, g_0 \in G_0$)
and $\mathrm{Out}(\mathfrak{g})$ for the quotient group
\[
 \mathrm{Out}(\mathfrak{g})~
 =~\mathrm{Aut}(\mathfrak{g})/\mathrm{Inn}(\mathfrak{g})
\]
which, by abuse of language, is called the group of outer automorphisms
of~$\mathfrak{g}$ even though its elements are not automorphisms but
rather equivalence classes of auto\-morphisms.
With this notation, we see that any short exact sequence of the
form~(\ref{eq:GREXT1}) gives rise to homomorphisms
\begin{equation} \label{eq:GREXT3}
 \Gamma~\longrightarrow~\mathrm{Out}(G_0)
\end{equation}
and
\begin{equation} \label{eq:GREXT4}
 \Gamma~\longrightarrow~\mathrm{Out}(\mathfrak{g})
\end{equation}
defined by taking the left coset $g G_0$ of an element $g$ of~$G$
to the equivalence class of the auto\-morphism of~$G_0$ given by
conjugation with~$g$ and of the auto\-morphism of~$\mathfrak{g}$
given by~$\mathrm{Ad}(g)$, respectively; its image is a subgroup
of the respective outer automorphism group $\mathrm{Out}(G_0)$ or
$\mathrm{Out}(\mathfrak{g})$ which we shall typically (and in
either case) denote by $\Sigma$.
Of course, a natural question to ask is whether these homomorphisms
can be lifted to homomorphisms
\begin{equation} \label{eq:GREXT5}
 \Gamma~\longrightarrow~\mathrm{Aut}(G_0)
\end{equation}
and
\begin{equation} \label{eq:GREXT6}
 \Gamma~\longrightarrow~\mathrm{Aut}(\mathfrak{g})
\end{equation}
respectively.
A particularly important case where this happens is when the extension~%
(\ref{eq:GREXT2}) is \emph{split}, that is, when there exists a homomorphism
of $\Gamma$ into~$G$ which is a right inverse to $\pi$: then identifying
$\Gamma$ with its image by this homomorphism allows to consider $\Gamma$
as a subgroup of~$G$ and~$G$ as the semidirect product of~$G_0$ and~$\Gamma$:
\begin{equation}
 G~=~G_0 : \Gamma
\end{equation}
Unfortunately, this is not always the case, even for compact $G$~\cite{Ha}.
However, there is a slightly weaker statement that does hold
for any compact Lie group, known as Lee's supplement theorem,%
\footnote{We are grateful to one of the referees for drawing our
attention to this theorem.}
according to which $G$ always contains a finite subgroup $\Lambda$,
called a \emph{Lee supplement} for~$G_0$ in~$G$, which is such that
(a)~$\Lambda$ and~$G_0$ generate $G$ and (b)~their intersection
$\, \Lambda_0 = \Lambda \cap G_0 \,$ is a normal subgroup of $G$
contained in the center of~$G_0$ \cite[p.~272]{HM}: then obviously,
there is a canonical isomorphism $\, \Gamma \cong \Lambda/\Lambda_0 \,$
obtained by restricting the projection $\, \pi: G \rightarrow \Gamma \,$
to $\Lambda$ and observing that this restriction is still surjective,
due to~(a), and has $\Lambda_0$ as its kernel.
(The split case corresponds to $\, \Lambda = \Gamma$, $\Lambda_0 = \{1\}$.)
As a result, we obtain, as in the split case, actions of $\, \Gamma \cong
\Lambda/\Lambda_0$ on~$G_0$ and on~$\mathfrak{g}$ by automorphisms,
defined by restricting the respective actions of~$G$ on~$G_0$ by conjugation
and on~$\mathfrak{g}$ by $\mathrm{Ad}$ to $\Lambda$ and observing that
$\Lambda_0$ acts trivially because it is contained in the center of~$G_0$.
Therefore, in the compact case, the homomorphisms~(\ref{eq:GREXT3}) and~%
(\ref{eq:GREXT4}) can always be lifted to homomorphisms of the form~%
(\ref{eq:GREXT5}) and~(\ref{eq:GREXT6}), respectively, and hence $\Sigma$
can be realized as a subgroup of~$\mathrm{Aut}(G_0)$ or~$\mathrm{Aut}%
(\mathfrak{g})$, rather than just of~$\mathrm{Out}(G_0)$ or~$\mathrm{Out}%
(\mathfrak{g})$, respectively.


One of the most important tools in the study of maximal subgroups
is the concept of the normalizer of a subgroup.
Suppose that $H$ is a Lie subgroup of~$G$ with connected one-component
$H_0$ and $\mathfrak{h}$ is the Lie algebra of both $H$ and~$H_0$.
We denote by
\begin{equation} \label{eq:NORM1}
 N_G(H_0)~=~\{ \, g \in G~|~g H_0 g^{-1} \subset H_0 \, \}
\end{equation}
the normalizer of $H_0$ in $G$ and similarly by
\begin{equation} \label{eq:NORM2}
 N_G(\mathfrak{h})~
 =~\{ \, g \in G~|~\mathrm{Ad}(g) \mathfrak{h} \subset \mathfrak{h} \, \}
\end{equation}
the normalizer of $\mathfrak{h}$ in $G$ with respect to the
adjoint representation.
Obviously, we have
\begin{equation} \label{eq:NORM3}
 N_G(H_0)~=~N_G(\mathfrak{h})
\end{equation}
and we shall denote this group simply by $N$ whenever there is no
danger of confusion.
It is easy to see that $N$ is always a closed subgroup of~$G$\,%
since using a decomposition of the Lie algebra $\mathfrak{g}$
into the direct sum of the subalgebra $\mathfrak{h}$ and an
arbitrarily chosen complementary subspace, $N$ can be written
as the inverse image of the ``block triangular subgroup''
$\, \{ \, T \in \mathrm{GL}(\mathfrak{g})~|~T(\mathfrak{h})
\subset \mathfrak{h} \, \} \,$ under the homomorphism
$\; \mathrm{Ad}: G \rightarrow \mathrm{GL}(\mathfrak{g})$.
The~Lie algebra of $N$, denoted by $\mathfrak{n}$, is
\begin{equation} \label{eq:NORM4}
 \mathfrak{n}~=~\{ \, X \in \mathfrak{g}~|~
                   \mathrm{ad}(X) \mathfrak{h} \subset \mathfrak{h} \, \}
\end{equation}
and is the normalizer of $\mathfrak{h}$ in $\mathfrak{g}$.
Of course, we may also consider the normalizer $N_G(H)$ of~$H$ in~$G$
which is defined analogously but will be needed only occasionally,%
\footnote{In passing, we note that the normalizer of a general Lie subgroup
$H$ may fail to be closed and hence it cannot even be guaranteed ``a priori''
that it is a Lie subgroup, but it is an easy exercise to check that this
cannot happen when $H$ is closed (or, as argued before, when $H$ is
connected).}
and we then have the following sequence of inclusions which, although
very simple to prove, is fundamental:
\begin{equation} \label{eq:INCL1}
 H_0 \subset H \subset N_G(H) \subset N_G(H_0)~.
\end{equation}
Indeed, the last inclusion follows by observing that for any $\, n \in N_G(H)$,
we have $\, n H n^{-1} \subset H \,$ and hence $\, n H_0 n^{-1}
\subset H_0 \,$ since conjugation by $n$ is a homeomorphism of $H$ and
thus preserves its connected one-component.
The corresponding inclusion at the level of Lie algebras is
\begin{equation} \label{eq:INCL2}
 \mathfrak{h} \subset \mathfrak{n}~,
\end{equation}
which is also obvious from equation~(\ref{eq:NORM4}).
\begin{remark}
 In general, $N$ may be larger than $H_0$ or $H$ and
 $\mathfrak{n}$ may be larger than $\mathfrak{h}$.
 An extreme case is when $N$ contains $G_0$, i.e., $H_0$ is a normal sub%
 group of~$G_0$, or equivalently, when $\, \mathfrak{n} = \mathfrak{g}$,
 i.e., $\mathfrak{h}$ is an ideal of $\mathfrak{g}$.
 The other extreme case occurs when $N$ has $H_0$ as its connected
 one-component, or equivalently, when~$N/H_0$ is a discrete group,
 which in turn is equivalent to the condition that $\, \mathfrak{n}
 = \mathfrak{h}$.
\end{remark}
\begin{Definition} \label{def:salgan}
 Let $G$ be a Lie group and $\mathfrak{g}$ be a Lie algebra.
 A closed subgroup $H$ of~$G$ is called \emph{self-normalizing}
 if it coincides with its own normalizer $N_G(H)$ in~$G$.
 Similarly, a subalgebra $\mathfrak{h}$ of $\mathfrak{g}$
 is called \emph{self-normalizing} if it coincides with its
 own normalizer $\mathfrak{n}$ in~$\mathfrak{g}$.
\end{Definition}
In order to complete the picture, consider the general case
where $\mathfrak{h}$ is any subalgebra of~$\mathfrak{g}$.
Denoting its normalizer in $\mathfrak{g}$ by $\mathfrak{n}_1$,
the normalizer of $\mathfrak{n}_1$ in $\mathfrak{g}$ by
$\mathfrak{n}_2$ and so on, we obtain an ascending sequence
of subalgebras of $\mathfrak{g}$ which, for dimensional reasons,
ends at some proper subalgebra $\mathfrak{n}_r$ of $\mathfrak{g}$:
\begin{equation}
 \mathfrak{h} \subset \mathfrak{n}_1 \subset \mathfrak{n}_2
              \subset \ldots \subset \mathfrak{n}_r \subset \mathfrak{g}~.
\end{equation}
Analogously, using equation~(\ref{eq:INCL1}), we obtain an ascending
sequence of subgroups of $G$,
\begin{equation}
 H \subset N_1 \subset N_2 \subset \ldots \subset N_r \subset G~,
\end{equation}
where $\, N_1 = N_G(\mathfrak{h}),\,N_2 = N_G(\mathfrak{n}_1),\ldots,
N_r = N_G(\mathfrak{n}_{r-1})$. 
There are then two possibilities: either $\mathfrak{n}_r$ is an ideal
of~$\mathfrak{g}$ or else~$\mathfrak{n}_r$ is self-normalizing.

\section{The three types of maximal subgroups}

The first main issue in this paper is to show that the problem of
classifying the maximal subgroups of an arbitrary Lie group $G$ can
be reduced to (a)~that of classifying the maximal subgroups of its
component group~$\Gamma$ together with (b)~that of classifying the
maximal $\Gamma$-invariant subgroups of its connected one-component~%
$G_0$, and these split naturally into two types: those whose connected
one-component is a normal subgroup and those who are the normalizer of
their own connected one-component.
The first step in this reduction, which in the context of Lie groups
has apparently never been formulated explicitly (although its analogue
in the context of finite groups seems to be well known), will be quite
elementary once we figure out the precise meaning of the term ``maximal
$\Gamma$-invariant subgroup'' of~$G_0$, the main obstacle being that
there is no canonical way to make $\Gamma$ act on~$G_0$.
(Here, the emphasis is on the word ``canonical'' because we have
already mentioned that, at least for compact Lie groups, the
homomorphism~(\ref{eq:GREXT3}) can always be lifted to a homomorphism of
the form~(\ref{eq:GREXT5}), so the main problem with understanding the
meaning of ``$\Gamma$-invariance'' is not lack of existence but rather
lack of uniqueness: the lifting is not canonical since it depends on the
choice of the splitting, in the split case, or of the Lee supplement,
in the general compact case.)
However, and this is the key observation, $\Gamma$ does act naturally
on conjugacy classes of elements and also of subgroups of~$G_0$, and
this will suffice since the property of maximality of subgroups is
invariant under conjugation.

To explain this in more detail, let us introduce some more notation.
First, given any element $\gamma$ of the component group~$\Gamma$, we shall
denote the corresponding connected component of~$G$ by~$G_\gamma$, i.e.,
$G_\gamma = \pi^{-1}(\{\gamma\})$, so that in particular, $G_1 \equiv G_0$.
Similarly, given any subset $H$ of~$G$, we write $\, H_\gamma = H \cap
G_\gamma \,$ and, in particular, $H_1 = H \cap G_1$.
Note that in the case of interest to us here, namely when $H$ is a
closed subgroup (or Lie subgroup) of~$G$ with connected one-component
$H_0$ and Lie algebra $\mathfrak{h}$, there is no reason for $H_1$ to
be connected, so in general, $H_1$ will be an intermediate closed subgroup 
(or Lie subgroup) of~$G_0$ that lies in between $H_0$ and~$H$:
just like $H_0$, it is necessarily normal in~$H$, and of course, $H$,
$H_1$ and $H_0$ all have the same Lie algebra $\mathfrak{h}$.
Next, we shall write $[g_0^{}]$ for the conjugacy class of an
element $g_0^{}$ of~$G_0^{}$ in~$G_0^{}$ and, similarly, $[H_0^{}]$
for the conjugacy class of a subset (subgroup\,/\,closed subgroup\,/\,%
Lie subgroup) $H_0^{}$ of~$G_0^{}$ in~$G_0^{}$; explicitly, as sets,
\[
 [g_0^{}]~=~\{ \, \tilde{g}_0^{} g_0^{} \tilde{g}_0^{-1}~|~
               \tilde{g}_0 \in G_0 \, \}~~~,~~~
 [H_0^{}]~=~\{ \, \tilde{g}_0^{} g_0^{} \tilde{g}_0^{-1}~|~
               \tilde{g}_0 \in G_0 \,,\, g_0^{} \in H_0^{} \, \}~.
\]
With this notation, given $\, g \in G_\gamma$, we define
\begin{equation} \label{eq:ACTGAMG}
 \gamma \cdot [g_0^{}]~=~[ g g_0^{} g^{-1} ]~~~,~~~
 \gamma \cdot [H_0^{}]~=~[ g H_0^{} g^{-1} ]~.
\end{equation}
Observe that this provides well defined actions of $\Gamma$ on the set
of conjugacy classes of elements of~$G_0^{}$ and on the set of conjugacy
classes of subsets (subgroups\,/\,closed subgroups\,/\,Lie subgroups)
of~$G_0^{}$, since if we replace $g$ by any other representative of~$\gamma$,
say $\, g \hat{g}_0^{}$, and $g_0^{}$ by any other representative of the
same conjugacy class, say $\, \tilde{g}_0^{} g_0^{} \tilde{g}_0^{-1}$,
the resulting element $\; (g \hat{g}_0^{})(\tilde{g}_0^{} g_0^{}
\tilde{g}_0^{-1})(g \hat{g}_0^{})^{-1} \;$ equals
$\; \bigl( g (\hat{g}_0^{} \tilde{g}_0^{}) g^{-1} \bigr) g g_0^{} g^{-1}
\bigl( g (\hat{g}_0^{} \tilde{g}_0^{}) g^{-1} \bigr)^{-1}$
which is conjugate to $\, g g_0^{} g^{-1}$.
Similarly, we shall write $[X]^{\vphantom{A}}$ for the $\mathrm{Ad}(G_0)$-orbit
of an element $X$ of~$\mathfrak{g}$ and, similarly, $[\mathfrak{h}]$ for the
$\mathrm{Ad}(G_0)$-orbit of a subset (subspace\,/\,subalgebra) $\mathfrak{h}$
of~$\mathfrak{g}$; explicitly, as sets,
\[
 [X]~=~\{ \, \mathrm{Ad}(\tilde{g}_0^{}) X~|~\tilde{g}_0 \in G_0 \, \}~~~,~~~
 [\mathfrak{h}]~=~\{ \, \mathrm{Ad}(\tilde{g}_0^{}) X~|~
                     \tilde{g}_0 \in G_0 \,,\, X \in \mathfrak{h} \, \}~.
\]
With this notation, given $\, g \in G_\gamma$, we define
\begin{equation} \label{eq:ACTGAMA}
 \gamma \cdot [X]~=~[ \mathrm{Ad}(g) X ]~~~,~~~
 \gamma \cdot [\mathfrak{h}]~=~[ \mathrm{Ad}(g) \mathfrak{h} ]~.
\end{equation}
Observe that this provides well defined actions of $\Gamma$ on the
set of $\mathrm{Ad}(G_0)$-orbits of elements of~$G_0^{}$ and on the
set of $\mathrm{Ad}(G_0)$-orbits of subsets (subspaces\,/\,subalgebras)
of~$\mathfrak{g}$, since if we replace $g$ by any other representative
of~$\gamma$, say $\, g \hat{g}_0^{}$, and $X$ by any other
representative of the same $\mathrm{Ad}(G_0)$-orbit, say
$\, \mathrm{Ad}(\tilde{g}_0^{}) X$, the resulting element
$\; \mathrm{Ad}(g \hat{g}_0^{}) \mathrm{Ad}(\tilde{g}_0^{}) X \;$ equals
$\; \mathrm{Ad}(g (\hat{g}_0^{} \tilde{g}_0^{}) g^{-1}) \mathrm{Ad}(g) X \;$
which is in the same $\mathrm{Ad}(G_0)$-orbit as $\, \mathrm{Ad}(g) X$.
All these actions are canonical since they do not depend on any additional
data such as the choice or even the existence of a splitting or, more
generally, a Lee supplement.

Concerning the relation between subgroups of~$G$ and subgroups of~$G_0$,
we note that of course every closed subgroup (or Lie subgroup) $H$ of~$G$
gives rise to a closed subgroup (or Lie subgroup) $H_1$ of $G_0$ by taking
its intersection with $G_0$, but the converse is not true: there may
be closed subgroups (or Lie subgroups) of~$G_0$ that cannot be obtained
from closed subgroups (or Lie subgroups) of~$G$ in this way. The situation
is clarified by the following simple observation.
\begin{proposition} \label{prop:invsg}
 Let\/ $G$ be a Lie group with connected one-component\/~$G_0$ and
 component group\/~$\Gamma$, and let\/ $H_1$ be a closed subgroup
 (or Lie subgroup) of\/~$G_0$.
 Then\/ $H_1$ extends to a closed subgroup (or Lie subgroup)\/ $H$ of\/~$G$
 such that $\, H \cap G_0 = H_1 \,$ and $\, \pi(H) = \Gamma \,$ if and
 only if its conjugacy class\/ $[H_1]$ in\/~$G_0$ is\/ $\Gamma$-invariant.
 In this case, both\/ $H$ and\/~$H_1$ have the same connected one-%
 component\/~$H_0$ whose conjugacy class\/ $[H_0]$ in\/~$G_0$ is also\/
 $\Gamma$-invariant, and all have the same Lie algebra\/~$\mathfrak{h}$
 whose\/ $\mathrm{Ad}(G_0)$-orbit is also\/ $\Gamma$-invariant.
\end{proposition}
\begin{proof}
 Assume first that $H$ is a closed subgroup (or Lie subgroup) of~$G$
 such that $\, H \cap G_0 = H_1 \,$ and $\, \pi(H) = \Gamma$.
 Then for any $\, \gamma \in \Gamma$, $H_\gamma \neq \emptyset$,
 and picking any $\, h_\gamma \in H_\gamma$, we get
 \[
  h_\gamma^{-1} H_\gamma^{}~=~H_1^{}~~~,~~~H_\gamma^{} h_\gamma^{-1}~=~H_1^{}
 \]
 since both left and right translation by $h_\gamma$ are
 bijections that map $G_0$ onto~$G_\gamma$ and $H$ onto itself.
 But this means that conjugation by $h_\gamma$ leaves $H_1$ invariant and
 therefore, being a homeomorphism, also leaves its connected one-component
 $H_0$ invariant, while $\mathrm{Ad}(h_\gamma)$ leaves $\mathfrak{h}$ inva%
 riant, so $[H_1]$, $[H_0]$ and $[\mathfrak{h}]$ are all $\Gamma$-invariant.
 (Of course, $h_\gamma$ is not unique but may be freely multiplied
 by elements of~$H_1$, from either side, and that is precisely the
 possible ambiguity in the choice of~$h_\gamma$.)
 Conversely, the condition that $\gamma$ leaves $[H_1]$ invariant means
 that, for any representative $g_\gamma$ of $\gamma$ in~$G_\gamma$, the
 subgroups $g_\gamma^{} H_1^{\phantom{\beta}} g_\gamma^{-1}$ and $H_1^{}$ are
 conjugate in~$G_0$, i.e., there exists $\, g_0^{(\gamma)} \in G_0^{} \,$
 such that $\; g_\gamma^{} H_1^{\phantom{\beta}} g_\gamma^{-1} = g_0^{(\gamma)}
 H_1^{\phantom{\beta}} g_0^{(\gamma)\,{-1}}$, which implies that conjugation
 by $\, h_\gamma^{} = g_\gamma^{} g_0^{(\gamma)\,{-1}} \,$ leaves $H_1^{}$
 invariant.
 Now since $\, h_\gamma \in G_\gamma$, we can define $H$ by
 \[
  H~=~\dot{\bigcup\limits_{\gamma \in \Gamma}} \, H_\gamma
  \qquad \mbox{where} \qquad
  H_\gamma~=~h_\gamma H_1~=~H_1 h_\gamma
 \]
 and check immediately that $H$ is a closed subgroup (or Lie subgroup)
 of~$G$ such that $\, H \cap G_0 = H_1$ (provided we choose $\, h_1 = 1$,
 which is always possible) and $\, \pi(H) = \Gamma$.
 (Of course, once again, $h_\gamma$ is not unique but may be freely
 multiplied by elements of~$H_1$, from either side: this does not
 change the definition of~$H$.
 However, it may even be multiplied by elements of the normalizer
 of~$H_1$ in~$G_0$, and this can change the definition of~$H$.
 In such a case, there are several subgroups of~$G$ having the
 same intersection with~$G_0$, but this is no problem since
 nowhere do we need, nor did we claim, that $H$ is uniquely
 determined by~$H_1$: this cannot be guaranteed in general.)
\end{proof}
\begin{remark} \label{rmk:gaminv}
 By abuse of language, we shall say that a closed subgroup (or Lie subgroup)
 $H_1$ of~$G_0$ is $\Gamma$-invariant if its conjugacy class $[H_1]$ in~$G_0$
 is $\Gamma$-invariant.%
 \footnote{If one wants to be more precise, one can say that $H_1$ is
 $\Gamma$-invariant up to conjugacy.}  
 Similarly, we shall say that a subalgebra $\mathfrak{h}$ of~$\mathfrak{g}$
 is $\Gamma$-invariant if its $\mathrm{Ad}(G_0)$-orbit $[\mathfrak{h}]$
 is $\Gamma$-invariant.
 Explicitly, this property can be formulated as the condition that every
 $\, \gamma \in \Gamma \,$ has at least one representative $\, g_\gamma \in
 G_\gamma \,$ such that conjugation by $g_\gamma$ leaves $H_1$ invariant and
 $\mathrm{Ad}(g_\gamma)$ leaves $\mathfrak{h}$ invariant, or in other words,
 as the condition that the respective normalizers $N_G(H_1)$ of~$H_1$ and
 $N_G(\mathfrak{h})$ of~$\mathfrak{h}$ in~$G$ meet every connected
 component of~$G$.
\end{remark}

With this language, we can complement Definitions~\ref{def:maxsa}
and~\ref{def:maxsg} as follows.
\begin{Definition} \label{def:maxisa}
 Let $G$ be a Lie group with Lie algebra~$\mathfrak{g}$
 and component group~$\Gamma$.
 A \emph{maximal $\Gamma$-invariant subalgebra} of~$\mathfrak{g}$ is a
 proper $\Gamma$-invariant sub\-algebra $\mathfrak{m}$ of~$\mathfrak{g}$
 such that if $\tilde{\mathfrak{m}}$ is any $\Gamma$-invariant subalgebra
 of~$\mathfrak{g}$ with $\, \mathfrak{m} \subset \tilde{\mathfrak{m}}
 \subset \mathfrak{g}$, then $\, \tilde{\mathfrak{m}} = \mathfrak{m} \,$
 or $\, \tilde{\mathfrak{m}} = \mathfrak{g}$.
 The~same terminology applies when the expression ``subalgebra''
 is everywhere replaced by the expression ``ideal''.
\end{Definition}
\begin{Definition} \label{def:maxisg}
 Let $G$ be a Lie group with connected one-component~$G_0$
 and component group~$\Gamma$.
 A \emph{maximal\/ $\Gamma$-invariant closed subgroup} of~$G_0$ is a
 proper $\Gamma$-invariant closed subgroup~$M_1$ of~$G_0$ such that if
 $\tilde{M}_1$ is any $\Gamma$-invariant closed subgroup of~$G_0$ with
 $\, M_1 \subset \tilde{M}_1 \subset G_0$, then $\, \tilde{M}_1 = M_1 \,$
 or $\, \tilde{M}_1 = G_0$.
 The same terminology applies when the expression ``closed subgroup''
 is everywhere replaced by the expression ``closed normal subgroup''.
\end{Definition}

In order for this concept to become useful, we should convince ourselves
that non-trivial proper $\Gamma$-invariant closed subgroups exist at all.
As it turns out, there is only one special situation where this may fail
to be true; we give it a name:
\begin{Definition} \label{def:irrat}
 Let $G$ be a Lie group with connected one-component~$G_0$
 and component group~$\Gamma$.
 We say that $G_0$ is \emph{$\Gamma$-irrational} if $G_0$ is abelian
 and simply connected, i.e., $G_0 \cong \mathbb{R}^n$, and $\Gamma$ acts
 on it irreducibly and so that it contains no $\Gamma$-invariant lattice.
\end{Definition}
\begin{proposition} \label{prop:invsg1}
 Let\/ $G$ be a Lie group with connected one-component\/~$G_0$ and
 component group\/~$\Gamma$, and suppose\/ $G_0$ is not\/ $\Gamma$-%
 irrational.
 Then $G_0$ contains non-trivial proper $\Gamma$-invariant closed
 subgroups.
\end{proposition}
\begin{proof}
 The only information available on the action of~$\Gamma$ on (conjugacy
 classes of elements and subgroups of) $G_0$ is that it preserves
 products; therefore, we shall search for non-trivial proper closed
 subgroups $H$ of~$G_0$ which are invariant under all automorphisms
 of~$G_0$, up to conjugacy in~$G_0$.
 \pagebreak
 \begin{itemize}
  \item If $G_0$ is neither semisimple nor solvable, take $H$ to be
        its radical, which is then a non-trivial proper connected
        closed subgroup and invariant under all automorphisms
        \cite[Proposition~10.12, p.~207]{SW}.
  \item If $G_0$ is solvable but not abelian, take $H$ to be one of the
        non-trivial proper subgroups $(G_0)_k$ of its closed commutator
        subgroup series, all of which are closed subgroups and invariant
        under all automorphisms.
        (Explicitly, let $G_0^{(1)}$ be its algebraic commutator subgroup
        $G_0'$, whose elements are products of commutators $\, g_1^{} \,
        g_2^{} \, g_1^{-1} g_2^{-1} \,$ with $\, g_1^{},g_2^{} \in G_0^{}$,
        define $G_0^{(k)}$ recursively to be the algebraic commutator subgroup
        $G_0^{(k-1)\,'}$ of~$G_0^{(k-1)}$, and let $(G_0)_k$ be the closure
        of~$G_0^{(k)}$ in~$G_0$.
        Then it is easy to see by induction on $k$ that each $G_0^{(k)}$ is
        an abstract normal subgroup of~$G_0$ and each $(G_0)_k$ is a closed
        normal subgroup of~$G_0$ such that both $G_0^{(k)}$ and $(G_0)_k$
        are invariant under all automorphisms of~$G_0$ and both
        $G_0^{(k-1)}/G_0^{(k)}$ and $(G_0)_{k-1}/(G_0)_k$ are abelian;
        moreover, it is also known that $G_0$ is solvable iff the
        descending abstract commutator subgroup series $\; G_0 \supset
        G_0^{(1)} \supset \ldots \supset G_0^{(k)} \supset \ldots$ terminates
        after a finite number of steps and also iff the descending closed
        commutator subgroup series $\; G_0 \supset (G_0)_1 \supset \ldots
        \supset (G_0)_k \supset \ldots \,$ terminates after a finite number
        of steps \cite[Chapter~10.1, pp.~201-204]{SW}.
        Then at least one of the $(G_0)_k$ must be non-trivial and proper.
        Note that the argument works even if $G_0^{(1)}$ should happen
        to be dense in~$G_0$, so that $(G_0)_1 = G_0$, but this just
        means that we have to look at the next step(s) in the series.)
  \item If $G_0$ is semisimple and compact but not abelian, take $H$
        to be any maximal torus, which is then a non-trivial proper
        connected closed subgroup and invariant under all automorphisms,
        up to conjugacy.
  \item If $G_0$ is semisimple and non-compact but not abelian, take $H$ to
        be any maximal compact subgroup, which is then a non-trivial proper
        connected closed subgroup and invariant under all automorphisms,
        up to conjugacy.
  \item If $G_0$ is abelian but not simply connected, so $\, G_0 \cong
        \mathbb{T}^p \times \mathbb{R}^q \,$ where $\mathbb{T}$ denotes
        the unit circle (one-dimensional torus) \cite[Theorem~6.20,
        p.~155]{SW}, with $\, p>0$, take $H$ to be $\, \mathbb{Z}_N^p
        \times \mathbb{R}_{\vphantom{N}}^q$, for any value of~$N$, which
        is then a non-trivial proper closed subgroup and invariant under
        all automorphisms since $\mathbb{Z}_N^p$ is the subgroup of~%
        $\mathbb{T}^p$ formed by the $N^{\mathrm{th}}$ roots of unity,
        and any automorphism of any abelian group preserves the set
        of $N^{\mathrm{th}}$ roots of unity.
  \item If $\, G_0 \cong \mathbb{R}^n \,$ is reducible under the action
        of~$\Gamma$, take $H$ to be the non-trivial proper $\Gamma$-%
        invariant subspace (and hence connected closed subgroup) of~%
        $\mathbb{R}^n$ whose existence is guaranteed by reducibility.
  \item If $\, G_0 \cong \mathbb{R}^n \,$ is irreducible under the
        action of~$\Gamma$ but contains a $\Gamma$-invariant lattice,
        take $H$ to be this lattice, which is a discrete closed sub%
        group of~$\mathbb{R}^n$ (and all discrete closed subgroups
        of~$\mathbb{R}^n$ are of this form \cite[Lemma~6.18, p.~155]{SW}).
        \\ \qedhere
 \end{itemize}
\end{proof}
\begin{remark} \label{rmk:irrat}
 From a practical point of view, the exceptional case of $G_0$ being
 $\Gamma$-irrational is completely irrelevant: it cannot occur if
 $G_0$ is compact, which is the case of main interest in this paper,
 or if $G_0$ is connected, or if $G_0$ is not abelian. And even when
 $\, G_0 \cong \mathbb{R}^n$, it cannot occur for standard group
 actions, such as those of reflection groups or crystallographic
 groups, which always leave some lattice invariant.
\end{remark}

As before, we could, in principle, consider the option of replacing
the concept of a maximal $\Gamma$-invariant closed subgroup by that
of a maximal $\Gamma$-invariant Lie subgroup, but again this will
not lead to anything interesting: closed maximal $\Gamma$-invariant
Lie subgroups are maximal $\Gamma$-invariant closed subgroups, and
on dense maximal $\Gamma$-invariant Lie subgroups, nothing is known,
not even regarding their existence. Thus once again, we shall adhere
to the general convention already mentioned in the footnote on the
title page:
\vspace{2ex}
\begin{center}
 \emph{Convention: ``maximal $\Gamma$-invariant subgroup'' always means \\
       ``maximal $\Gamma$-invariant closed subgroup''.}
\end{center}
\vspace{2ex}

Now we are ready to formulate our first main theorem, which provides
a splitting of the maximal subgroups of an arbitrary Lie group~$G$
into three distinct classes:
\begin{theorem} \label{theo:maxsg1}
 Let\/ $G$ be a Lie group with connected one-component\/~$G_0$ and
 component group\/~$\Gamma$, and let\/ $M$ be a maximal subgroup
 of\/~$G$ with connected one-component\/~$M_0$.
 Set $\, M_1 = M \,\cap\, G_0$, so $\, M_0 \subset M_1 \subset M$.
 Then one of the following three alternatives holds.
 \begin{enumerate}
  \item $M_0 = G_0$, $M_1 = G_0 \,$ and\/ $M$ is an (upwards) extension
        of\/~$G_0$ by a maximal subgroup\/~$\Xi$ of\/~$\Gamma$:
        \begin{equation} 
         M~=~G_0 \,.\, \Xi~.
        \end{equation}
        We then say that\/ $M$ is of\, \textbf{trivial type}.
  \item $M_0$ is a proper\/ $\Gamma$-invariant closed connected normal 
        subgroup of\/~$G_0$, or equi\-valently, a proper connected closed
        normal subgroup of\/~$G$, $M_1$ is an (upwards) extension of\/~$M_0$
        by a discrete maximal $\Gamma$-invariant subgroup\/ $D_0$ of the
        quotient group\/~$G_0/M_0$ and\/ $M$ is an (upwards) extension
        of\/~$M_1$ by\/~$\Gamma$:
        \begin{equation} 
         M_1~=~M_0 \,.\, D_0~~~,~~~M~=~M_1 \,.\, \Gamma~.
        \end{equation}
        We then say that\/ $M$ and\/ $M_1$ are of\, \textbf{normal type}.
        Here, we may have $\, M_1 = M_0$, i.e., $D_0$ is trivial, but
        according to Proposition~\ref{prop:invsg1}, this can only occur
        when the quotient group\/ $G_0/M_0$ is $\Gamma$-irrational,
        and hence we shall refer to this exceptional situation%
        \footnote{See Remark~\ref{rmk:irrat}.}\ 
        as the\/ $\Gamma$-irrational type.

 \pagebreak

  \item $M_0$ is a proper\/ $\Gamma$-invariant closed connected non-normal
        subgroup of\/~$G_0$, $M_1$ is its normalizer in\/~$G_0$, $M$ is its
        normalizer in\/~$G$ and\/ $M$ is an (upwards) extension of\/~$M_1$
        by\/~$\Gamma$:
        \begin{equation} 
         M_1~=~N_{G_0}(M_0)~~~,~~~M~=~N_G(M_0)~~~,~~~M~=~M_1 \,.\, \Gamma~.
        \end{equation}
        We then say that\/ $M$ and\/~$M_1$ are of\, \textbf{normalizer type}.%
        \footnote{Here, the term ``normalizer'' may be thought of as an
        abbreviation for ``normalizer of its own connected one-component''
        or ``normalizer of its own Lie algebra''.}
 \end{enumerate}
 In the last two cases, $M$ meets every connected component of\/~$G$,
 and\/ $M_1$ is a maximal\/ $\Gamma$-invariant subgroup of\/~$G_0$
 which, except in the case of the\/ $\Gamma$-irrational type,
 is self-normalizing and contains the center of\/~$G_0$.
 Conversely, the condition that\/ $\Xi$ should be maximal in\/~$\Gamma$,
 in the first case, and that\/ $M_1$ should be maximal\/ $\Gamma$-invariant
 in\/~$G_0$, in the last two cases, guarantees that\/ $M$ will be maximal
 in\/~$G$.
\end{theorem}
\begin{remark}
 Note that Theorem~\ref{theo:maxsg1} does not provide a full ``if and
 only if'' statement because one aspect is still missing: the restrictions
 that must be imposed on~$M_0$ in order to guarantee that $M$ is really
 maximal, or equivalently, that $M_1$ is really maximal $\Gamma$-invariant.
 This question will be dealt with later; see~Theorem~\ref{theo:maxsg2}.
\end{remark}
\begin{proof}
 Assume $M$ is a maximal subgroup of~$G$, set $\, \Xi = \pi(M) \,$ and let
 $\tilde{M}$ be its preimage in~$G$ (i.e., $\tilde{M} = \pi^{-1}(\pi(M))$),
 $N$ be the normalizer of~$M_0$ in~$G$ and $N_1$ be the normalizer of~$M_0$
 in~$G_0$: then obviously, $N_1 = N \cap G_0 \,$ and, by Proposition~%
 \ref{prop:invsg}, $N_1$~is $\Gamma$-invariant.
 Since $\tilde{M}$ is a closed subgroup of~$G$ such that $\, M \subset
 \tilde{M} \subset G \,$, maximality of~$M$ in~$G$ implies that either
 $\, \tilde{M} = M \,$ or $\, \tilde{M} = G$.
 In~the first case ($\tilde{M} = M$), we conclude that $M$ is the
 (disjoint) union of all connected components of~$G$ that project
 to elements of~$\Xi$ (in particular, $G_0 \subset M$), and $M$
 being maximal in~$G$ is then equivalent to $\Xi$ being maximal
 in~$\Gamma$: this is the first alternative of the theorem.
 In~the second case ($\tilde{M} = G$), we conclude that $M$ meets
 every connected component of~$G$ (in~particular, $\Xi = \Gamma$),
 so $\, \Gamma = M/M_1$ and $M$ is an extension of~$M_1$ by $\Gamma$.
 It then follows from Proposition~\ref{prop:invsg} that $M_1$ is
 $\Gamma$-invariant and that $M$ being maximal in~$G$ is equivalent
 to $M_1$ being maximal $\Gamma$-invariant in~$G_0$.
 Note, however, that $M_0$ is just $\Gamma$-invariant and proper
 but may or may not be maximal $\Gamma$-invariant.
 Indeed, the fact that $M$ is maximal in~$G$ only implies that
 either $\, N = M$, $N_1 = M_1 \,$ or $\, N = G$, $N_1 = G_0$,
 corresponding to the third and to the second alternative of
 the theorem, respectively.
 Similarly, letting $N_2$ be the normalizer of~$M_1$ in~$G_0$ and
 noting that $N_2$ is closed and $\Gamma$-invariant since $M_1$ is,
 the fact that $M_1$ is maximal $\Gamma$-invariant in~$G_0$ implies
 that either $\, N_2 = M_1 \,$ or $\, N_2 = G_0$.
 In the first case, we conclude that $M_1$ is self-normalizing and
 contains the center of~$G_0$ (since the normalizer of any subgroup
 always contains the center).
 Therefore, we are left with the task of analyzing in more detail
 what happens when $M_0$ and/or $M_1$ are normal.
 To handle part of this analysis in one stroke, assume that $\bar{M}$
 is any $\Gamma$-invariant closed subgroup of~$G_0$ with connected
 one-component~$M_0$ (in particular, it can be either $M_0$ or~$M_1$),
 and suppose that $\bar{M}$ is normal.
 Then we may consider the quotient groups $G/\bar{M}$ and~%
 $G_0/\bar{M}$, which will fit into a short exact sequence 
 \[
  \{1\}~\longrightarrow~G_0/\bar{M}~\longrightarrow~G/\bar{M}~
        \stackrel{\bar{\pi}}{\longrightarrow}~\Gamma~\longrightarrow~\{1\}~,
 \]
 where $\bar{\pi}$ is induced from $\pi$ by passing to the quotient.
 Repeating the arguments at the beginning of this section, we see
 that the concept of $\Gamma$-invariant subgroups of~$G_0/\bar{M}$
 is well-defined, and it is then obvious that the non-trivial proper
 $\Gamma$-invariant closed subgroups of~$G_0/\bar{M}$ will, by taking
 the inverse image under the canonical projection from~$G_0$ to~%
 $G_0/\bar{M}$, correspond to the $\Gamma$-invariant closed subgroups
 of~$G_0$ properly containing $\bar{M}$ and properly contained in~$G_0$,
 with the additional property that the former are discrete if and only
 if the latter have connected one-component~$M_0$.
 In the case of the second alternative of the theorem, we can use
 this construction with $\, \bar{M} = M_0 \,$ to conclude that
 since $M_1$ is a maximal $\Gamma$-invariant subgroup of~$G_0$,
 $D_0 = M_1/M_0 \,$ must be a discrete maximal $\Gamma$-invariant
 subgroup of~$G_0/M_0$.
 For the final statements of the theorem about $M_1$, we can use
 the same construction with $\, \bar{M} = M_1$, combined with
 Proposition~\ref{prop:invsg1} applied to the quotient group~%
 $G_0/M_1$, to conclude that if $M_1$ is normal, then $G_0/M_1$
 must be $\Gamma$-irrational.
 In particular, $G_0/M_1 \cong \mathbb{R}^n$, which is simply
 connected, and an elementary application of the exact homotopy
 sequence to the bundle $\, G_0 \longrightarrow G_0/M_1 \,$ now
 shows that $M_1$ must be connected; hence we are in the case
 of the normal type, with $\, M_1 = M_0$.
\end{proof}

\begin{remark} \label{rmk:maxsg3}
 A substantial part of Theorem~\ref{theo:maxsg1}, namely that concerning the
 rela\-tion between $M$ and~$M_1$, has been motivated by a corresponding
 statement for finite groups relating maximal subgroups of finite groups
 to maximal subgroups of their extensions by groups of outer automorphisms
 (this is what we are dealing with here if the homomorphism~(\ref{eq:GREXT3})
 is supposed to have trivial kernel): a brief summary can be found in
 Section~5 of the Introduction to the \textsc{atlas}~\cite[p.~xix]{ATL}.
 (Note, however, that the notion of ``maximal $\Gamma$-invariant''
 subgroup does not appear there.)
 In particular, in Ref.~\cite{ATL}, a maximal subgroup $M$ of~$G$, viewed
 as an extension of a subgroup of~$G_0$, is called ``trivial'' if it is
 of the form $G_0 \,.\, \Xi$ with $\Xi$ a maximal subgroup of~$\Gamma$,
 as in the first alternative of Theorem~\ref{theo:maxsg1}, and we have
 decided to follow this terminology, while it is called ``ordinary'' or
 ``novel'' if it is of the form $M_1 \,.\, \Gamma$ with $M_1$ a maximal
 $\Gamma$-invariant subgroup of~$G_0$, as in the last two alternatives
 of Theorem~\ref{theo:maxsg1}, according to whether $M_1$ is a maximal
 or non-maximal subgroup of~$G_0$.
 Note that ``novel'' maximal subgroups~$M$ (for which $M_1$ is non-maximal)
 can only appear when $M_1^{}$ is contained in several maximal subgroups
 $M_1^{(i)}$ of~$G_0^{}$ that are not conjugate in~$G_0^{}$ but are abstractly
 isomorphic (as Lie groups) and whose conjugacy classes in~$G_0^{}$ are
 permuted transitively by~$\Gamma$: then none of them can be extended to
 a maximal subgroup of~$G$ and $M_1$ will be a maximal $\Gamma$-invariant
 subgroup of~$G_0$ which is ``smaller than expected'' since it is properly
 contained in each of the $M_1^{(i)}$.
\end{remark}

\begin{remark} \label{rmk:siginv}
 Note that the condition of $\Gamma$-invariance of a closed subgroup
 (or Lie subgroup) of~$G_0$ or of a subalgebra of~$\mathfrak{g}$ that
 plays a central role in this paper depends only on the ``outer automorphism
 part'' of~$\Gamma$, that is, on the image $\Sigma$ of~$\Gamma$ under the
 homomorphism~(\ref{eq:GREXT3}) or~(\ref{eq:GREXT4}), respectively.
 The full component group $\Gamma$ is an extension of $\Sigma$ by
 the kernel of the respective homomorphism, but for the purpose of
 deciding whether a closed subgroup (or Lie subgroup) $H_1$ of~$G_0$
 or a subalgebra $\mathfrak{h}$ of~$\mathfrak{g}$ is $\Gamma$-invariant,
 this kernel is irrelevant, so throughout the discussion in this section,
 we could equally well have used the term ``$\Sigma$-invariant'', instead
 of ``$\Gamma$-invariant'', where $\Sigma$ is the aforementioned subgroup
 of~$\mathrm{Out}(G_0)$ or~$\mathrm{Out}(\mathfrak{g})$, respectively. 
\end{remark}

\section{Quasiprimitive and primitive subalgebras}

In Theorem~\ref{theo:maxsg1}, the maximality condition for a closed
subgroup $M$ of a Lie group~$G$ which meets every connected component
of~$G$ is reduced to a maximality condition on its intersection $M_1$
with the connected one-component $G_0$ of~$G$, but what is still missing
is to translate that into an appropriate maximality condition on the
corresponding connected one-component~$M_0$, or equivalently, on the
corresponding Lie algebra~$\mathfrak{m}$~-- which, as we shall see,
is weaker than that of being maximal $\Gamma$-invariant or maximal
$\Sigma$-invariant and also somewhat more intricate.
This is the problem we shall turn to~now.
\vspace{2ex}

We begin with a preliminary definition that should be more or less obvious.
\begin{Definition} \label{def:sigsim1}
 Given a Lie algebra $\mathfrak{g}$ and a subgroup $\Sigma$
 of its outer auto\-morphism group $\mathrm{Out}(\mathfrak{g})$,
 we say that $\mathfrak{g}$ is \emph{$\Sigma$-simple} if it is not
 abelian and contains no non-trivial proper $\Sigma$-invariant ideal.
 When $\Sigma$ is the image under the homomorphism~(\ref{eq:GREXT4})
 of the component group $\Gamma$ of a Lie group~$G$ with Lie algebra
 $\mathfrak{g}$, we also use the term ``$\Gamma$-simple'' as a
 synonym for ``$\Sigma$-simple''.
\end{Definition}
For a correct understanding of this terminology, it is useful to
note that an ideal of~$\mathfrak{g}$ is, by definition, invariant
under inner automorphisms of~$\mathfrak{g}$ and so its invariance
under an arbitrary automorphism of~$\mathfrak{g}$ depends only on
the equivalence class of the latter, as an element of~$\mathrm{Out}%
(\mathfrak{g})$.
\begin{proposition} \label{prop:sigsim}
 Let\/ $\mathfrak{g}$ be a Lie algebra and\/ $\Sigma$ be a subgroup
 of its outer auto\-morphism group\/ $\mathrm{Out}(\mathfrak{g})$.
 Then\/ $\mathfrak{g}$ is\/ $\Sigma$-simple if and only if it is the
 direct sum of various copies of the same simple Lie algebra which
 under the action of\/~$\Sigma$ are permuted transitively among
 themselves.
\end{proposition}
\begin{proof}
 The ``if'' part is obvious, whereas the arguments for the ``only if''
 part are very similar to those used in the last part of the proof of
 Theorem~\ref{theo:maxsg1} above:
 \begin{itemize}
  \item If $\mathfrak{g}$ is neither solvable nor semisimple, then
        its radical is a non-trivial proper $\Sigma$-invariant ideal
        of~$\mathfrak{g}$. \vspace{-1ex}
  \item If $\mathfrak{g}$ is solvable but not abelian, then its
        commutator subalgebra is a non-trivial proper $\Sigma$-%
        invariant ideal of~$\mathfrak{g}$. \vspace{-1ex}
  \item If $\mathfrak{g}$ is semisimple, then decomposing it into the
        direct sum of its simple ideals and assembling these into
        $\Sigma$-orbits shows that if there is more than one such
        orbit, any one of them will be a non-trivial proper $\Sigma$-%
        invariant ideal of~$\mathfrak{g}$. \vspace{-4ex}
 \end{itemize}
\end{proof}
\noindent
A similar terminology is used for connected Lie groups.
\begin{Definition} \label{def:sigsim2}
 Given a connected Lie group $G_0$ and a subgroup $\Sigma$ of its
 outer automorphism group $\mathrm{Out}(G_0)$, we say that $G_0$
 is \emph{$\Sigma$-simple} if it is not abelian and contains no
 non-trivial proper connected $\Sigma$-invariant normal Lie subgroup.
 When $\Sigma$ is the image under the homomorphism~(\ref{eq:GREXT3})
 of the component group $\Gamma$ of a Lie group~$G$ with connected
 one-component $G_0$, we also use the term ``$\Gamma$-simple'' as a
 synonym for ``$\Sigma$-simple''.
\end{Definition}

To describe the maximality condition for subalgebras that we are after,
let us assume that $\mathfrak{g}$ is a Lie algebra and $\Sigma$ is a
subgroup of its outer automorphism group $\mathrm{Out}(\mathfrak{g})$,
as before, and that $\mathfrak{h}$ is a subalgebra of~$\mathfrak{g}$.
These data give rise to a subgroup
\begin{equation} \label{eq:AUTSIG}
 \mathrm{Aut}_\Sigma^{}(\mathfrak{g},\mathfrak{h})~
 =~\mathrm{Aut}_\Sigma^{}(\mathfrak{g}) \cap
   \mathrm{Aut}(\mathfrak{g},\mathfrak{h})
\end{equation}
of the automorphism group $\mathrm{Aut}(\mathfrak{g})$ of~$\mathfrak{g}$,
where $\mathrm{Aut}_\Sigma^{}(\mathfrak{g})$ denotes the inverse image of
$\Sigma$ under the canonical projection from $\mathrm{Aut}(\mathfrak{g})$
to~$\mathrm{Out}(\mathfrak{g})$ (i.e., the group of automorphisms of~%
$\mathfrak{g}$ that project to~$\Sigma$) and
\begin{equation} \label{eq:AUTSUB}
 \mathrm{Aut}(\mathfrak{g},\mathfrak{h})~
 =~\{ \, \phi \in \mathrm{Aut}(\mathfrak{g}) \; |~
         \phi(\mathfrak{h}) = \mathfrak{h} \, \}~.
\end{equation}
denotes the stability group of $\mathfrak{h}$ in~$\mathrm{Aut}(\mathfrak{g})$.
When $\Sigma$ is trivial ($\Sigma = \{1\}$), this becomes
\begin{equation} \label{eq:INTSUB}
 \mathrm{Inn}(\mathfrak{g},\mathfrak{h})~
 =~\{ \, \phi \in \mathrm{Inn}(\mathfrak{g}) \; |~
         \phi(\mathfrak{h}) = \mathfrak{h} \, \}~.
\end{equation}
In what follows, we shall only be interested in $\Sigma$-invariant subalgebras
and thus note that, according to Remarks~\ref{rmk:gaminv} and~\ref{rmk:siginv},
the condition that $\mathfrak{h}$ should be $\Sigma$-invariant means that every
element of $\Sigma$ should have a representative in~$\mathrm{Aut}_\Sigma^{}%
(\mathfrak{g})$ which leaves $\mathfrak{h}$ invariant, or in other words,
that under the canonical projection from $\mathrm{Aut}(\mathfrak{g})$
to~$\mathrm{Out}(\mathfrak{g})$, the subgroup $\mathrm{Aut}_\Sigma^{}%
(\mathfrak{g},\mathfrak{h})$ should be mapped \emph{onto} the sub%
group~$\Sigma\,$: this will ensure that, for any other subalgebra
$\tilde{\mathfrak{h}}$ of~$\mathfrak{g}$, $\mathrm{Aut}_\Sigma^{}%
(\mathfrak{g},\mathfrak{h})$-invariance is a stronger condition
than $\Sigma$-invariance.
(Of~course, by the very definition of~$\mathrm{Aut}_\Sigma^{}%
(\mathfrak{g},\mathfrak{h})$, $\mathfrak{h}$ itself is always
$\mathrm{Aut}_\Sigma^{}(\mathfrak{g},\mathfrak{h})$-invariant.)
Note also that in the special case where $\mathfrak{h}$ is a
$\Sigma$-invariant ideal, we have $\mathrm{Aut}_\Sigma^{}%
(\mathfrak{g},\mathfrak{h}) = \mathrm{Aut}_\Sigma(\mathfrak{g}) \,$
(independently of~$\mathfrak{h}$), so that for any other subalgebra
$\tilde{\mathfrak{h}}$ of~$\mathfrak{g}$, $\mathrm{Aut}_\Sigma^{}%
(\mathfrak{g},\mathfrak{h})$-invariance of~$\tilde{\mathfrak{h}}$
is then equivalent to $\tilde{\mathfrak{h}}$~being a $\Sigma$-%
invariant ideal as well.
\begin{Definition} \label{def:primsa}
 Let $\mathfrak{g}$ be a Lie algebra and $\Sigma$ be a subgroup
 of its outer auto\-morphism group $\mathrm{Out}(\mathfrak{g})$.
 A \emph{$\Sigma$-quasiprimitive subalgebra} of~$\mathfrak{g}$ is a
 proper $\Sigma$-invariant subalgebra $\mathfrak{m}$ of~$\mathfrak{g}$
 which is maximal among all $\mathrm{Aut}_\Sigma^{}(\mathfrak{g},%
 \mathfrak{m})$-invariant subalgebras of~$\mathfrak{g}$, that is,
 such that if $\tilde{\mathfrak{m}}$ is any $\mathrm{Aut}_\Sigma^{}%
 (\mathfrak{g},\mathfrak{m})$-invariant subalgebra of~$\mathfrak{g}$
 with $\, \mathfrak{m} \subset \tilde{\mathfrak{m}} \subset
 \mathfrak{g}$, then $\, \tilde{\mathfrak{m}} = \mathfrak{m} \,$
 or $\, \tilde{\mathfrak{m}} = \mathfrak{g}$.
 A \emph{$\Sigma$-primitive subalgebra} of~$\mathfrak{g}$ is a $\Sigma$-%
 quasiprimitive sub\-algebra of~$\mathfrak{g}$ which contains no non-trivial
 proper $\Sigma$-invariant ideal of~$\mathfrak{g}$.
 When $\Sigma$ is the image under the homomorphism~(\ref{eq:GREXT4})
 of the component group $\Gamma$ of a Lie group~$G$ with Lie algebra
 $\mathfrak{g}$, we also use the term ``$\Gamma$-(quasi)primitive''
 as a synonym for ``$\Sigma$-(quasi)primitive'', and when $\Sigma$
 is trivial ($\Sigma = \{1\}$), we omit the prefix and simply speak
 of a (quasi)primitive subalgebra.
\end{Definition}
\begin{remark}
 Note that the normalizer $\mathfrak{n}$ of~$\mathfrak{h}$ in~$\mathfrak{g}$
 is $\mathrm{Aut}_\Sigma^{}(\mathfrak{g},\mathfrak{h})$-invariant, for any
 choice of~$\Sigma$, because it is even $\mathrm{Aut}(\mathfrak{g},%
 \mathfrak{h})$-invariant. Indeed, for $\, \phi \in \mathrm{Aut}%
 (\mathfrak{g},\mathfrak{h})$ and $\, X \in \mathfrak{n}$, we have
 \[
  Y \in \mathfrak{h}~~\Longrightarrow~~
  \mathrm{ad}(\phi(X))(Y)~=~[\phi(X),\phi(\phi^{-1}(Y))]~
  =~\phi \big( \mathrm{ad}(X)(\phi^{-1}(Y)) \big) \in \mathfrak{h}~,
 \]
 that is, $\phi(X) \in \mathfrak{n}$.
 Therefore, a $\Sigma$-quasiprimitive subalgebra is either self-normalizing
 or is an ideal.
 In the second case, the observation immediately preceding Definition~%
 \ref{def:primsa} shows that the condition of being $\Sigma$-quasiprimitive
 is equivalent to that of being a \mbox{maximal} $\Sigma$-invariant ideal.
 Thus a $\Sigma$-primitive subalgebra is always self-normalizing.
 Moreover, the concepts of $\Sigma$-quasiprimitive subalgebra and
 $\Sigma$-primitive subalgebra \mbox{coincide} when the ambient
 Lie algebra is $\Sigma$-simple. 
\end{remark}
\begin{remark}
 The position of a $\Sigma$-quasiprimitive subalgebra~$\mathfrak{m}$
 of a Lie algebra~$\mathfrak{g}$ relative to the center $\mathfrak{z}$
 of~$\mathfrak{g}$ is strongly restricted, since using the fact
 that their sum $\, \mathfrak{m} + \mathfrak{z} \,$ is obviously
 an $\mathrm{Aut}_\Sigma^{}(\mathfrak{g},\mathfrak{m})$-invariant
 subalgebra of~$\mathfrak{g}$, we immediately conclude that either
 $\mathfrak{m}$ contains $\mathfrak{z}$ (this is necessarily the
 case when $\mathfrak{m}$ is self-normalizing, since any self-%
 normalizing subalgebra contains the center) or else $\mathfrak{m}$
 is a maximal $\Sigma$-invariant ideal of~$\mathfrak{g}$ which,
 together with $\mathfrak{z}$ which is also a $\Sigma$-invariant
 ideal of~$\mathfrak{g}$, spans all of~$\mathfrak{g}$, i.e.,
 $\mathfrak{g} = \mathfrak{m} + \mathfrak{z} \,$ (note that
 this sum is not necessarily direct).
 However, it follows from Theorem~\ref{theo:maxsg1} that the second
 case (which may occur even when $\mathfrak{m}$ is primitive and in
 this context corresponds to the so-called affine primitive examples
 mentioned in~\cite[Proposition~2.3]{Go}) does not arise in the study
 of maximal subgroups, and hence it will be discarded in what follows.%
 \footnote{A nice example is provided by considering $\, \mathfrak{m}
 = \mathfrak{su}(n) \,$ as a primitive subalgebra of $\, \mathfrak{g}
 = \mathfrak{u}(n) \,$: there is no maximal subgroup of $U(n)$ which
 has $SU(n)$ as its connected one-component, since we may extend the
 latter to closed subgroups with arbitrarily large finite component
 groups $\mathbb{Z}_p$ by multiplying matrices in $SU(n)$ with
 $p$-th roots of unity.}
\end{remark}

So far, we have not imposed any ``a priori'' restrictions on the choice
of~$\Sigma$, but the idea is of course that $\Sigma$ should be discrete.
(For example, if a Lie algebra $\mathfrak{g}$ is $\Sigma$-simple, then
due to Proposition~\ref{prop:sigsim}, $\mathfrak{g}$ is necessarily
semisimple, so $\mathrm{Out}(\mathfrak{g})$ is discrete, and hence
so is~$\Sigma$.)
This is also what is required in order to make contact with the situation
encountered in the previous section, where $G$ is a Lie group with connected
one-component $G_0$, component group $\Gamma$ and Lie algebra~$\mathfrak{g}$
and where $\, \Sigma \subset \mathrm{Out}(\mathfrak{g})$ is the
image of~$\Gamma$ under the homomorphism~(\ref{eq:GREXT4}). Indeed, take
the adjoint representation of~$G$, viewed as a Lie group homomorphism
from~$G$ to $GL(\mathfrak{g})$ which we shall denote by~$\mathrm{Ad}$
or sometimes by~$\mathrm{Ad}_G$ in order to avoid ambiguities whenever
this seems con\-venient: then the very definition of the homomorphism~%
(\ref{eq:GREXT4}), together with the fact that, also by definition,
the group of inner automorphisms of~$\mathfrak{g}$ is
\begin{equation} \label{eq:GREXT7}
 \mathrm{Inn}(\mathfrak{g})~=~\mathrm{Ad}(G_0)~,
\end{equation}
implies that applying $\mathrm{Ad}$ to the short exact sequence~%
(\ref{eq:GREXT1}) gives the corresponding

\pagebreak

\noindent
adjoint short exact sequence in the form
\begin{equation} \label{eq:GREXT8}
 \{1\}~\longrightarrow~\mathrm{Ad}(G_0)~\longrightarrow~\mathrm{Ad}(G)~
       \longrightarrow~\Sigma~\longrightarrow~\{1\}~.
\end{equation}
This shows that there is a Lie group (even a linear Lie group) which has
$\Sigma$ as its component group, namely the adjoint group $\mathrm{Ad}(G)$,
and also implies that
\begin{equation} \label{eq:GREXT9}
 \mathrm{Aut}_\Sigma^{}(\mathfrak{g})~=~\mathrm{Ad}(G)~.
\end{equation}
Suppose further that $\mathfrak{h}$ is a $\Sigma$-invariant sub%
algebra of~$\mathfrak{g}$, and let $N$ be its normalizer in~$G$.
Then
\begin{equation} \label{eq:GREXT10}
 \mathrm{Aut}_\Sigma^{}(\mathfrak{g},\mathfrak{h})~=~\mathrm{Ad}_G(N)~.
\end{equation}
What seems to be more difficult is to formulate, for an arbitrary Lie
algebra $\mathfrak{g}$, reasonably simple conditions that would allow
to state precisely which subgroups $\Sigma$ of the group $\mathrm{Out}%
(\mathfrak{g})$ can arise from a Lie group $G$ having $\mathfrak{g}$
as its Lie algebra, especially in those cases where $\mathrm{Inn}%
(\mathfrak{g})$ fails to be closed in $\mathrm{Aut}(\mathfrak{g})$
and hence $\mathrm{Out}(\mathfrak{g})$ is not even a Lie group
(see~\cite[Chap.~2, Exercise D3]{He} for an example). But we
have not attempted to do so since we are ultimately interested
in compact Lie groups, for which the answer is immediate:
$\Sigma$ must be finite.
\vspace{2ex}

The importance of the concept of a quasiprimitive subalgebra in the
present context stems from its relation to that of a maximal subgroup.
However, a precise formulation of this relation depends on additional
hypotheses, one of which is that the Lie algebra $\mathfrak{g}$
be reductive.

\section{Reductive Lie algebras and groups}

In order to continue our analysis, we shall now impose an important
restriction: we consider only reductive Lie algebras and groups.
Recall that, according to one of various possible definitions (all
of which are equivalent), a Lie algebra $\mathfrak{g}$ is reductive
if it can be written as the direct sum of its center $\mathfrak{z}$
and a semisimple ideal; the latter is then necessarily equal to
the derived subalgebra $\mathfrak{g}'$ of~$\mathfrak{g}$, or what
is the same, the commutator subalgebra $[\mathfrak{g},\mathfrak{g}]$
of~$\mathfrak{g}$.
As is well known, any compact Lie algebra (i.e., any Lie algebra
that can be obtained as the Lie algebra of some compact Lie group)
is reductive, and of course so is any semisimple Lie algebra.
\linebreak
For Lie groups, the situation is somewhat more complicated
since there is a certain amount of ambiguity in the literature
as to what are the precise conditions required of a reductive
Lie group, over and above the obvious one that its Lie algebra
should be reductive.
(Another standard requirement is that its component group
should be finite.)
We shall return to this subject below and then give a precise
definition (see Definition~\ref{def:redlgr}), but before that,
we want to state and prove the basic theorem about the relation
between quasiprimitive subalgebras and maximal subgroups in the
reductive setting, which does not depend on any of the more
sophisticated aspects of that definition.
\vspace{2ex}

As a preliminary step, we specify criteria to ensure that certain
connected normal subgroups of certain connected Lie groups are
automatically closed.

\pagebreak

\begin{lemma} \label{lem:CLNSG}
 Let\/ $G_0$ be a connected Lie group with Lie algebra\/~$\mathfrak{g}$
 and let\/ $H_0$ be a connected Lie subgroup of\/~$G_0$ with Lie
 algebra\/~$\mathfrak{h}$. Suppose that\/ $\mathfrak{h}$ is an ideal
 of\/~$\mathfrak{g}$ and, in addition, either one of the following
 two hypotheses is satisfied:
 \begin{enumerate}
  \item $G_0$ is simply connected.
  \item $\mathfrak{g}$ is reductive and\/ $\mathfrak{h}$ contains
        the center\/~$\mathfrak{z}$ of\/~$\mathfrak{g}$.
 \end{enumerate}
 Then\/ $H_0$ is closed in\/~$G_0$.
\end{lemma}
\begin{proof}
 Following the idea outlined in~\cite[p.~127]{Che}, consider the quotient
 algebra $\, \mathfrak{l} = \mathfrak{g}/\mathfrak{h} \,$ and assume $L$
 to be any connected Lie group which has Lie algebra~$\mathfrak{l}$.
 Suppose we can lift the canonical projection from~$\mathfrak{g}$
 to~$\mathfrak{l}$, which is a Lie algebra homomorphism $\, f:
 \mathfrak{g} \longrightarrow \mathfrak{l}$, to a Lie group
 homomorphism $\, F: G_0 \longrightarrow L$.
 Then it is obvious that $H_0$ is equal to the connected one-component
 of the kernel of~$F$ (since both are generated by the same subalgebra
 of~$\mathfrak{g}$, namely~$\mathfrak{h}$), and hence is closed in~$G_0$.
 This proves the first statement, since the desired lifting
 of~$f$ to~$F$ always exists if~$G_0$ is simply connected.
 If not, let $\tilde{G}_0$ be the universal covering group
 of~$G_0$, with covering homomorphism $\, \pi: \tilde{G}_0
 \longrightarrow G_0 \,$ whose kernel will be denoted by
 $\tilde{N}$, so that $\, G_0 = \tilde{G}_0/\tilde{N}$,
 and lift the Lie algebra homomorphism $\, f: \mathfrak{g}
 \longrightarrow \mathfrak{l} \,$ to a Lie group homomorphism
 $\, \tilde{F}: \tilde{G}_0 \longrightarrow L$.
 In order for the previous argument to continue working,
 we must show that $\tilde{F}$ is trivial on~$\tilde{N}$,
 since this guarantees that it factors to yield a Lie group
 homomorphism $\, F: G_0 \longrightarrow L$, as before.
 Now $\tilde{N}$ is a discrete normal subgroup of~$\tilde{G}_0$
 contained in its center~$Z(\tilde{G}_0)$, and $\tilde{F}$, being
 surjective since $f$ is surjective and~$L$ is connected, maps the
 center~$Z(\tilde{G}_0)$ of~$\tilde{G}_0$ to the center $Z(L)$ of~$L$.
 Thus it suffices to show that~$L$ can be chosen to be centerfree.
 This is guaranteed by taking $L$ to be the adjoint group
 of~$\mathfrak{l}$, provided $\mathfrak{l}$ is centerfree.
 But if $\mathfrak{g}$ is reductive and $\mathfrak{h}$ is an
 ideal of~$\mathfrak{g}$ containing its center $\mathfrak{z}$,
 then $\; \mathfrak{l} = \mathfrak{g}/\mathfrak{h} \cong
 [\mathfrak{g},\mathfrak{g}]/([\mathfrak{g},\mathfrak{g}]
 \cap \mathfrak{h}) \;$ is semisimple and hence centerfree.
\end{proof}
\begin{remark} \label{rmk:CLNSG}
 It may be worthwhile to point out that the conclusion of the second
 statement in Lemma~\ref{lem:CLNSG} would become wrong if we were to
 omit any one of its hypotheses.
 The source of the counterexamples is always the same,
 namely, the irrational flow on the torus.
 \begin{enumerate}
  \item Take $G_0$ to be the two-dimensional torus~$\mathbb{T}^2$ itself.
        Of course, $G_0$ is abelian and hence any subgroup $H_0$ of~$G_0$
        is normal.
        Thus taking~$H_0$ to be the irrational flow on this torus, we
        obtain a connected normal Lie subgroup of~$G_0$ which is dense
        in~$G_0$, rather than closed, even though $\mathfrak{g}$, being
        abelian, is reductive.
        However, $\mathfrak{h}$ does not contain the center $\mathfrak{z}$
        of~$\mathfrak{g}$, since $\, \mathfrak{z} = \mathfrak{g}$.
  \item Take $G_0$ to be the group of upper triangular complex
        $(3 \times 3)$-matrices of determinant~$1$ and diagonal
        entries from the unit circle~$\mathbb{T}$: its commutator
        subgroup $G_0^\prime$ is the group of upper triangular
        complex $(3 \times 3)$-matrices with diagonal entries
        all equal to~$1$.
        Here, $G_0$ is solvable, with discrete center $\, Z(G_0)
        \cong \mathbb{Z}_3$, whereas $G_0^\prime$ is nilpotent,
        and by the very definition of the commutator subgroup, any
        subgroup $H_0$ of~$G_0$ containing~$G_0^\prime$ is normal. 
        In~fact, the qoutient group $G_0/G_0^\prime$ is the two-%
        dimensional torus~$\mathbb{T}^2$, so taking $H_0$ to be
        the inverse image of the irrational flow on this torus
        under the canonical projection from~$G_0$ to~$G_0/G_0^\prime$,
        we obtain a connected normal Lie subgroup of~$G_0$ which
        is dense in~$G_0$, rather than closed, even though
        $\mathfrak{h}$ contains the center~$\mathfrak{z}$
        of~$\mathfrak{g}$, since $\, \mathfrak{z} = 0$.
        However, $\mathfrak{g}$ is not reductive, since it
        is solvable without being abelian.
 \end{enumerate}
\end{remark}
\noindent
Now we are ready to formulate the first main theorem of this section.
\begin{theorem} \label{theo:maxsg2}
 Let\/ $G$ be a Lie group with connected one-component\/~$G_0$, \linebreak
 component group\/~$\Gamma$ and Lie algebra\/~$\mathfrak{g}$, and let\/ $M$
 be a closed subgroup of\/~$G$ with connected one-component\/~$M_0$ and Lie
 algebra\/~$\mathfrak{m}$ such that\/ $M$ meets every connected component
 of\/~$G$.
 Set $\, M_1 = M \cap G_0$, so $\, M_0 \subset M_1 \subset M$.
 Moreover, suppose that the Lie algebra\/ $\mathfrak{g}$ is reductive
 and the quotient group\/ $M_1/M_0$ is finite. \linebreak
 Then\/ $M$ is a maximal subgroup of\/~$G$, or equivalently, $M_1$ is
 a maximal\/ $\Gamma$-in\-variant subgroup of\/~$G_0$, if and only
 if\/ $\mathfrak{m}$ is a\/ $\Gamma$-quasiprimitive subalgebra
 of\/~$\mathfrak{g}$ containing its center~$\mathfrak{z}$ and,
 more specifically, one of the following two alternatives holds:
 \begin{itemize}
  \item \textbf{Normal type:\,} $\mathfrak{m}$ is a maximal\/ $\Gamma$-%
        invariant ideal of\/~$\mathfrak{g}$, or equivalently, $M_0$~is~a
        maximal\/ $\Gamma$-invariant connected normal subgroup of\/~$G_0$,
        and\/ $M_1/M_0$ is a finite maximal\/ $\Gamma$-invariant subgroup
        of\/~$G_0/M_0$: from these data, $M_1$ is recovered by taking
        the inverse image under the canonical projection from\/~$G_0$
        to\/~$G_0/M_0$ while\/ $M$ is recovered using Proposition~%
        \ref{prop:invsg}.
  \item \textbf{Normalizer type:\,} $\mathfrak{m}$ is a self-normalizing\/
        $\Gamma$-quasiprimitive subalgebra of\/~$\mathfrak{g}$: then\/ $M_1$
        is its normalizer in\/~$G_0$ while\/ $M$ is its normalizer in\/~$G$.
 \end{itemize}
\end{theorem}
\begin{proof}
 This is naturally divided into two parts. \\[1mm]
 The first part of the proof consists in showing that if $M$ is
 maximal, or equivalently, $M_1$ is maximal $\Gamma$-invariant,
 then $\mathfrak{m}$ must be $\Gamma$-quasiprimitive.
 First of all, \linebreak it is clear from Theorem~\ref{theo:maxsg1}
 that $\mathfrak{m}$ is $\Gamma$-invariant and must either be an ideal
 of~$\mathfrak{g}$ or else be self-normalizing, depending on whether
 $M_0$ is normal or not.
 Now~suppose that $\tilde{\mathfrak{m}}$ is an
 $\mathrm{Aut}_\Sigma^{}(\mathfrak{g},\mathfrak{m})$-%
 invariant subalgebra of~$\mathfrak{g}$ containing~%
 $\mathfrak{m}$, where $\Sigma$ is the image of~$\Gamma$
 under the homomorphism~(\ref{eq:GREXT4}), as usual:
 we want to show that either $\, \tilde{\mathfrak{m}}
 = \mathfrak{m} \,$ or $\, \tilde{\mathfrak{m}} = \mathfrak{g}$.
 This will be done in three steps.
 The first step consists in noting that in the case of the $\Gamma$-%
 irrational type, where $G_0/M_0$ is $\mathbb{R}^n$ (as an abelian
 Lie group) and hence $\mathfrak{g}/\mathfrak{m}$ is $\mathbb{R}^n$
 (as an abelian Lie algebra), this is trivial, since an intermediate
 $\Gamma$-invariant subalgebra between $\mathfrak{m}$ and $\mathfrak{g}$
 would correspond to a non-trivial proper $\Gamma$-invariant subspace
 of~$\mathfrak{g}/\mathfrak{m}$, and there are no such subspaces
 because $\Gamma$ acts irreducibly on $G_0/M_0$ and hence also
 on~$\mathfrak{g}/\mathfrak{m}$.
 Therefore, according to Theorem~\ref{theo:maxsg1}, we may assume,
 without loss of generality, that $M_1$ is self-normalizing and
 contains the center of~$G_0$, implying that $\mathfrak{m}$
 contains the center~$\mathfrak{z}$ of~$\mathfrak{g}$. 
 The second step consists in showing that we may, without
 loss of generality, assume $\tilde{\mathfrak{m}}$ to be
 a $\Gamma$-invariant ideal of~$\mathfrak{g}$.
 Indeed, if $\mathfrak{m}$ is a $\Gamma$-invariant ideal of~%
 $\mathfrak{g}$, $\mathrm{Aut}_\Sigma^{}(\mathfrak{g},\mathfrak{m})$
 is equal to $\mathrm{Ad}(G)$ and contains $\mathrm{Ad}(G_0)$, so the
 conclusion is immediate.
 If on the other hand $\mathfrak{m}$ is a $\Gamma$-invariant self-%
 normalizing subalgebra of~$\mathfrak{g}$, $\mathrm{Aut}_\Sigma^{}%
 (\mathfrak{g},\mathfrak{m})$ is equal to $\mathrm{Ad}_G(M)$ and
 contains $\mathrm{Ad}_G(M_1)$ since in this case, according to
 Theorem~\ref{theo:maxsg1}, $M$ is the normalizer of~$\mathfrak{m}$
 in~$G$ and $M_1$ is the normalizer of~$\mathfrak{m}$ in~$G_0$,
 so $\mathrm{Aut}_\Sigma^{}(\mathfrak{g},\mathfrak{m})$-invariance
 of~$\tilde{\mathfrak{m}}$ is equivalent to the statement that its
 normalizer $\tilde{N}$ in~$G$ contains~$M$ and also that its normalizer
 $\tilde{N}_1$ in~$G_0$ is $\Gamma$-invariant and contains $M_1$.
 Maximality of~$M$ now implies that either $\, \tilde{N} = M$,
 $\tilde{N}_1 = M_1 \,$ or $\, \tilde{N} = G$, $\tilde{N}_1 = G_0$.
 Thus if $\tilde{\mathfrak{n}}$ denotes the normalizer of~%
 $\tilde{\mathfrak{m}}$ in~$\mathfrak{g}$, which is just the
 Lie algebra of~$\tilde{N}$ and of~$\tilde{N}_1$, we conclude
 that in the first case, $\tilde{\mathfrak{m}} \subset
 \tilde{\mathfrak{n}} = \mathfrak{m}$, which is one of the two
 possible conclusions that we ultimately want to arrive at, while
 in the second case, $\tilde{\mathfrak{n}} = \mathfrak{g}$, i.e.,
 $\tilde{\mathfrak{m}}$ is an ideal of~$\mathfrak{g}$. \linebreak
 The third step consists in showing that, independently of
 whether $\mathfrak{m}$ is an ideal or is self-normalizing,
 $\tilde{\mathfrak{m}}$ being an ideal containing $\mathfrak{m}$,
 which in turn contains $\mathfrak{z}$, \linebreak leads to the
 conclusion that either $\, \tilde{\mathfrak{m}} = \mathfrak{m} \,$
 or $\, \tilde{\mathfrak{m}} = \mathfrak{g}$.
 To this end, let $\tilde{M}_0$ be the connected Lie subgroup of~$G_0$
 corresponding to~$\tilde{\mathfrak{m}}$ and let $\tilde{M}_1$ be the
 abstract subgroup of~$G_0$ generated by~$M_1$ and~$\tilde{M}_0$;
 explicitly, $\tilde{M}_1 = M_1 \tilde{M}_0 \,$ since $\tilde%
 {\mathfrak{m}}$ is an ideal and hence $\tilde{M}_0$ is normal.
 Now since $\mathfrak{g}$ is assumed to be reductive, we can
 apply Lemma~\ref{lem:CLNSG} to conclude that $\tilde{M}_0$
 is a closed subgroup of~$G_0$.
 Moreover, it follows that $\tilde{M}_1$ is also a closed subgroup
 of~$G_0$ and that it has connected one-component~$\tilde{M}_0$.
 [\,Indeed, since as we have just seen, $\tilde{M}_0$ is a closed
 normal subgroup of~$G_0$, the quotient $G_0/\tilde{M}_0$ is a
 well-defined connected Lie group and the image of~$M_1$ under
 the canonical projection from~$G_0$ to $G_0/\tilde{M}_0$ is
 the quotient group $M_1/(M_1 \cap \tilde{M}_0)$; let us denote
 it by $Q$, for the sake of brevity. Now according to the first
 and second isomorphism theorems of group theory, we have
 $\; \tilde{M}_1/\tilde{M}_0 \cong Q \cong (M_1/M_0)\,/\,
 ((M_1 \cap \tilde{M}_0)/M_0)$: the second isomorphism shows
 that $Q$ is a quotient group of~$M_1/M_0$, which is assumed
 to be finite, and hence $Q$ itself must be finite, whereas
 the first isomorphism shows that $\tilde{M}_1$ is precisely
 the inverse image of~$Q$ under the canonical projection
 from~$G_0$ to $G_0/\tilde{M}_0$. Therefore, $\tilde{M}_1$
 is a closed subgroup of~$G_0$ with connected one-component
 $\tilde{M}_0$ and component group~$Q$.]
 Maximality of~$M$ now implies that either $\, \tilde{M}_1 = M_1 \,$
 or $\, \tilde{M}_1 = G_0 \,$ and hence either $\, \tilde{\mathfrak{m}}
 = \mathfrak{m} \,$ or $\, \tilde{\mathfrak{m}} = \mathfrak{g}$, thus
 completing the proof that $\mathfrak{m}$ is $\Gamma$-quasi\-primitive.
 In particular, note that if $M$ and $M_1$ are of normal type, this
 means that $\mathfrak{m}$ is a maximal $\Gamma$-invariant ideal
 of~$\mathfrak{g}$ and hence $M_0$ is a maximal $\Gamma$-invariant
 connected normal subgroup of~$G_0$.
 \\[1mm]
 The second part of the proof is devoted to showing that the conditions
 stated in Theorem~\ref{theo:maxsg2} are not only necessary but also
 sufficient to guarantee that $M$ is a maximal subgroup of~$G$, or
 equivalently, $M_1$ is a maximal $\Gamma$-invariant subgroup of~$G_0$.
 This will be done separately for each of the two types.
 \begin{itemize}
  \item Normal type: First of all, it follows from Lemma~\ref{lem:CLNSG}
        that the connected normal Lie subgroup $\tilde{M}_0$ of~$G_0$
        corresponding to an ideal $\tilde{\mathfrak{m}}$ of~$\mathfrak{g}$
        containing its center~$\mathfrak{z}$ is automatically closed: this
        implies that as soon as $\, \mathfrak{z} \subset \mathfrak{m}$,
        $\mathfrak{m}$ is maximal $\Gamma$-invariant if and only if
        if $M_0$ is maximal $\Gamma$ invariant.
        Assuming this to be the case, and assuming $D_0$ to be a finite
        maximal $\Gamma$-invariant subgroup of~$G_0/M_0$, let $M_1$ be the
        inverse image of~$D_0$ under the canonical projection from~$G_0$
        to~$G_0/M_0$ and $M$ be the closed subgroup of~$G$ obtained
        from~$M_1$ according to Proposition~\ref{prop:invsg}.
        Now if $\tilde{M}_1$ is a $\Gamma$-invariant closed subgroup of~$G_0$
        containing $M_1$, it is clear that $\tilde{M_1}/M_0$ is a $\Gamma$-%
        invariant closed subgroup of $G_0/M_0$ containing $\, D_0 = M_1/M_0$.
        (Here, we just use the definition of the quotient topology, which
        states that a subset $A$ ($U$) of~$G_0/M_0$ is closed (open)
        in~$G_0/M_0$ if and only if its inverse image $\pi^{-1}(A)$
        ($\pi^{-1}(U)$) under the canonical projection $\pi$ from~$G_0$
        to~$G_0/M_0$ is closed (open) in~$G_0$. But $\tilde{M}_1$ is a
        subgroup containing~$M_1$ and hence~$M_0$, so that $\, \pi^{-1}
        (\tilde{M}_1/M_0) = \tilde{M}_1$.)
        Maximality of~$D_0$ now implies that either $\, \tilde{M_1}/M_0
        = M_1/M_0 \,$ or $\, \tilde{M_1}/M_0 = G_0/M_0 \,$, so either
        $\, \tilde{M_1} = M_1 \,$ or $\, \tilde{M_1} = G_0 \,$.
  \item Normalizer type: Suppose that $\mathfrak{m}$ is a self-normalizing
        $\Gamma$-quasiprimitive sub\-algebra of~$\mathfrak{g}$, $M_1$ is
        its normalizer in~$G_0$ and $M$ is its normalizer in~$G$.
        Now if $\tilde{M}$ is a closed subgroup of~$G$ containing~$M$ and
        $\tilde{\mathfrak{m}}$ is its Lie algebra, it is obvious that
        $\tilde{\mathfrak{m}}$ is an $\mathrm{Ad}_G(M)$-invariant
        subalgebra of~$\mathfrak{g}$ containing $\mathfrak{m}$.
        Thus in view of equation~(\ref{eq:GREXT10}), $\Gamma$-%
        quasiprimitivity of~$\mathfrak{m}$ implies that either
        $\, \tilde{\mathfrak{m}} = \mathfrak{m} \,$ or
        $\, \tilde{\mathfrak{m}} = \mathfrak{g} \,$ and
        hence either $\, \tilde{M} \subset M$, according to equations~%
        (\ref{eq:NORM3}) and~(\ref{eq:INCL1}), so in fact $\, \tilde{M} = M$,
        or else $\, \tilde{M} = G$.
        \\ \qedhere
 \end{itemize}
\end{proof}

\begin{remark}
 We emphasize here that the conclusions stated in Theorem~ \ref{theo:maxsg2}
 for the normal type are not independent, in the sense that the requirement
 on~$\mathfrak{m}$~-- namely the condition that $\mathfrak{m}$ should be a
 maximal $\Gamma$-invariant ideal of~$\mathfrak{g}$, which is equivalent
 to the condition that $\mathfrak{g}/\mathfrak{m}$ should not admit any
 non-trivial proper $\Gamma$-invariant ideals~-- is a necessary pre%
 requisite for the existence of \emph{finite} maximal $\Gamma$-invariant
 subgroups of~$G_0/M_0$.
 Indeed, if $\mathfrak{m}$ is not maximal, we can choose a non-trivial
 proper $\Gamma$-invariant ideal of~$\mathfrak{g}/\mathfrak{m}$, say
 $\mathfrak{h}$, and consider the connected normal Lie subgroup $H_0$
 of~$G_0/M_0$ corresponding to~$\mathfrak{h}$, noting that according
 to Lemma~\ref{lem:CLNSG}, $H_0$ is closed because $\mathfrak{g}$ is
 reductive and $\mathfrak{m}$ contains its center~$\mathfrak{z}$, so
 $\mathfrak{g}/\mathfrak{m}$ is semisimple, i.e., reductive with
 trivial center.
 Therefore, if $D_0$ is any given finite $\Gamma$-invariant subgroup
 of~$G_0/M_0$, the abstract subgroup $H$ of~$G_0/M_0$ generated by
 $D_0$ and~$H_0$, explicitly given by $\, H = D_0 H_0 \,$ since
 $\mathfrak{h}$ is an ideal and hence $H_0$ is normal, is a closed
 subgroup of $G_0/M_0$ with connected one-component~$H_0$ properly
 containing $D_0$ and properly contained in~$G_0/M_0$, which shows
 that $D_0$ cannot be maximal.
 (The fact that $H$ is also a closed subgroup of~$G_0/M_0$ and that
 it has connected one-component~$H_0$ is proved in exactly the same
 way as the corresponding statement in the first part of the proof
 of Theorem~\ref{theo:maxsg2}, replacing $G_0$ by~$G_0/M_0$, $M_0$
 by~$\{1\}$, $M_1$ by~$D_0$ and $\tilde{M}_0$ by~$H_0$.)
\end{remark}

Returning to the discussion of reductive Lie algebras and groups
initiated at the very beginning of this section, we proceed to
formulate the precise definition of the concept of a reductive
Lie group announced there:
\begin{Definition} \label{def:redlgr}
 Let $G$ be a Lie group with connected one-component $G_0$,
 com\-ponent group $\Gamma$ and Lie algebra $\mathfrak{g}$.
 Assuming that $\mathfrak{g}$ is a reductive Lie algebra, we say
 that $G$ is a \emph{reductive Lie group} if it comes equipped
 with the following additional structure: a compact subgroup~$K$
 of~$G$, an involutive automorphism $\theta$ of~$\mathfrak{g}$
 and an $\mathrm{Ad}(G)$-invariant as well as $\theta$-invariant
 non-degenerate symmetric bilinear form $B$ on~$\mathfrak{g}$
 such that if we write
 \begin{equation} \label{eq:CARDEC1}
  \mathfrak{g}~=~\mathfrak{k} \oplus \mathfrak{p}~,
 \end{equation}
 where $\mathfrak{k}$ and $\mathfrak{p}$ denote the eigenspaces
 of~$\theta$ for eigenvalue $+1$ and $-1$, respectively~-- which
 implies the commutation relations
 \begin{equation} \label{eq:CARDEC2}
  [\mathfrak{k},\mathfrak{k}] \subset \mathfrak{k}~~~,~~~
  [\mathfrak{k},\mathfrak{p}] \subset \mathfrak{p}~~~,~~~
  [\mathfrak{p},\mathfrak{p}] \subset \mathfrak{k}~,
 \end{equation}
 expressing the requirement that $\theta$ should be an automorphism,
 and forces $\mathfrak{k}$ and $\mathfrak{p}$ to be orthogonal
 under~$B$~-- we have
 \begin{enumerate}[(i)]
  \item $\mathfrak{k}$ is precisely the subalgebra of~$\mathfrak{g}$
        corresponding to the closed subgroup $K$ of~$G$,
  \item $B$ is negative definite on $\mathfrak{k}$ and positive definite
        on~$\mathfrak{p}$,
  \item the map
        \begin{equation}
         \begin{array}{ccc}
          K \times\, \mathfrak{p} & \longrightarrow &      G      \\[1mm]
                   (k,X)          &   \longmapsto   & k \, \exp(X)
         \end{array}
        \end{equation}
        is a global diffeomorphism,
  \item the derived subgroup $G_0^{\,\prime}$ of~$G_0$ (defined as the
        connected Lie subgroup of~$G_0$ corresponding to the derived
        subalgebra $\mathfrak{g}'$ of~$\mathfrak{g}$) has finite center.
 \end{enumerate}
 As in the semisimple case, $K$ is called the \emph{maximal compact
 subgroup} of~$G$ (its maxi\-mality follows easily from condition~(iii)
 above), $\theta$ is called the \emph{Cartan involution}, the direct
 decomposition~(\ref{eq:CARDEC1}) is called the \emph{Cartan
 decomposition} and $B$ is called the \emph{invariant bilinear form}.
\end{Definition}
Note that this definition is slightly weaker than that of a ``reductive
Lie group in the Harish-Chandra class'' given in Chapter~7, Section~2 of
Ref.~\cite{Kn}, in that among the conditions~(i)-(vi) listed there,
we impose conditions~(i)-(iv) and~(vi) but discard condition~(v).
There are various advantages that speak in favor of this omission.
The basic one is that any compact Lie group and any semisimple Lie
group with \mbox{finite} component group\,%
\footnote{We do not adhere to the convention adopted in Ref.~\cite{Kn}
according to which semisimple Lie groups are automatically supposed to be
connected.}
and finite center is reductive in our sense: in particular, this avoids
an annoying flaw of the definition given in Ref.~\cite{Kn}, according
to which compact Lie groups $G$ whose component group $\Gamma$ acts on
their Lie algebra~$\mathfrak{g}$ by outer automorphisms, such as the
full orthogonal groups in even dimensions, are not reductive.
More precisely, it excludes all compact Lie groups for which the
image $\Sigma$ of~$\Gamma$ under the homomorphism~(\ref{eq:GREXT4}),
which according to equation~(\ref{eq:GREXT8}) is also the component
group of~$\mathrm{Ad}(G)$, is non-trivial, i.e., $\Sigma \neq \{1\}$:
such groups are of great interest for our work and we cannot see any
convincing reason for excluding them from the category of reductive
Lie groups. \linebreak
Fortunately, many of the important structural results on \mbox{reductive}
Lie groups stated in Ref.~\cite{Kn}, especially Propositions~7.19,
7.20 and~7.21, remain valid in our slightly broader context.
In particular, $G_0^{\,\prime}$ is a closed subgroup of~$G_0$.
More~generally, we have the following analogue of Lemma~\ref{lem:CLNSG},
whose proof is however much less elementary.
\begin{proposition} \label{prop:CLSUBG}
 Let\/ $G$ be a reductive Lie group with connected one-component\/~$G_0$,
 component group\/~$\Gamma$ and Lie algebra\/~$\mathfrak{g}$. Then for
 any semisimple ideal of\/~$\mathfrak{g}$ and, more generally, for any
 $\theta$-invariant semisimple subalgebra of\/~$\mathfrak{g}$, the
 corresponding connected Lie subgroup of\/~$G_0$ is closed and has
 finite center.
\end{proposition}
\begin{proof}
 For the second statement, see Proposition~7.20~(b) of Ref.~\cite{Kn}.
 Using this fact, the proof of the first statement for the general case
 is identical with that for the special case of the derived subgroup
 $G_0^{\,\prime}$ of~$G_0$, which is the content of Proposition~7.20~(a)
 of Ref.~\cite{Kn}.
\end{proof}

\noindent
In particular, this proposition can be applied to the simple ideals of~%
$\mathfrak{g}$ as well as to the $\Gamma$-simple ideals of~$\mathfrak{g}$
(each of which is, according to Proposition~\ref{prop:sigsim}, the direct
sum of a certain number of copies of the same simple ideal, transitively
permuted among themselves under the action of~$\Gamma$) to deduce that
the connected Lie subgroups of~$G_0^{}$ corresponding to these ideals
are closed subgroups of~$G_0^{}$ with finite center: they will be called
the \emph{simple factors} and the \emph{$\Gamma$-simple factors} of~$G_0$,
respectively. Note that the latter, being $\Gamma$-invariant, extend to
closed subgroups of~$G$, each of which meets every connected component
of~$G$.

In what follows, our main goal will be to show how, very roughly
speaking, the study of the maximal subgroups of a reductive Lie group
can be reduced to that of the maximal subgroups of its simple factors.
We start at the Lie algebra level, where this reduction can be achieved
through the following two theorems.
\begin{theorem} \label{theo:maxisa1}
 Let\/ $\mathfrak{g}$ be a Lie algebra and\/ $\Sigma$ be a subgroup
 of its outer auto\-morphism group\/ $\mathrm{Out}(\mathfrak{g})$.
 Assuming that\/ $\mathfrak{g}$ is reductive, with canonical
 decomposition
 \begin{equation} \label{eq:DECRLA1}
  \mathfrak{g}~=~\mathfrak{z} \; \oplus \, \mathfrak{g}'~~~,~~~
  \mathfrak{g}'~=~\mathfrak{g}_1^\Sigma \, \oplus \ldots \oplus \,
                  \mathfrak{g}_r^\Sigma
 \end{equation}
 into its center\/~$\mathfrak{z}$, its derived algebra\/~$\mathfrak{g}'$
 and its\/ $\Sigma$-simple ideals\/ $\mathfrak{g}_1^\Sigma \,,\ldots,\,
 \mathfrak{g}_r^\Sigma$, let\/ $\mathfrak{m}$ be a\/ $\Sigma$-quasi%
 primitive subalgebra of\/~$\mathfrak{g}$ containing\/~$\mathfrak{z}$.
 Then denoting the intersection of\/~$\mathfrak{m}$ with\/~$\mathfrak{g}'$
 by\/~$\mathfrak{m}'$, so that\/ $\, \mathfrak{m} = \mathfrak{z} \oplus
 \mathfrak{m}'$, one of the following two alternatives holds.
 \begin{enumerate}
  \item $\mathfrak{m}'$ is the direct sum of all\/ $\Sigma$-simple ideals
        of\/~$\mathfrak{g}$ except one, say\/~$\mathfrak{g}_i^\Sigma$, plus
        a\/ $\Sigma$-primitive subalgebra\/ $\mathfrak{m}_i^\Sigma$
        of\/~$\mathfrak{g}_i^\Sigma$:
        \begin{equation} \label{eq:MAXISA2}
         \mathfrak{m}'~=~\mathfrak{g}_1^\Sigma \, \oplus \ldots
                         \oplus \, \mathfrak{g}_{i-1}^\Sigma \,
                         \oplus \, \mathfrak{m}_i^\Sigma \,
                         \oplus \, \mathfrak{g}_{i+1}^\Sigma \,
                         \oplus \ldots \oplus \, \mathfrak{g}_r^\Sigma~.
        \end{equation}
        We then say that\/ $\mathfrak{m}$ is of\,
        \textbf{$\Sigma$-simple type}.
  \item $\mathfrak{m}'$ is the direct sum of all\/ $\Sigma$-simple
        ideals of\/~$\mathfrak{g}$ except two isomorphic ones, say\/
        $\mathfrak{g}_i^\Sigma$ and\/~$\mathfrak{g}_j^\Sigma$, with the
        ``diagonal'' subalgebra\/ $\mathfrak{g}_{ij}^\Sigma$ of\/~%
        $\, \mathfrak{g}_i^\Sigma \oplus \mathfrak{g}_j^\Sigma$,
        \begin{equation} \label{eq:MAXISA3}
         \mathfrak{m}'~=~\bigoplus_{\genfrac{}{}{0pt}{}{k=1}{k \neq i,j}}^r
                         \mathfrak{g}_k^\Sigma \, \oplus \,
                         \mathfrak{g}_{ij}^\Sigma~,
        \end{equation}

 \pagebreak

        where ``diagonal'' means that under suitable\/ $\Sigma$-equivariant
        isomorphisms $\, \mathfrak{g}_i^\Sigma \cong \mathfrak{g}_s^\Sigma \,$
        and $\, \mathfrak{g}_j^\Sigma \cong \mathfrak{g}_s^\Sigma$, the sub%
        algebra\/ $\mathfrak{g}_{ij}^\Sigma$ of\/~$\, \mathfrak{g}_i^\Sigma
        \oplus \mathfrak{g}_j^\Sigma \,$ corresponds to the subalgebra
        \[
         \mathrm{diag} \; \mathfrak{g}_s^\Sigma~
         =~\{ (X,X)~|~X \in \mathfrak{g}_s^\Sigma \}
        \]
        of~$\, \mathfrak{g}_s^\Sigma \oplus \mathfrak{g}_s^\Sigma$.
        We then say that\/ $\mathfrak{m}$ is of\,
        \textbf{$\Sigma$-diagonal type}.
        In this case, $\mathfrak{m}$ is a maximal\/
        $\Sigma$-invariant subalgebra of\/~$\mathfrak{g}$.
 \end{enumerate}
 Moreover, $\mathfrak{m}$ will be a maximal\/ $\Sigma$-invariant ideal
 of\/~$\mathfrak{g}$ if and only if\/ $\mathfrak{m}$ is of\, $\Sigma$-%
 simple type and $\; \mathfrak{m}_i^\Sigma = \{0\}$; in all other cases,
 $\mathfrak{m}$ is self-normalizing.
\end{theorem}

\begin{proof}
 Letting $\pi_i^\Sigma$ denote the projection of $\mathfrak{g}$
 onto $\mathfrak{g}_i^\Sigma$, which is an $\mathrm{Aut}_\Sigma^{}%
 (\mathfrak{g})$-equivariant Lie algebra homomorphism, and
 $\mathfrak{m}_i^\Sigma$ denote the image of~$\mathfrak{m}$
 under this projection, we note first of all that there can be
 at most one index, say $i$, for which $\mathfrak{m}_i^\Sigma$
 is properly contained in~$\mathfrak{g}_i^\Sigma$.
 Indeed, if there were two such indices, say $i$ and $j$, then
 \[
  \mathfrak{z} \; \oplus \, \mathfrak{g}_1^\Sigma \, \oplus \ldots \oplus \,
  \mathfrak{m}_i^\Sigma \, \oplus \ldots \oplus \,
  \mathfrak{g}_j^\Sigma \, \oplus \ldots \oplus \, \mathfrak{g}_r^\Sigma
 \]
 and
 \[
  \mathfrak{z} \; \oplus \, \mathfrak{g}_1^\Sigma \, \oplus \ldots \oplus \,
  \mathfrak{g}_i^\Sigma \, \oplus \ldots \oplus \,
  \mathfrak{m}_j^\Sigma \, \oplus \ldots \oplus \, \mathfrak{g}_r^\Sigma
 \]
 would be proper $\mathrm{Aut}_\Sigma^{}(\mathfrak{g},\mathfrak{m})$-%
 invariant subalgebras of~$\mathfrak{g}$ and
 \[
  \mathfrak{z} \; \oplus \, \mathfrak{g}_1^\Sigma \, \oplus \ldots \oplus \,
  \mathfrak{m}_i^\Sigma \, \oplus \ldots \oplus \,
  \mathfrak{m}_j^\Sigma \, \oplus \ldots \oplus \, \mathfrak{g}_r^\Sigma
 \]
 would be another $\mathrm{Aut}_\Sigma^{}(\mathfrak{g},\mathfrak{m})$-%
 invariant subalgebra of~$\mathfrak{g}$ properly contained in any of the
 two preceding ones and containing~$\mathfrak{m}$ (though perhaps not
 properly), which contradicts maximality of~$\mathfrak{m}$.
 Similarly, if there is one such index, say $i$, so that $\mathfrak{m}'$ is of
 the form given in equation~(\ref{eq:MAXISA2}), then $\mathfrak{m}_i^\Sigma$
 must be a $\Sigma$-primitive subalgebra of~$\mathfrak{g}_i^\Sigma$
 since otherwise, we could find a proper $\mathrm{Aut}_\Sigma^{}%
 (\mathfrak{g}_i^\Sigma,\mathfrak{m}_i^\Sigma)$-invariant subalgebra
 $\tilde{\mathfrak{m}}_i^\Sigma$ of~$\mathfrak{g}_i^\Sigma$ properly
 containing $\mathfrak{m}_i^\Sigma$, and then
 \[
  \mathfrak{z} \; \oplus \,
  \mathfrak{g}_1^\Sigma \, \oplus \ldots \oplus \,
  \mathfrak{g}_{i-1}^\Sigma \, \oplus \,
  \tilde{\mathfrak{m}}_i^\Sigma \, \oplus \,
  \mathfrak{g}_{i+1}^\Sigma \, \oplus \ldots \oplus \, \mathfrak{g}_r^\Sigma
 \]
 would be a proper $\mathrm{Aut}_\Sigma^{}(\mathfrak{g},\mathfrak{m})$-%
 invariant subalgebra of~$\mathfrak{g}$ properly containing~$\mathfrak{m}$,
 which contradicts maximality of~$\mathfrak{m}$.
 (Note that there are then two possibilities: either 
 $\, \mathfrak{m}_i^\Sigma = \{0\}$, which means that $\mathfrak{m}$
 is an ideal, or $\, \mathfrak{m}_i^\Sigma \neq \{0\}$, which implies
 that both $\mathfrak{m}$, as a subalgebra of~$\mathfrak{g}$, and
 $\mathfrak{m}_i^\Sigma$, as a subalgebra of~$\mathfrak{g}_i^\Sigma$,
 must be self-normalizing.)
 Finally, we must deal with the case that there is no such index, i.e.,
 that $\, \mathfrak{m}_i^\Sigma = \mathfrak{g}_i^\Sigma$ for all $i$.
 (This of course implies that $\mathfrak{m}$ cannot be an ideal
 and hence must be self-normalizing.
 Note also that this case cannot occur when $\, r=1$, so we
 may assume without loss of generality that $\, r \geqslant 2$.)
 Then for any fixed value of~$i$, $\mathfrak{m} + \mathfrak{g}_i^\Sigma$ is
 an $\mathrm{Aut}_\Sigma^{}(\mathfrak{g},\mathfrak{m})$-invariant subalgebra
 of~$\mathfrak{g}$ containing~$\mathfrak{m}$, which due to maximality of~%
 $\mathfrak{m}$ implies that either $\, \mathfrak{m} + \mathfrak{g}_i^\Sigma
 = \mathfrak{m}$, i.e., $\mathfrak{g}_i^\Sigma \subset \mathfrak{m}$, or
 $\, \mathfrak{m} + \mathfrak{g}_i^\Sigma = \mathfrak{g}$, and the second
 alternative must occur for at least two different values of~$i$ since
 if it did not occur at all, $\mathfrak{m}$ would have to be all of~%
 $\mathfrak{g}$, and if it occurred for only one value of~$i$,
 $\mathfrak{m}$ would contain the kernel of $\pi_i^\Sigma$ and
 hence, in order to satisfy the condition $\, \mathfrak{m}_i^\Sigma
 = \mathfrak{g}_i^\Sigma$, would once again have to be all of~$\mathfrak{g}$.
 \linebreak
 (As~we shall see soon, the alternative in question must in fact occur for
 exactly two different values of~$i$.)
 This in turn implies that for any fixed value of~$i$, $\mathfrak{m} \cap
 \mathfrak{g}_i^\Sigma$ is a $\Sigma$-invariant ideal of~$\mathfrak{g}$,
 since we may choose $\, j \neq i \,$  such that $\, \mathfrak{g}
 = \mathfrak{m} + \mathfrak{g}_j^\Sigma \,$ and then use the obvious
 inclusions $\, [\mathfrak{m} \cap \mathfrak{g}_i^\Sigma \,, \mathfrak{m}]
 \subset \mathfrak{m} \cap \mathfrak{g}_i^\Sigma \,$ and $\, [\mathfrak{m}
 \cap \mathfrak{g}_i^\Sigma \,, \mathfrak{g}_j^\Sigma] = \{0\} \,$ to
 derive that $\, [\mathfrak{m} \cap \mathfrak{g}_i^\Sigma \,, \mathfrak{g}]
 \subset \mathfrak{m} \cap \mathfrak{g}_i^\Sigma \,$.
 But $\mathfrak{g}_i^\Sigma$ is $\Sigma$-simple, so we must have either
 $\, \mathfrak{m} \cap \mathfrak{g}_i^\Sigma = \mathfrak{g}_i^\Sigma$,
 i.e., $\mathfrak{g}_i^\Sigma \subset \mathfrak{m}$, or $\, \mathfrak{m}
 \cap \mathfrak{g}_i^\Sigma = \{0\}$, i.e., $\mathfrak{g} = \mathfrak{m}
 \oplus \mathfrak{g}_i^\Sigma$.
 Now let $\mathfrak{g}''$ denote the $\Sigma$-invariant semisimple
 subalgebra of~$\mathfrak{g}$ obtained by taking the direct sum
 of all those $\Sigma$-simple ideals $\mathfrak{g}_i^\Sigma$ of~%
 $\mathfrak{g}$ that intersect $\mathfrak{m}$ trivially, and let
 $\mathfrak{m}''$ denote the intersection of $\mathfrak{m}$ and~%
 $\mathfrak{g}''$.
 Moreover, let $\, \mathfrak{g}'' = \mathfrak{g}_1'' \oplus
 \mathfrak{g}_2'' \,$ be any non-trivial direct decomposition
 of~$\mathfrak{g}''$ into two complementary $\Sigma$-invariant
 ideals of~$\mathfrak{g}$ (each of which must therefore be the
 direct sum of some of the $\Sigma$-simple ideals
 $\mathfrak{g}_i^\Sigma$ of~$\mathfrak{g}$ which intersect
 $\mathfrak{m}$ trivially), and let $\, \pi_1'': \mathfrak{g}''
 \rightarrow \mathfrak{g}_1'' \,$ and $\, \pi_2'': \mathfrak{g}''
 \rightarrow \mathfrak{g}_2'' \,$ be the corresponding projections.
 Then we claim that, for $k\!=\!1$ as well as for $k\!=\!2$, the
 restriction of~$\pi_k''$ to~$\mathfrak{m}''$ establishes a $\Sigma$-%
 equivariant Lie algebra isomorphism between $\mathfrak{m}''$ and
 $\mathfrak{g}_k''$.
 \begin{quote}
  To prove this for $k\!=\!1$, say, note that for any $\Sigma$-simple ideal
  $\mathfrak{g}_j^\Sigma$ of~$\mathfrak{g}$ contained in $\mathfrak{g}_2''$,
  we have $\, \mathfrak{g} = \mathfrak{m} \oplus \mathfrak{g}_j^\Sigma \,$
  and hence $\, \mathfrak{g}'' = \mathfrak{m}'' \oplus \mathfrak{g}_j^\Sigma$.
  Therefore, we may decompose any element $X''$ of~$\mathfrak{g}''$ in the
  form $\, X'' = X_{\mathfrak{m}}'' + X_j'' \,$ with $\, X_{\mathfrak{m}}''
  \in \mathfrak{m}''$ and $X_j'' \in \mathfrak{g}_j^\Sigma \,$
  to get $\, \pi_1''(X'') = \pi_1''(X_{\mathfrak{m}}'')$, which shows that
  $\pi_1''$ maps $\mathfrak{m}''$ onto $\mathfrak{g}_1''$. On the other hand,
  we can repeat the argument employed above to show that $\, \mathfrak{m}''
  \cap \mathfrak{g}_2'' = \mathfrak{m} \cap \mathfrak{g}_2''$ \linebreak
  is a $\Sigma$-invariant ideal of~$\mathfrak{g}$ (choose any 
  $\Sigma$-simple ideal $\mathfrak{g}_i^\Sigma$ of~$\mathfrak{g}$
  contained in $\mathfrak{g}_1''$ such that $\, \mathfrak{g}
  = \mathfrak{m} \oplus \mathfrak{g}_i^\Sigma \,$ and use
  the obvious inclusions \linebreak $\, [\mathfrak{m} \cap \mathfrak{g}_2''
  \,,\mathfrak{m}] \subset \mathfrak{m} \cap \mathfrak{g}_2'' \,$ and
  $\, [\mathfrak{m} \cap  \mathfrak{g}_2'' \,, \mathfrak{g}_i^\Sigma]
  = \{0\} \,$ to derive that $\, [\mathfrak{m} \cap \mathfrak{g}_2'' \,,
  \mathfrak{g}] \subset \mathfrak{m} \cap \mathfrak{g}_2''$),
  so it must be the direct sum of some of the $\Sigma$-simple
  ideals $\mathfrak{g}_j^\Sigma$ of~$\mathfrak{g}$ contained
  in~$\mathfrak{g}_2''$, which is impossible since by the definition
  of~$\mathfrak{g}''$, any such $\mathfrak{g}_j^\Sigma$ intersects
  $\mathfrak{m}$ trivially. This shows that $\, \mathfrak{m}'' \cap
  \mathfrak{g}_2'' = \{0\}$, and since $\mathfrak{g}_2''$ is the kernel
  of~$\pi_1''$, that the restriction of~$\pi_1''$ to~$\mathfrak{m}''$ is
  one-to-one.
 \end{quote}
 But this can only happen if both $\mathfrak{g}_1''$ and~$\mathfrak{g}_2''$
 are $\Sigma$-simple and isomorphic to $\mathfrak{m}''$, which completes
 the proof of the main statement of the theorem. \\[1mm]
 What remains to be shown is that $\, \mathrm{diag} \;
 \mathfrak{g}_s^\Sigma \,$ is really a maximal $\Sigma$-invariant
 subalgebra of~$\, \mathfrak{g}_s^\Sigma \oplus \mathfrak{g}_s^\Sigma$.
 To this end, assume that $\mathfrak{h}$ is a $\Sigma$-invariant sub%
 algebra of~$\, \mathfrak{g}_s^\Sigma \oplus \mathfrak{g}_s^\Sigma \,$
 containing $\, \mathrm{diag} \; \mathfrak{g}_s^\Sigma \,$ and define
 $\, \mathfrak{h}_1^{} = \mathrm{pr}_1^{} \bigl( \mathfrak{h} \,\cap
 (\mathfrak{g}_s^\Sigma \oplus \{0\}) \bigr) \,$ and $\, \mathfrak{h}_2^{}
 = \mathrm{pr}_2^{} \bigl( \mathfrak{h} \,\cap (\{0\} \oplus
 \mathfrak{g}_s^\Sigma) \bigr) \,$ where $\mathrm{pr}_1^{}$
 and $\mathrm{pr}_2^{}$ denote the first and second projection
 from~$\, \mathfrak{g}_s^\Sigma \oplus \mathfrak{g}_s^\Sigma \,$
 to~$\, \mathfrak{g}_s^\Sigma$, respectively.
 Obviously, $\mathfrak{h}_1^{}$ and $\mathfrak{h}_2^{}$ are
 $\Sigma$-invariant subalgebras of $\mathfrak{g}_s^\Sigma$,
 but they are even ideals: this follows by observing that since
 $\mathfrak{h}$ contains $\, \mathrm{diag} \; \mathfrak{g}_s^\Sigma$,
 the projections $\mathrm{pr}_1^{}$ and $\mathrm{pr}_2^{}$ both map
 $\mathfrak{h}$ onto~$\mathfrak{g}_s^\Sigma$, and hence
 \begin{eqnarray*}
  X_1^{} \in \mathfrak{g}_s^\Sigma~,~Y_1^{} \in \mathfrak{h}_1^{}
  &\Longrightarrow& \mbox{$\exists$ $X_2^{} \in \mathfrak{g}_s^\Sigma$
                          such that } (X_1^{},X_2^{}) \in \mathfrak{h}~,~
                    (Y_1^{},0) \in \mathfrak{h} \\
  &\Longrightarrow& ([X_1^{},Y_1^{}],0) = [(X_1^{},X_2^{}),(Y_1^{},0)]
                    \in \mathfrak{h} \\
  &\Longrightarrow& [X_1^{},Y_1^{}] \in \mathfrak{h}_1^{} \\[2mm]
  X_2^{} \in \mathfrak{g}_s^\Sigma~,~Y_2^{} \in \mathfrak{h}_2^{}
  &\Longrightarrow& \mbox{$\exists$ $X_1^{} \in \mathfrak{g}_s^\Sigma$
                          such that } (X_1^{},X_2^{}) \in \mathfrak{h}~,~
                    (0,Y_2^{}) \in \mathfrak{h} \\
  &\Longrightarrow& (0,[X_2^{},Y_2^{}]) = [(X_1^{},X_2^{}),(0,Y_2^{})]
                    \in \mathfrak{h} \\
  &\Longrightarrow& [X_2^{},Y_2^{}] \in \mathfrak{h}_2^{}
 \end{eqnarray*}
 But $\mathfrak{g}_s^\Sigma$ is $\Sigma$-simple, so it follows that either
 $\, \mathfrak{h}_1^{} = \mathfrak{g}_s^\Sigma \,$ or $\, \mathfrak{h}_2^{}
 = \mathfrak{g}_s^\Sigma \,$ or $\, \mathfrak{h}_1^{} = \{0\}$ \linebreak
 and $\, \mathfrak{h}_2^{} = \{0\}$.
 In the first two cases, we can use the hypothesis that $\mathfrak{h}$
 contains $\, \mathrm{diag} \; \mathfrak{g}_s^\Sigma \,$ to conclude that
 $\, \mathfrak{h} = \mathfrak{g}_s^\Sigma \oplus \mathfrak{g}_s^\Sigma$,
 whereas in the third case, it follows that $\, \mathfrak{h} =
 \mathrm{diag} \; \mathfrak{g}_s^\Sigma$.
\end{proof}
\begin{theorem} \label{theo:maxisa2}
 Let\/ $\mathfrak{g}$ be a Lie algebra and\/ $\Sigma$ be a subgroup
 of its outer auto\-morphism group\/ $\mathrm{Out}(\mathfrak{g})$.
 Assuming that\/ $\mathfrak{g}$ is\/ $\Sigma$-simple, with canonical
 decomposition
 \begin{equation} \label{eq:DECSLA1}
  \mathfrak{g}~
  =~\mathfrak{g}_s^{} \, \oplus \ldots \oplus \, \mathfrak{g}_s^{}
  \quad \mbox{($n$ summands)}
 \end{equation}
 into the direct sum of\/ $n$ copies of the same simple Lie algebra\/~%
 $\mathfrak{g}_s^{}$, let\/ $\mathfrak{m}$ be a\/ $\Sigma$-primitive
 subalgebra of\/~$\mathfrak{g}$.
 Then one of the following two alternatives holds.
 \begin{enumerate}
  \item Up to ``twisting'' with an appropriate automorphism of\/~%
        $\mathfrak{g}$, $\mathfrak{m}$ is the direct sum of\/ $n$
        copies of the same\/ $\Sigma_s^{}$-primitive subalgebra\/
        $\mathfrak{m}_s^{}$ of\/~$\mathfrak{g}_s^{}$,
        \begin{equation} \label{eq:MAXISA4}
         \mathfrak{m}~
         =~\mathfrak{m}_s^{} \, \oplus \ldots \oplus \, \mathfrak{m}_s^{}
         \quad \mbox{($n$ summands)}~,
        \end{equation}
        where the subgroup\/ $\Sigma_s^{}$ of\/~$\mathrm{Out}%
        (\mathfrak{g}_s^{})$ is obtained from the subgroup\/
        $\Sigma$ of\/~$\mathrm{Out}(\mathfrak{g})$ by projection.
        We then say that\/ $\mathfrak{m}$ is of\, \textbf{simple type}.
  \item Up to ``twisting'' with an appropriate automorphism of\/~%
        $\mathfrak{g}$, $\mathfrak{m}$ is the direct sum of a certain
        number (say\/ $q$) of copies of the diagonal subalgebra of
        the direct sum of a certain number (say\/ $p$) of copies
        of\/~$\mathfrak{g}_s^{}$,
        \begin{equation} \label{eq:MAXISA5}
         \begin{array}{cc}
          \mathfrak{m}~=~\mathrm{diag}_p \; \mathfrak{g}_s^{} \,
                         \oplus \ldots \oplus \,
                         \mathrm{diag}_p \; \mathfrak{g}_s^{}
          \quad \mbox{($q$ summands)} \\[2mm]
          \mathrm{diag}_p \; \mathfrak{g}_s^{}~
          =~\{ (X_s^{},\ldots,X_s^{})~|~X_s^{} \in \mathfrak{g}_s^{} \}
          \quad \mbox{($p$ summands)}
         \end{array}~,
        \end{equation}
        where\/ $p$ and\/ $q$ are divisors of~$n$, with $\, p>1 \,$ and
        $\, q<n$, chosen such that\/ $p$ is the minimum and\/ $q$ the
        maximum possible value for which the resulting partition
        \begin{equation} \label{eq:MAXISA6}
        \{1,\ldots,n\}~=~\big\{ \{1,\ldots,p\} , \ldots ,
                                \{n-p+1,\ldots,n\} \bigr\}
        \end{equation}
        is\/ $\Sigma$-invariant.
        We then say that\/ $\mathfrak{m}$ is of\, \textbf{diagonal type}.
 \end{enumerate}
\end{theorem}
Before turning to the proof, let us comment on the structure of the
automorphism group of a semisimple Lie algebra $\mathfrak{g}$ which,
as in the above theorem, is the direct sum of a certain number of
copies (say $n$) of one and the same simple Lie algebra~$\mathfrak{g}_s$:
this information will be needed for properly understanding some of the
statements of the theorem.
In this situation, the automorphism group $\mathrm{Aut}(\mathfrak{g})$
of~$\mathfrak{g}$ is the \emph{wreath product} of the automorphism group
$\mathrm{Aut}(\mathfrak{g}_s)$ of~$\mathfrak{g}_s$ with the
permutation group $S_n$, that is, the semidirect product
$\, \mathrm{Aut}(\mathfrak{g}_s)^{\,n} : S_n \,$,
and since inner automorphisms preserve the summands in the
direct sum~(\ref{eq:DECSLA1}), so $\, \mathrm{Inn}(\mathfrak{g})
= \mathrm{Inn}(\mathfrak{g}_s)^{\,n}$, we conclude that similarly,
the outer auto\-morphism group $\mathrm{Out}(\mathfrak{g})$ of~%
$\mathfrak{g}$ is the \emph{wreath product} of the outer automorphism
group $\mathrm{Out}(\mathfrak{g}_s)$ of~$\mathfrak{g}_s$ with the
permutation group $S_n$, that is, the semidirect product
$\, \mathrm{Out}(\mathfrak{g}_s)^{\,n} : S_n \,$.
Note that the projection from this semidirect product to~$S_n$ is
a homomorphism while that to~$\mathrm{Out}(\mathfrak{g}_s)^{\,n}$
is not, nor is its composition with any one of the projections
to~$\mathrm{Out}(\mathfrak{g}_s)$.
Remarkably, however, its composition with the $i$-th projection to~%
$\mathrm{Out}(\mathfrak{g}_s)$ does become a homomorphism when restricted
to the semidirect pro\-duct $\, \mathrm{Out}(\mathfrak{g}_s)^{\,n} :
 (S_n)_i \,$, where $(S_n)_i$ is the stability group of $i$ in~$S_n$:
\begin{equation} \label{eq:STABGR}
 (S_n)_i~=~\{ \, \pi \in S_n~|~\pi(i) = i \, \}~.
\end{equation}
Next, let $\Sigma$ be a given subgroup of~$\mathrm{Out}(\mathfrak{g})$,
and denote by $S_n^\Sigma$ the subgroup of $S_n^{}$ obtained as the
image of~$\Sigma$ under the projection from~$\mathrm{Out}(\mathfrak{g})$
to~$S_n$.
Thus elements $\sigma$ of~$\Sigma$ are certain $(n+1)$-tuples $(\sigma_1,
\ldots,\sigma_n;\pi)$, where $\,\sigma_1,\ldots,\sigma_n \in \mathrm{Out}%
(\mathfrak{g}_s)$ and $\, \pi \in S_n \,$: they can be represented
by $(n+1)$-tuples $(\varphi_1,\ldots,\varphi_n;\pi)$, where $\, \varphi_1,
\ldots,\varphi_n \in \mathrm{Aut}(\mathfrak{g}_s)$ and $\, \pi \in S_n \,$
(each $\varphi_i$ being defined up to multiplication by inner automorphisms
of~$\mathfrak{g}_s$), acting on $n$-tuples $(X_1,\ldots,X_n)$ of elements
$X_1,\ldots,X_n$ of~$\mathfrak{g}_s$ according to
\begin{equation} \label{eq:WREATH1}
 (\varphi_1,\ldots,\varphi_n;\pi) \cdot (X_1,\ldots,X_n)~
 =~(\varphi_1(X_{\pi(1)}),\ldots,\varphi_n(X_{\pi(n)}))~.
\end{equation}
As noticed above, the projection taking $\sigma$ to $\pi$ is a homo%
morphism whereas the ones taking $\sigma$ to $(\sigma_1,\ldots,\sigma_n)$
or to $\sigma_i$ are not, but the latter does become a homomorphism when
restricted to the intersection of~$\Sigma$ with the semi\-direct product
$\, \mathrm{Out}(\mathfrak{g}_s)^{\,n} : (S_n)_i$.
Its image is a subgroup of $\mathrm{Out}(\mathfrak{g}_s)$ that we shall
denote by~$\Sigma_{s,i}$ and whose definition can be described explicitly
as follows: an element $\sigma_s$ of~$\mathrm{Out}(\mathfrak{g}_s)$
belongs to~$\Sigma_{s,i}$ if and only if there exist elements
$\sigma_1,\ldots,\sigma_n$ of~$\mathrm{Out}(\mathfrak{g}_s)$ such
that $\, \sigma_i = \sigma_s \,$ and a permutation $\, \pi \in S_n \,$
satis\-fying $\pi(i) = i$ such that, when each $\sigma_k$ is interpreted
as representing an (equi\-valence class of) isomorphism(s) from~%
$\mathfrak{g}_{s,\pi(k)}$ to~$\mathfrak{g}_{s,k}$, $(\sigma_1,
\ldots,\sigma_n;\pi)$ belongs to~$\Sigma$. \linebreak
(Here and in what follows, we shall, for $\, 1 \leqslant
i \leqslant n$, refer to the ideal
\[
 \mathfrak{g}_{s,i} = \{0\} \oplus \ldots \oplus \mathfrak{g}_s
                            \oplus \ldots \oplus \{0\}
\]
of~$\mathfrak{g}$, with $\mathfrak{g}_s$ in the $i$-th position, as its
$i$-th simple summand.)
Substantial further simplifications can be achieved by using the freedom of
modifying the identification~(\ref{eq:DECSLA1}) above, ``twisting'' it with
an appropriate automorphism of~$\mathfrak{g}$, which amounts to replacing the
subgroup $\Sigma$ of~$\mathrm{Out}(\mathfrak{g})$ with a conjugate subgroup.
For~example, we can perform such a ``twist'' with an automorphism preserving
each of the simple summands $\mathfrak{g}_{s,i}$ of~$\mathfrak{g}$
(i.e., belonging to the subgroup $\mathrm{Aut}(\mathfrak{g}_s)^{\,n}$
of~$\mathrm{Aut}(\mathfrak{g})$), in order to guarantee that,
for $\, 1 \leqslant i,j \leqslant n$, there exists an element
$\, \sigma_{i,j} \in \Sigma \,$ which can be represented by an auto%
morphism $\, \phi_{i,j} \in \mathrm{Aut}_\Sigma^{}(\mathfrak{g}) \,$
taking the $j$-th simple summand $\mathfrak{g}_{s,j}$ to the $i$-th
simple summand $\mathfrak{g}_{s,i}$ and acting as the identity
on~$\mathfrak{g}_s$, i.e.,%
\footnote{Unfortunately, there is no information about how $\phi_{i,j}$
acts on the remaining simple summands, which is why we use quotation
marks when referring to these automorphisms as ``shift operators''.}
\[
 \phi_{i,j} \cdot
 (\underbrace{0,\ldots,0}_{j-1},X,\underbrace{0,\ldots,0}_{n-j})~
 =~(\underbrace{0,\ldots,0}_{i-1},X,\underbrace{0,\ldots,0}_{n-i})~.
\]
Note that once this is achieved, the fact that $\Sigma_{s,i}$ and
$\Sigma_{s,j}$ are conjugate under $\phi_{i,j}$ implies that the
subgroups $\Sigma_{s,1},\ldots,\Sigma_{s,n}$ of~$\mathrm{Out}%
(\mathfrak{g}_s)$ will all be equal, so we may drop the index~$i$
and simply denote them by $\Sigma_s$.
Similarly, given any $\Sigma$-invariant subalgebra $\mathfrak{h}$
of~$\mathfrak{g}$, its image under the $i$-th projection from~%
$\mathfrak{g}$ to~$\mathfrak{g}_s$ will not depend on~$i$ and
will be a $\Sigma_s$-invariant subalgebra of~$\mathfrak{g}_s$
which we may simply denote by~$\mathfrak{h}_s$.
Moreover, we obtain the inclusion
\begin{equation} \label{eq:WREATH2}
 \mathrm{Aut}_\Sigma^{}(\mathfrak{g})~\subset~
 \mathrm{Aut}_{\Sigma_s}^{}(\mathfrak{g}_s^{})^n : S_n^\Sigma
\end{equation}
which implies the inclusion
\begin{equation} \label{eq:WREATH3}
 \Sigma~\subset~\Sigma_s^n : S_n^\Sigma
\end{equation}
and, for any $\Sigma$-invariant subalgebra $\mathfrak{h}$ of~$\mathfrak{g}$,
the inclusion
\begin{equation} \label{eq:WREATH4}
 \mathrm{Aut}_\Sigma^{}(\mathfrak{g},\mathfrak{h})~\subset~
 \mathrm{Aut}_{\Sigma_s}^{}(\mathfrak{g}_s^{},\mathfrak{h}_s^{})^n : S_n^\Sigma
\end{equation}
Indeed, to see that given an element $(\varphi_1,\ldots,\varphi_n;\pi)$
of~$\mathrm{Aut}_\Sigma^{}(\mathfrak{g})$, its $j$-th component $\varphi_j$
really belongs to $\mathrm{Aut}_{\Sigma_s}(\mathfrak{g}_s)$, we use that it
is also the $j$-th component of the product $\, (\varphi_1,\ldots,\varphi_n;
\pi) \, \phi_{\pi(j),j}$, which by definition maps to~$(S_n)_j$ under the
last projection and still belongs to~$\mathrm{Aut}_\Sigma^{}(\mathfrak{g})$,
and it is then clear that if $\mathfrak{h}$ is $\Sigma$-invariant,
$(\varphi_1,\ldots,\varphi_n;\pi)$ and $\phi_{\pi(j),j}$ both preserve
$\mathfrak{h}$ and hence $\varphi_j$ preserves $\mathfrak{h}_s$.
\begin{quote}
 The construction of such a ``twisting automorphism'' is elementary,
 using the hypothesis that $\mathfrak{g}$ is $\Sigma$-simple, i.e.,
 that $\Sigma$ permutes the simple summands of~$\mathfrak{g}$
 transitively: it implies that for $\, 1 \leqslant j \leqslant n$,
 there exists an element $\, \sigma_j \in \Sigma \,$ which can be
 represented by an automorphism $\, \varphi_j \in
 \mathrm{Aut}_\Sigma^{}(\mathfrak{g}) \,$ taking $\mathfrak{g}_{s,j}$
 to~$\mathfrak{g}_{s,1}$, say, so it will act on~$\mathfrak{g}_{s,j}$ by
 shifting elements of $\mathfrak{g}_s$ from the $j$-th to the first
 position and then applying some automorphism $(\varphi_j)_1$
 of~$\mathfrak{g}_s$. Of course, we have $\, \sigma_1 = 1 \,$ and hence
 $(\varphi_1)_1$ can be taken as the identity of~$\mathfrak{g}_s$.
 Then the subgroup $\, \Sigma' = \sigma \Sigma \sigma^{-1} \,$ of~%
 $\mathrm{Out}(\mathfrak{g})$, where $\, \sigma \in \mathrm{Out}%
 (\mathfrak{g}) \,$ is the class of $\, \varphi = (\varphi_1)_1
 \times \ldots \times (\varphi_n)_1 \in \mathrm{Aut}(\mathfrak{g})$,
 has the required property, with $\, \sigma_{1,j}' = \sigma
 \sigma_j^{} \sigma^{-1} \,$ and $\, \sigma_{i,j}' = \sigma
 \sigma_i^{-1} \sigma_j^{} \sigma^{-1}$.
\end{quote}
Similarly, given any $\Sigma$-invariant equivalence relation in the
set $\{1,\ldots,n\}$, we can perform a further ``twisting'' with an
automorphism which is a ``pure permutation'' of~$\mathfrak{g}$ (i.e.,
belongs to the subgroup $S_n$ of~$\mathrm{Aut}(\mathfrak{g})$) such
that the corresponding partition into equivalence classes takes the
``block'' form given by equation~(\ref{eq:MAXISA6}) with $\, 1 \leqslant
p,q \leqslant n \,$ satisfying $\, pq = n$; then $S_n^\Sigma$ becomes
a subgroup of the wreath product of the permutation group $S_p^{}$
with the permutation group $S_q^{}$, that is, of the semidirect product
$\, S_p^{\,q} : S_q^{} \,$.
\begin{quote}
 Again, the construction of such a ``twisting automorphism'' is elementary,
 using the hypothesis that $\mathfrak{g}$ is $\Sigma$-simple, i.e., that
 $\Sigma$ permutes the simple summands of~$\mathfrak{g}$ transitively:
 it implies that all equivalence classes must have the same cardinality,
 say $p$, which must therefore divide $n$.
\end{quote}
With these preliminaries out of the way, we can proceed to the
proof of \linebreak Theorem~\ref{theo:maxisa2}.
\begin{proof}
 To begin with, we observe that, once the aforementioned simplifications
 by means of appropriate ``twistings'' have been performed, it follows
 from the inclusion~(\ref{eq:WREATH4}) that if $\mathfrak{h}$ is any
 $\Sigma$-invariant subalgebra of~$\mathfrak{g}$ and $\mathfrak{h}_s$ is
 the $\Sigma_s$-invariant subalgebra of~$\mathfrak{g}_s$ obtained by any
 one of the projections from~$\mathfrak{g}$ to~$\mathfrak{g}_s$, then the
 direct sum $\, \mathfrak{h}_s \oplus \ldots \oplus \mathfrak{h}_s \,$
 will be a subalgebra of~$\mathfrak{g}$ which is not only again $\Sigma$-%
 invariant but is also $\mathrm{Aut}_\Sigma^{}(\mathfrak{g},\mathfrak{h})$-%
 invariant.
 (The first part of this statement can also be inferred directly by noting
 that, essentially, invariance under $\Sigma$ is equi\-valent to invariance
 under $\Sigma_s$ together with invariance under the ``shift operators''
 introduced above, and for the direct sum $\, \mathfrak{h}_s \oplus \ldots
 \oplus \mathfrak{h}_s \,$, the latter is manifest.)
 Thus if $\mathfrak{m}$ is a $\Sigma$-primitive subalgebra of~$\mathfrak{g}$,
 there are precisely two distinct possibilities: either $\mathfrak{m}_s$ is
 a proper subalgebra of~$\mathfrak{g}_s$, $\mathfrak{m} = \mathfrak{m}_s
 \oplus \ldots \oplus \mathfrak{m}_s \,$ and, as is then easy to see,
 $\mathfrak{m}_s$ is a $\Sigma_s$-primitive subalgebra of~$\mathfrak{g}_s$,
 or else $\, \mathfrak{m}_s = \mathfrak{g}_s$.
 To handle the second case, consider, for $\, 1 \leqslant i,j \leqslant n$,
 the image $\, \mathrm{pr}_j \bigl( \mathfrak{m} \cap \ker \, \mathrm{pr}_i
 \bigr) \,$ of the intersection of~$\mathfrak{m}$ with the kernel of the
 $i$-th projection $\mathrm{pr}_i$ from~$\mathfrak{g}$ to~$\mathfrak{g}_s$
 under the $j$-th projection $\mathrm{pr}_j$ from~$\mathfrak{g}$ to~%
 $\mathfrak{g}_s$ and note that this is an ideal of~$\mathfrak{g}_s$.
 \begin{quote}
  Indeed, suppose that $\, X \in \mathrm{pr}_j \bigl( \mathfrak{m} \cap
  \ker \, \mathrm{pr}_i \bigr) \,$ and $\, X' \in \mathfrak{g}_s \,$.
  This means that there exists $\, (X_1,\ldots,X_n) \in \mathfrak{m} \,$
  such that $\, X_j = X \,$ and $\, X_i = 0$, and since $\mathrm{pr}_j$
  maps $\mathfrak{m}$ onto $\mathfrak{g}_s$, there also exists
  $\, (X_1',\ldots,X_n') \in \mathfrak{m}$ \linebreak such that
  $\, X_j' = X'$.
  Since $\mathfrak{m}$ is a subalgebra and $\ker \, \mathrm{pr}_i$ is
  an ideal of~$\mathfrak{g}$, it follows that
  \[
   ([X_1',X_1],\ldots,[X_n',X_n])~=~[(X_1',\ldots,X_n'),(X_1,\ldots,X_n)]
  \]
  belongs to $\, \mathfrak{m} \cap \ker \, \mathrm{pr}_i \,$ and
  therefore, $\, [X',X] \in \mathrm{pr}_j \bigl( \mathfrak{m} \cap
  \ker \, \mathrm{pr}_i \bigr)$.
 \end{quote}
 Since $\mathfrak{g}_s$ is simple, there are only two possibilites: either
 $\; \mathrm{pr}_j \bigl( \mathfrak{m} \cap \ker \, \mathrm{pr}_i \bigr)
 = \{0\}$, which means that $\, \mathfrak{m} \cap \ker \, \mathrm{pr}_i
 \subset \mathfrak{m} \cap \ker \, \mathrm{pr}_j \,$, or else
 $\; \mathrm{pr}_j \bigl( \mathfrak{m} \cap \ker \, \mathrm{pr}_i \bigr)
 = \mathfrak{g}_s \,$.
 Using that not only $\mathrm{pr}_j$ but also $\mathrm{pr}_i$ maps
 $\mathfrak{m}$ onto~$\mathfrak{g}_s$ and that $\mathfrak{m}$ is a
 subalgebra of~$\mathfrak{g}^{}$, we can sharpen the second alternative
 to the statement that for any two elements of~$\mathfrak{g}_s$, there
 exists $\, (X_1,\ldots,X_n) \in \mathfrak{m} \,$ such that $X_i$
 equals the first and $X_j$ equals the second.
 As a result, it becomes evident that both alternatives are symmetric
 under the exchange of $i$ and~$j$ and therefore, the first alternative
 provides an equivalence relation $\sim$ in $\{1,\ldots,n\}$, defined by
 \[
  i \sim j~~\Longleftrightarrow~~
  \mathfrak{m} \cap \ker \, \mathrm{pr}_i~
  =~\mathfrak{m} \cap \ker \, \mathrm{pr}_j~.
 \]
 When $i$ and $j$ are in the same equivalence class, one projection factors
 over the kernel of the other to provide mutually inverse automorphisms
 $\varphi_{i,j}$ and $\varphi_{j,i}$ of~$\mathfrak{g}_s$ such that
 $\; \varphi_{ji} \circ \mathrm{pr}_i|_{\mathfrak{m}} = \mathrm{pr}_j|%
 _{\mathfrak{m}} \;$ and $\; \varphi_{ij} \circ \mathrm{pr}_j|_{\mathfrak{m}}
 = \mathrm{pr}_i|_{\mathfrak{m}} \,$.
 Together with an appropriate permutation to~bring the resulting partition
 of~$\{1,\ldots,n\}$ into the form~(\ref{eq:MAXISG6}), these can be used
 to define the ``twist'' automorphism that maps $\mathfrak{m}$ into the
 subalgebra given by equation~(\ref{eq:MAXISG5}).
 The requirement that $\mathfrak{m}$ should be $\Sigma$-primitive
 (which incudes the condition that it should not be equal to all
 of~$\mathfrak{g}$) is then guaranteed by the condition that $p$
 should be chosen as small as possible and $q$ as large as possible,
 with the restriction that $\, p>1$, $q<n$.
\end{proof}

In order to perform the same reduction at the Lie group level, we proceed
in several steps, assuming as before that $G$ is a Lie group with connected
one-component~$G_0$, component group~$\Gamma$ and Lie algebra $\mathfrak{g}$
(all steps refer to properties of~$G_0$):
\begin{itemize}
 \item \emph{Step 1:}\/ from general reductive to centerfree
       semisimple Lie groups,
 \item \emph{Step 2:}\/ from centerfree semisimple to centerfree
       $\Gamma$-simple Lie groups,
 \item \emph{Step 3:}\/ from centerfree $\Gamma$-simple to
       centerfree simple Lie groups.
\end{itemize}

The first step is elementary: it simply consists in dividing out the
center $Z(G_0)$ of the connected one-component $G_0$ of~$G$, which
is a $\Gamma$-invariant closed normal subgroup of~$G_0$ and which,
according to Theorem~\ref{theo:maxsg1}, is contained in any maximal
$\Gamma$-invariant subgroup of~$G_0$.
Thus considering the quotient Lie group $\, \hat{G} = G/Z(G_0)$,
with connected one-component $\, \hat{G}_0 = G_0/Z(G_0)$, component
group $\, \hat{\Gamma} = \Gamma \,$ and Lie algebra $\, \mathfrak{g}/
\mathfrak{z} \,$, it becomes obvious that every maximal subgroup $M$ of~$G$
is the inverse image of a maximal subgroup $\hat{M}$ of~$\hat{G}$ (namely
$\, \hat{M} = M/Z(G_0)$) under the projection from~$G$ to~$\hat{G}$ and
simi\-larly every maximal $\Gamma$-invariant subgroup $M_1$ of~$G_0$ is
the inverse image of a maximal $\Gamma$-invariant subgroup $\hat{M}_1$
of~$\hat{G}_0$ (namely $\, \hat{M}_1 = M_1/Z(G_0)$) under the projection
from~$G_0$ to~$\hat{G}_0$: in fact, this prescription establishes a one-%
to-one correspondence between maxi\-mal subgroups $M$ of~$G$ and maximal
subgroups $\hat{M}$ of~$\hat{G}$ and similarly between maxi\-mal $\Gamma$-%
invariant subgroups $M_1$ of~$G_0$ and maximal $\Gamma$-invariant subgroups
$\hat{M}_1$ of~$\hat{G}_0$.
Note that this procedure is completely general (it works for any Lie group),
so we must only convince ourselves that if $G$ is reductive, $\hat{G}$ will
be semisimple and will have trivial center.
To this end, assume that $\mathfrak{g}$ is reductive and consider the
derived subgroup $G_0^{\,\prime}$ of~$G_0$, which is semisimple and whose
center $Z(G_0^{\,\prime})$ is equal to its intersection $\, Z(G_0) \cap
G_0^{\,\prime} \,$ with the center $Z(G_0)$ of~$G_0$.
(The only non-trivial statement in this equality is the inclusion
$\, Z(G_0^{\,\prime}) \subset Z(G_0)$, which follows by observing
that $\, g_0 \in Z(G_0) \,$ is equivalent to $\, \mathrm{Ad}(g_0)
= \mathrm{id}_{\mathfrak{g}} \,$ while $\, g_0 \in Z(G_0^{\,\prime}) \,$
is equivalent to $\, \mathrm{Ad}(g_0)|_{\mathfrak{g}'} = \mathrm{id}_%
{\mathfrak{g}'}$, since $G_0$ and $G_0^{\,\prime}$ are connected, but
$\, \mathfrak{g} = \mathfrak{z} \oplus \mathfrak{g}' \,$ and obviously
$\, \mathrm{Ad}(g_0)|_{\mathfrak{z}} = \mathrm{id}_{\mathfrak{z}} \,$
for any $\, g_0 \in G_0$, since $G_0$ is connected.)
Therefore, $\hat{G}_0$ can be identified with the quotient group
$G_0^{\,\prime}/Z(G_0^{\,\prime})$, implying that $\hat{G}_0$ has
trivial center (see Proposition~6.30 of Ref.~\cite{Kn}): it
is the adjoint group of~$G_0^{\,\prime}$.

The second and third step will be summarized in the form of two
theorems which result from transferring Theorems~\ref{theo:maxisa1}
and~\ref{theo:maxisa2} from the Lie algebra to the Lie group context:
this can be done, e.g., by invoking Theorems~\ref{theo:maxsg1} and~%
\ref{theo:maxsg2}. Direct proofs can also be given but will be left
to the reader since they are completely analogous to the proofs of
Theorems~\ref{theo:maxisa1} and~\ref{theo:maxisa2} given above.
\begin{theorem} \label{theo:maxsg3}
 Let\/ $G$ be a Lie group with connected one-component\/~$G_0$ and component
 group\/~$\Gamma$, and let $\, \Sigma \subset \mathrm{Out}(G_0) \,$ be the
 image of\/~$\Gamma$ under the homomorphism~(\ref{eq:GREXT4}).
 Assuming that\/ $G_0$ is semisimple with trivial center and that
 \begin{equation} \label{eq:DECRLG1}
  G_0^{}~=~G_{1,0}^\Sigma \times \ldots \times G_{r,0}^\Sigma
 \end{equation}
 is its canonical decomposition into its\/ $\Sigma$-simple factors\/
 $G_{1,0}^\Sigma,\ldots,G_{r,0}^\Sigma$, let\/ $M$ be a maxi\-mal subgroup
 of\/~$G$ with connected one-component\/~$M_0$ such that\/ $M$ meets
 every connected component of\/~$G$.
 Set $\, M_1 = M \cap G_0$, so $\, M_0 \subset M_1 \subset M$.
 Then one of the following two alternatives holds:
 \begin{itemize}
  \item $M_1^{}$ is the direct product of all\/ $\Sigma$-simple factors of\/~%
        $G_0^{}$ except one, say\/~$G_{i,0}^\Sigma$, with a maximal\/ $\Sigma$-%
        invariant subgroup\/ $M_{i,1}^\Sigma$ of\/~$G_{i,0}^\Sigma \,$:
        \begin{equation} \label{eq:MAXISG1}
         M_1^{}~=~G_{1,0}^\Sigma \times \ldots \times G_{i-1,0}^\Sigma
                  \times M_{i,1}^\Sigma \times
                  G_{i+1,0}^\Sigma \times \ldots \times G_{r,0}^\Sigma
        \end{equation}
        We then say that\/ $M$ and\/ $M_1^{}$ are of\,
        \textbf{$\Sigma$-simple type}.
  \item $M_1^{}$ is the direct product of all\/ $\Sigma$-simple factors
        of\/~$G_0^{}$ except two isomorphic ones, say\/ $G_{i,0}^\Sigma$
        and\/~$G_{j,0}^\Sigma$, with the ``diagonal'' subgroup\/
        $G_{ij,0}^\Sigma$ of $\, G_{i,0}^\Sigma \times G_{j,0}^\Sigma \,$,
        \begin{equation} \label{eq:MAXISG2}
         M_1^{}~=~\prod_{\genfrac{}{}{0pt}{}{k=1}{k \neq i,j}}^r G_{k,0}^\Sigma \,
                  \times \, G_{ij,0}^\Sigma~,
        \end{equation}
        where ``diagonal'' means that under suitable\/ $\Sigma$-equivariant
        isomorphisms \linebreak $G_{i,0}^\Sigma \cong G_{s,0}^\Sigma \,$
        and $\, G_{j,0}^\Sigma \cong G_{s,0}^\Sigma \,$, the subgroup\/
        $G_{ij,0}^\Sigma$ of $\, G_{i,0}^\Sigma \times G_{j,0}^\Sigma \,$
        corresponds to the subgroup
        \[
         \mathrm{diag} \; G_{s,0}^\Sigma~
         =~\{ (g_0^{},g_0^{})~|~g_0^{} \in G_{s,0}^\Sigma \}
        \]
        of~$\, G_{s,0}^\Sigma \times G_{s,0}^\Sigma \,$.
        We then say that\/ $M$ and\/ $M_1^{}$ are of\,
        \textbf{$\Sigma$-diagonal type}.
 \end{itemize}
 Moreover, $M$ and\/~$M_1$ will be of normal type if and only if they are
 of\/ $\Sigma$-simple type and\/ $M_{i,1}^\Sigma$ is discrete; in all other
 cases, $M$ and\/~$M_1$ will be of normalizer type.
\end{theorem}
\begin{theorem} \label{theo:maxsg4}
 Let\/ $G$ be a Lie group with connected one-component\/~$G_0$ and component
 group\/~$\Gamma$, and let $\, \Sigma \subset \mathrm{Out}(G_0) \,$ be the
 image of\/~$\Gamma$ under the homomorphism~(\ref{eq:GREXT4}).
 Assuming that\/ $G_0$ is $\Sigma$-simple with trivial center and that
 \begin{equation} \label{eq:DECSLG1}
  G_0^{}~=~G_s^{} \times \ldots \times G_s^{}
  \quad \mbox{($n$ factors)}
 \end{equation}
 is its canonical decomposition into the direct product of\/ $n$ copies
 of the same simple Lie group\/~$G_s^{}$, let\/ $M$ be a maximal subgroup
 of\/~$G$ with connected one-component\/ $M_0$ such that\/ $M$ meets
 every connected component of\/~$G$.
 Set $\, M_1 = M \cap G_0$, so $\, M_0 \subset M_1 \subset M$.
 Then one of the following two alternatives holds:
 \begin{enumerate}
  \item Up to ``twisting'' with an appropriate automorphism of\/~$G_0^{}$,
        $M_1^{}$ is the direct \mbox{product} of\/ $n$ copies of the same
        maximal\/ $\Sigma_s^{}$-invariant subgroup\/ $M_{s,1}^{}$ of\/~$G_s^{}$,
        \begin{equation} 
         M_1^{}~
         =~M_{s,1}^{} \, \times \ldots \times \, M_{s,1}^{}
         \quad \mbox{($n$ factors)}~,
        \end{equation}
        where the subgroup\/ $\Sigma_s^{}$ of\/~$\mathrm{Out}(G_s^{})$ is
        obtained from the subgroup\/ $\Sigma$ of\/~$\mathrm{Out}(G_0^{})$
        by projection.
        We then say that\/ $M$ and\/~$M_1^{}$ are of\, \textbf{simple type}.
  \item Up to ``twisting'' with an appropriate automorphism of\/~$G_0^{}$,
        $M_1^{}$ is the direct \mbox{product} of a certain number (say\/ $q$)
        of copies of the diagonal subgroup of the direct product of a certain
        number (say\/ $p$) of copies of\/~$G_s^{}$,
        \begin{equation} \label{eq:MAXISG5}
         \begin{array}{cc}
          M_1^{}~=~\mathrm{diag}_p \; G_s^{} \times \ldots \times
                   \mathrm{diag}_p \; G_s^{}
          \quad \mbox{($q$ factors)} \\[2mm]
          \mathrm{diag}_p \; G_s^{}~
          =~\{ (g_s^{},\ldots,g_s^{})~|~g_s^{} \in G_s^{} \}
          \quad \mbox{($p$ factors)}
         \end{array}~,
        \end{equation}
        where\/ $p$ and\/ $q$ are divisors of~$n$, with $\, p>1 \,$ and
        $\, q<n$, chosen such that\/ $p$ is the minimum and\/ $q$ the
        maximum possible value for which the resulting partition
        \begin{equation} \label{eq:MAXISG6}
        \{1,\ldots,n\}~=~\big\{ \{1,\ldots,p\} , \ldots ,
                                \{n-p+1,\ldots,n\} \bigr\}
        \end{equation}
        is\/ $\Sigma$-invariant.
        We then say that\/ $M$ and\/~$M_1^{}$ are of\, \textbf{diagonal type}.
 \end{enumerate}
 Moreover, $M$ and\/~$M_1$ will be of normal type if and only if they are of
 simple type and\/ $M_{s,1}$ is discrete; in all other cases, $M$ and\/~$M_1$
 will be of normalizer type.
\end{theorem}
Regarding the proof of Theorem~\ref{theo:maxsg4}, there is only one case
which is not fully covered by Theorem~\ref{theo:maxisa2}, namely the one
where $M$ and~$M_1$ and~$M_{s,1}$ are discrete.
However, it is a simple exercise to perform the necessary adaptations,
observing that when $G_0$ is a connected semisimple Lie group with Lie
algebra~$\mathfrak{g}$, then if $G_0$ is simply connected or~-- as in
the case of interest here~-- if $G_0$ is centerfree, every automorphism
of~$\mathfrak{g}$ can be lifted to an automorphism of~$G_0$ and hence the
groups $\mathrm{Aut}(G_0)$ and $\mathrm{Aut}(\mathfrak{g})$ and similarly
the groups $\mathrm{Out}(G_0)$ and $\mathrm{Out}(\mathfrak{g})$ are
canonically isomorphic.
\vspace{1ex}

It is instructive to illustrate the phenomena that appear in
Theorems~\ref{theo:maxsg3} and~\ref{theo:maxsg4} by examples.
The first of them will also show that the maximal subgroups of a
non-connected compact Lie group may differ considerably from those
of its connected one-component.
\begin{example} \label{ex:SO4FIN}
 Consider the group $O(4)$ of all orthogonal transformations
 in~$\mathbb{R}^4$, whose connected one-component is the special
 orthogonal group $SO(4)$. The latter has a non-trivial center $Z$
 and a non-trivial outer automorphism group $\Sigma$, both isomorphic
 to~$\mathbb{Z}_2$: the non-trivial element of $Z$ is the matrix $-1_4$,
 whereas the non-trivial element of $\Sigma$ can be represented by
 conjugation with the reflection matrix $\, \mathrm{diag}(1,-1,%
 -1,-1)$, say: thus $SO(4)$ is not simple but is $\Sigma$-simple.
 As explained before Theorem~\ref{theo:maxsg3}, the first step
 consists in dividing by $Z$, descending to the adjoint group
 $\, SO(3) \times SO(3)$, on which the non-trivial element
 of~$\Sigma$ acts by switching the factors.
 (In passing, we note that $SO(4)$ lies ``in between'' its
 universal covering group $\, SU(2) \times SU(2) \,$ and its
 adjoint group $\, SO(3) \times SO(3)$; more precisely,
 $SO(4) \cong (SU(2) \times SU(2))/\mathbb{Z}_2 \,$
 with $\, \mathbb{Z}_2 = \{\,(1_2,1_2),(-1_2,-1_2)\} \,$
 and $\; SO(4)/Z \cong SO(3) \times SO(3)$.)%
 \footnote{Note that in this example, we encounter no less than three
 different $\mathbb{Z}_2$ groups, which must be clearly distinguished
 since they play very different roles.}
 Therefore, using the classification of the maximal subgroups of $SO(3)$
 given in Example~\ref{ex:SO3}~-- according to which there are two of
 normal type (which, as we recall, must be finite, since $SO(3)$ is
 simple), namely the cubic/octahedral group~$\mathrm{O}$ and the
 dodecahedral/icosahedral group~$\mathrm{I}$, while there is only
 one of normalizer type, namely~$O(2)$~-- we see that Theorems~%
 \ref{theo:maxsg3} and~\ref{theo:maxsg4} provide the following
 list of maximal subgroups and of maximal $\Sigma$-invariant
 subgroups of~$SO(4)/Z$ (which give the corresponding ones
 of~$SO(4)$ by taking the inverse image under the quotient
 homomorphism):

 \begin{table}[!htb]
  \begin{center}
  \begin{tabular}{|c|c|c|} \hline
   \begin{minipage}{4em}
    \begin{center}
     \rule{0mm}{3ex} Type \rule{0mm}{3ex}
    \end{center}
   \end{minipage}
   &
   \begin{minipage}{6em}
    \begin{center}
     \rule{0mm}{3ex} Maximal \rule{0mm}{3ex} \\ subgroups
    \end{center}
   \end{minipage}
   &
   \begin{minipage}{6em}
    \begin{center}
     \rule{0mm}{3ex} Maximal \rule{0mm}{3ex} \\
     $\Sigma$-invariant \\ subgroups
    \end{center}
   \end{minipage}
   \\[4ex] \hline\hline
   \rule{0mm}{5ex} simple/normal \rule{0mm}{5ex} &
   $\begin{array}{ccc}
     \mathrm{O} \times SO(3) \!\!&,&\!\! SO(3) \times \mathrm{O} \\
     \mathrm{I} \times SO(3) \!\!&,&\!\! SO(3) \times \mathrm{I}
    \end{array}$
   & 
   $\mathrm{O} \times \mathrm{O}$, $\mathrm{I} \times \mathrm{I}$
   \\[2ex] \hline
   \rule{0mm}{4ex} simple/normalizer \rule{0mm}{4ex} &
   $\begin{array}{ccc}
     O(2) \times SO(3) \!\!&,&\!\! SO(3) \times O(2)
    \end{array}$
   &
   $O(2) \times O(2)$
   \\[1ex] \hline
   \rule{0mm}{4ex} diagonal \rule{0mm}{4ex} &
   $\mathrm{diag}_2 \; SO(3)$
   &
   $\mathrm{diag}_2 \; SO(3)$
   \\[1ex] \hline
  \end{tabular}
  \end{center}
  \begin{center}
   \caption{\label{tab:SO4} \hspace*{\fill}
            Maximal subgroups and maximal $\Sigma$-invariant subgroups of
            $\, SO(4)/Z$ \newline
            \centerline{($SO(4)/Z \cong SO(3) \times SO(3)$)} \newline
            \centerline{($Z \cong \mathbb{Z}_2 \,$ = center of $SO(4)$,
                         $\Sigma \cong \mathbb{Z}_2 \,$ = outer
                         automorphism group of $SO(4)$)}}
  \end{center}
 \end{table}
 \vspace*{-4ex}
\end{example}

\pagebreak

\begin{example}
 Assume $G_0$ to be a centerfree connected semisimple Lie group which is
 the direct product of four copies of the same centerfree connected simple
 Lie group~$G_s$,
 \[
  G_0~=~G_s \times G_s \times G_s \times G_s~,
 \]
 and suppose that $\, \Gamma = \Sigma \,$ is the cyclic group $\mathbb{Z}_4$.
 Then $G_0$ is $\mathbb{Z}_4$-simple since $\mathbb{Z}_4$ acts transitively
 on the set $\{1,2,3,4\}$, and it is easily checked that the maximal
 $\mathbb{Z}_4$-invariant subgroup of~$G_0$ of diagonal type is
 \[
  \{ \, (g_1,g_2,g_1,g_2)~|~g_1,g_2 \in G_s \, \}
 \]
 which is isomorphic to $\, G_s \times G_s \,$ and conjugate to the
 subgroup $\, \mathrm{diag}_2 \, G_s^{} \,$ of equation~(\ref{eq:MAXISG5})
 under the transposition $\, 2 \leftrightarrow 3$.
\end{example}

\section{Real and complex primitive subalgebras}

As pointed out in the introduction, we shall for the rest of this paper
concentrate on maximal subgroups of compact Lie groups which are of
normalizer type, that is, maximal subgroups which are the normalizer
of their own Lie subalgebra.
From the theorems proved in the previous section, it follows that the
problem of determining all such subgroups can be reduced to that of
classifying the $\Sigma$-primitive subalgebras of compact simple Lie
algebras $\mathfrak{g}$, where $\Sigma$ is a subgroup of the (finite)
outer automorphism group $\mathrm{Out}(\mathfrak{g})$ of~$\mathfrak{g}$.
The best way to do so is by transferring the corresponding classification
for complex simple Lie algebras, which (at least for the case of trivial
$\Sigma$) is known, to their compact real forms.
This requires studying the behavior of primitive subalgebras of a compact
simple Lie algebra under complexification.

It is well known that the complexification of a real Lie algebra
defines a functor from the category of compact Lie algebras to the
category of complex reductive Lie algebras, both with homomorphisms
of Lie algebras as morphisms.
Moreover, this functor establishes a bijective correspondence
\begin{equation}
 \textrm{compact Lie algebras}~\longleftrightarrow~
 \textrm{complex reductive Lie algebras}
\end{equation}
at the level of isomorphism classes, since according to the Weyl
existence and conjugacy theorem for compact real forms~\cite{Kn},
every complex reductive Lie algebra admits compact real forms and
any two of these are conjugate.

If we fix a compact real form $\mathfrak{g}$ of a complex reductive
Lie algebra, denoted by~$\mathfrak{g}^\C$ to indicate that it is the
complexification of~$\mathfrak{g}$, then the above functor induces a
bijective correspondence between the lattice of conjugacy classes of
subalgebras of~$\mathfrak{g}$ and the lattice of conjugacy classes of
reductive subalgebras of~$\mathfrak{g}^\C$:
\begin{equation}
 \begin{minipage}{54mm}
  \begin{center}
   \textrm{lattice of conjugacy classes of} \\
   \textrm{subalgebras of the} \\
   \textrm{compact Lie algebra} $\mathfrak{g}$
  \end{center}  
 \end{minipage}
 ~\longleftrightarrow~
 \begin{minipage}{60mm}
  \begin{center}
   \textrm{lattice of conjugacy classes of} \\
   \textrm{reductive subalgebras of the} \\
   \textrm{complex reductive Lie algebra} $\mathfrak{g}^\C$
  \end{center}  
 \end{minipage}
\end{equation}
Now note that complex reductive Lie algebras contain non-reductive
subalgebras, such as the parabolic subalgebras.
Taking into account that every compact Lie algebra is reductive and
so is its complexification, these non-reductive subalgebras cannot be
obtained as the complexification of any subalgebra of any compact real
form.
Thus in order to obtain the primitive subalgebras of a compact simple
Lie algebra from the primitive subalgebras of its complexification,
we must restrict ourselves to reductive subalgebras and therefore
we should relativize the notions of maximal and (quasi)primitive
to the lattice of reductive subalgebras.
This leads to the following modification of Definitions~\ref{def:maxisa}
and~\ref{def:primsa}.
\begin{Definition} \label{def:maxrsal}
 Let $\mathfrak{g}$ be a reductive Lie algebra and $\Sigma$ be a sub%
 group of its outer auto\-morphism group $\mathrm{Out}(\mathfrak{g})$.
 A \emph{maximal $\Sigma$-invariant reductive subalgebra} of~$\mathfrak{g}$
 is a proper $\Sigma$-invariant reductive subalgebra $\mathfrak{m}$ of~%
 $\mathfrak{g}$ such that if $\tilde{\mathfrak{m}}$ is any $\Sigma$-%
 invariant reductive subalgebra of~$\mathfrak{g}$ with $\, \mathfrak{m}
 \subset \tilde{\mathfrak{m}} \subset \mathfrak{g}$, then
 $\, \tilde{\mathfrak{m}} = \mathfrak{m} \,$ or
 $\, \tilde{\mathfrak{m}} = \mathfrak{g}$.
 A \emph{$\Sigma$-quasiprimitive reductive subalgebra} of~$\mathfrak{g}$
 is a proper $\Sigma$-invariant reductive subalgebra $\mathfrak{m}$ of~%
 $\mathfrak{g}$ which is maximal among all $\mathrm{Aut}_\Sigma^{}%
 (\mathfrak{g},\mathfrak{m})$-invariant reductive subalgebras of~%
 $\mathfrak{g}$, that is, such that if $\tilde{\mathfrak{m}}$ is
 any $\mathrm{Aut}_\Sigma^{}(\mathfrak{g},\mathfrak{m})$-invariant
 reductive subalgebra of~$\mathfrak{g}$ with $\, \mathfrak{m} \subset
 \tilde{\mathfrak{m}} \subset \mathfrak{g}$, then $\, \tilde{\mathfrak{m}}
 = \mathfrak{m} \,$ or $\, \tilde{\mathfrak{m}} = \mathfrak{g}$.
 A \emph{$\Sigma$-primitive reductive subalgebra} of~$\mathfrak{g}$
 is a $\Sigma$-quasiprimitive reductive subalgebra of~$\mathfrak{g}$
 which contains no non-trivial proper $\Sigma$-invariant ideal
 of~$\mathfrak{g}$.
 When $\Sigma$ is the image under the homomorphism~(\ref{eq:GREXT4})
 of the component group $\Gamma$ of a Lie group~$G$ with Lie algebra
 $\mathfrak{g}$, we also use the term ``maximal $\Gamma$-invariant''
 as a synonym for ``maximal $\Sigma$-invariant'' and the term
 ``$\Gamma$-(quasi)primitive reductive'' as a synonym for
 ``$\Sigma$-(quasi)primitive reductive'', and when $\Sigma$
 is trivial ($\Sigma = \{1\}$), we omit the reference to this
 group and simply speak of a maximal reductive subalgebra and
 of a (quasi)primitive reductive subalgebra, respectively.
\end{Definition}
\noindent
With this terminology, we have the following
\begin{proposition} \label{thm:COMPACT}
 Let $\mathfrak{g}$ be a compact Lie algebra and $\mathfrak{g}^\C$
 be its complexification, and let $\Sigma$ be a subgroup of their
 outer automorphism group $\, \mathrm{Out}(\mathfrak{g}) \cong
 \mathrm{Out}(\mathfrak{g})^\C$.
 Then there is a bijective correspondence between the $\Sigma$-%
 (quasi)primitive subalgebras of~$\mathfrak{g}$ and the $\Sigma$-%
 (quasi)primitive reductive subalgebras of~$\mathfrak{g}^\C$ which
 preserves conjugacy classes of subalgebras.
 Moreover, this correspondence takes maximal $\Sigma$-invariant
 subalgebras of~$\mathfrak{g}$ to maximal $\Sigma$-invariant
 reductive subalgebras of~$\mathfrak{g}^\C$.
\end{proposition}
\begin{proof}
 The proof is straightforward and is left to the reader.
\end{proof}

\section{Primitive subalgebras of classical Lie algebras}

In this section, we summarize the classification of the primitive subalgebras
and, more generally, the $\Sigma$-primitive subalgebras of compact classical
Lie algebras.
Our presentation is based on a combination of results obtained by various
authors \cite{Ch,Dy1,Dy2,Go,Ko3}.
As has already been mentioned in the introduction, all these papers, with
the exception of parts of Ref.~\cite{Ko3}, refer to the complex case.
Since we are ultimately interested in compact Lie groups, we shall
restate them by using Proposition~\ref{thm:COMPACT} to translate
to the respective compact real forms.
Correspondingly, the term ``classical Lie algebra'' will in what follows
mean one of the compact classical Lie algebras in the standard representation
as given by equation~(\ref{eq:CLASLA}) below.

When working with classical Lie algebras the adequate method for analyzing
the inclusion of subalgebras is to use their standard realization as matrix
Lie algebras~-- more precisely, as Lie algebras of complex $(n \times n)$-%
matrices.
Concretely,
\begin{equation} \label{eq:CLASLA}
\begin{array}{ccc}
 \mbox{$A$-series} & \quad & \mathfrak{su}(n)~
 =~\{ \, X \in \mathfrak{u}(n)~|~\mathrm{tr}(X) = 0 \, \} \\
 \mbox{($A_r$, $r \geqslant 1$)} & \quad & \mbox{($n = r+1$)} \\[2mm]
 \mbox{$B$-series} & \quad & \mathfrak{so}(n)~
 =~\{ \, X \in \mathfrak{u}(n)~|~X^T + X = 0 \, \} \\
 \mbox{($B_r$, $r \geqslant 2$)} & \quad & \mbox{($n = 2r+1$ odd)} \\[2mm]
 \mbox{$C$-series} & \quad & \mathfrak{sp}(n)~
 =~\{ \, X \in \mathfrak{u}(n)~|~X^T J + J X = 0 \, \} \\
 \mbox{($C_r$, $r \geqslant 3$)} & \quad & \mbox{($n = 2r$ even)} \\[1mm]
 \mbox{$D$-series} & \quad & \mathfrak{so}(n)~
 =~\{ \, X \in \mathfrak{u}(n)~|~X^T + X = 0 \, \} \\
 \mbox{($D_r$, $r \geqslant 4$)} & \quad & \mbox{($n = 2r$ even)}
\end{array}
\end{equation}
where
\begin{equation}
 \mathfrak{u}(n)~=~\{ \, X \in \mathfrak{gl}(n,\mathbb{C})~|~
                         X^\dagger + X = 0 \, \}
\end{equation}
is the Lie algebra of antihermitean complex $(n \times n)$-matrices (as usual,
the symbols $.^T$ and $.^\dagger$ denote transpose and hermitean adjoint,
respectively),
\begin{equation}
 J~=~\left( \begin{array}{cc}
             0 & -1_r \\
             1_r & 0
            \end{array} \right) \mbox{ or }
     \left( \begin{array}{cc}
             0 & 1_r \\
             -1_r & 0
            \end{array} \right)
\end{equation}
is the standard symplectic $(2r \times 2r)$-matrix (the position of the
$-$ sign being a matter of taste), and the constraints on the values
of~$r$, or~$n$, in equation~(\ref{eq:CLASLA}) are imposed in order to
avoid repetitions due to the well-known canonical isomorphisms between
classical Lie algebras of low rank ($A_1 = B_1 = C_1 \,$, $D_1$ is abelian,
$B_2 = C_2 \,$, $D_2 = A_1 \oplus A_1 \,$, $A_3 = D_3$).
We also note that the outer automorphism groups of the simple Lie
algebras (including the five exceptional ones) are
\begin{equation} \label{eq:OUTSIM}
 \mathrm{Out}(\mathfrak{g})~=~\left\{ \begin{array}{ccc}
            \{1\} & \mbox{for $A_1$, $B_r$, $C_r$,
                              $E_7$, $E_8$, $F_4$, $G_2$} \\[1mm]
            \mathbb{Z}_2 & \mbox{for $A_r$ ($r \geqslant 2$),
                                     $D_r$ ($r \geqslant 5$), $E_6$} \\[1mm]
            S_3 & \mbox{for $D_4$}
           \end{array} \right\}~.
\end{equation}
This implies that, apart from the exceptional case of the algebra $D_4$
($\mathfrak{so}(8)$), $\Sigma$-primitive subalgebras of compact simple Lie
algebras come in just two types: primitive and $\mathbb{Z}_2$-primitive. 
Clearly, the latter exist only within the algebras $A_r$ ($r \geqslant 2$),
$D_r$ ($r \geqslant 4$) and $E_6$, and except for $D_4$, they coincide with
the almost primitive subalgebras of Ref.~\cite{Ko3}.
In the exceptional case of the algebra $D_4$ ($\mathfrak{so}(8)$), there are
other possibilities, since $\Sigma$ can be chosen to be one of the other two
$\mathbb{Z}_2$-subgroups of $S_3$, the cyclic subgroup $\mathbb{Z}_3$ or all
of~$S_3$; here, it is the $S_3$-primitive subalgebras that coincide with the
almost primitive subalgebras of Ref.~\cite{Ko3}.
In what follows, we shall not deal with these other types of $\Sigma$-%
primitive subalgebras of~$D_4$, nor with the $\mathbb{Z}_2$-primitive
subalgebras of~$E_6$: they require a separate analysis.
However, we shall present results for the remaining cases, which are the
great majority, corresponding to the algebras $A_r$ ($r \geqslant 2$)
and $D_r$ ($r \geqslant 5$): they will be handled by noting that any
automorphism $\, \phi \in \mathrm{Aut}(\mathfrak{g}) \,$ representing the
non-trivial element of $\, \mathrm{Out}(\mathfrak{g}) = \mathbb{Z}_2 \,$
can be realized as conjugation with an antiunitary transformation on~%
$\mathbb{C}^n$, in the case of~$\mathfrak{su}(n)$ ($n \geqslant 3$),
and as conjugation with an orthogonal transformation on~$\mathbb{R}^n$
of determinant $-1$, in the case of~$\mathfrak{so}(n)$ ($n$ even,
$n \geqslant 10$).
This procedure also works for~$\mathfrak{so}(8)$, except that it will in
this case not generate all of $\, \mathrm{Out}(\mathfrak{g}) = S_3 \,$
but only one of its three $\mathbb{Z}_2$-subgroups.

With these preliminaries out of the way, we shall divide the
set of sub\-algebras $\mathfrak{s}$ of a given classical Lie algebra
$\mathfrak{g}$ into several types, according to two natural criteria.
The~first criterion refers to the intrinsic nature of~$\mathfrak{s}$:
it can be
\begin{itemize}
 \item abelian, \vspace{-0.5ex}
 \item simple, \vspace{-0.5ex}
 \item truly semisimple, i.e., not simple, \vspace{-0.5ex}
 \item truly reductive, i.e., with non-trivial center and non-trivial
       derived sub\-algebra.
\end{itemize}
The second criterion refers to the nature of the inclusion of~$\mathfrak{s}$
in~$\mathfrak{g}$.
Since $\mathfrak{g}$ is a matrix algebra acting on an $n$-dimensional
(complex) vector space~$\, V \cong \mathbb{C}^n$, we can also think of
this inclusion as a (faithful) representation
\begin{equation} \label{eq:INCREP}
 \pi : \mathfrak{s}~\longrightarrow~
       \mathfrak{gl}(V) \cong \mathfrak{gl}(n,\mathbb{C})
\end{equation}
of~$\mathfrak{s}$ on~$V$, which must be of one of the following three types:
\begin{itemize}
 \item Type $0$: such representations are said not to be
       self-conjugate, or to be truly complex, and correspond
       to inclusions into $\, \mathfrak{g} = \mathfrak{su}(n)$.
       \vspace{-0.5ex}
 \item Type $+1$: such representations are said to be real, or orthogonal,
       and correspond to inclusions into $\, \mathfrak{g} = \mathfrak{so}(n)$;
       abstractly, elements of~$\mathfrak{so}(n)$ are characterized by the
       property of commuting with a given antilinear transformation $\tau$
       on~$V$ of square $+1$ (complex conjugation).
       \vspace{-0.5ex}
 \item Type $-1$: such representations are said to be pseudo-real, or
       quaternionic, or symplectic, and correspond to inclusions into
       $\, \mathfrak{g} = \mathfrak{sp}(n)$; abstractly, elements of~%
       $\mathfrak{sp}(n)$ are characterized by the property of commuting
       with a given antilinear transformation $\tau$ on~$V$ of square
       $-1$ (complex conjugation combined with application of~$J$).
\end{itemize}
There are then two distinct possibilities:
\begin{itemize}
 \item $\pi$ is irreducible, \vspace{-1ex}
 \item $\pi$ is reducible.
\end{itemize}

As was first observed by Dynkin in the case of maximal subalgebras,
the procedure for classifying the simple ones is rather different from
that for classifying the remaining ones.
As it turns out, the same goes for primitive and $\mathbb{Z}_2$-%
primitive subalgebras.

To explain Dynkin's strategy, note first of all that when $\mathfrak{s}\,$
is a simple maximal subalgebra of a classical Lie algebra $\mathfrak{g}$,
then with only one exception, $\mathfrak{s}\,$ is always irreducible, and
the same statement holds for primitive and $\mathbb{Z}_2$-primitive
subalgebras.
(The proof follows as a corollary from the general classification
of reducible primitive and $\mathbb{Z}_2$-primitive subalgebras of
classical Lie algebras carried out later in this section; see the
discussion following Lemma~\ref{lem:BFTP} below.)
What Dynkin realized was that conversely, almost every irreducible
representation of a simple Lie algebra $\mathfrak{s}\,$ provides an
inclusion of~$\mathfrak{s}\,$ as a maximal subalgebra of the pertinent
classical Lie algebra $\mathfrak{g}$, where ``almost all'' means that
there is only a handful of exceptions, which can be listed explicitly:
\begin{theorem}%
 \emph{(Dynkin~\cite{Dy1,Dy2}, Chekalov~\cite{Ch}, Komrakov~\cite{Ko3})}
 \label{theo:DYNEXC}
 Let\/ $\mathfrak{g}$ be a compact \mbox{classical} Lie algebra.
 Then every simple maximal subalgebra, every simple primitive
 subalgebra and every simple\/ $\mathbb{Z}_2$-primitive subalgebra\/
 $\mathfrak{s}$ of\/~$\mathfrak{g}$ is irreducible, with the only
 exception of the inclusions
 \begin{equation} \label{eq:RSMPSA}
  \mathfrak{so}(n-1) \subset \mathfrak{so}(n) \quad (n \geqslant 6)~.
 \end{equation}
 Conversely, every simple irreducible subalgebra $\mathfrak{s}$ is maximal,
 which implies that every simple irreducible subalgebra $\mathfrak{s}$ is
 primitive and every $\mathbb{Z}_2$-invariant simple irreducible subalgebra
 $\mathfrak{s}$ is $\mathbb{Z}_2$-primitive, unless the inclusion
 $\, \mathfrak{s} \subset \mathfrak{g} \,$ is one of the 18 exceptions
 listed in Table~1  of\/~\cite[p.~364]{Dy1} or Table~7 of\/~\cite[p.~236]{OV}.
 Among these exceptions, there is only one inclusion $\, \mathfrak{s}
 \subset \mathfrak{g} \,$ such that\/ $\mathfrak{s}$ is primitive,
 given by\/~\cite[p.~279]{Ch}, \cite[p.~200]{Ko3}
 \begin{equation} \label{eq:ISNPSA}
  \mathfrak{so}(12) \subset \mathfrak{so}(495)~,
 \end{equation}
 and none such that $\mathfrak{s}$ is $\mathbb{Z}_2$-primitive.
\end{theorem}
\begin{remark}
 Although the theorem does not provide an explicit list of simple maximal
 or primitive or $\mathbb{Z}_2$-primitive subalgebras of $\mathfrak{g}$,
 it is the cornerstone for finding all such subalgebras.
 Namely, fixing the Lie algebra $\mathfrak{g}$ with its natural representation
 on~$\mathbb{C}^n$, we apply representation theory~-- more specifically, the
 Weyl dimension formula~-- to first find all irreducible representations of
 all simple Lie algebras $\mathfrak{s}$ (classical and exceptional) of
 dimension~$n$.
 Next, we determine which among these provide inclusions of $\mathfrak{s}$
 in~$\mathfrak{g}$ by deciding whether the irreducible representation
 of~$\mathfrak{s}$ under consideration is self-conjugate, and if so,
 whether it is real (orthogonal) or pseudo-real/quaternionic (symplectic).
 Finally, we apply the theorem to eliminate the exceptions.
\end{remark}

To handle the general case of subalgebras which are not simple, we
must distinguish between two different constructions of the pertinent
representation, depending on whether the subalgebra in question is
reducible or irreducible:
\begin{Definition} \label{def:EXTDSTP}
 Let $\, \pi_1: \mathfrak{g}_1 \longrightarrow \mathfrak{gl}(V_1) \,$
 and $\, \pi_2: \mathfrak{g}_2 \longrightarrow \mathfrak{gl}(V_2) \,$
 be representations of Lie algebras $\mathfrak{g}_1$ and $\mathfrak{g}_2$
 on vector spaces $V_1$ and $V_2$, respectively.
 The \emph{external direct sum} of $\pi_1$ and $\pi_2$ is the
 representation $\, \pi_1 \boxplus \pi_2 \,$ of the Lie algebra
 $\, \mathfrak{g}_1 \oplus \mathfrak{g}_2 \,$ on the direct sum
 $V_1 \oplus V_2$ of the vector spaces $V_1$ and $V_2$ given by
 \begin{equation} \label{eq:EXTDS}
  \begin{array}{cccc}
   \pi_1 \boxplus \pi_2 :
   & \mathfrak{g}_1 \oplus \mathfrak{g}_2 & \longrightarrow &
     \mathfrak{gl}(V_1 \oplus V_2) \\[2mm]
   &              X_1 + X_2               &   \longmapsto   &
     \pi_1(X_1) \oplus \pi_2(X_2)
  \end{array}~.
 \end{equation}
 The \emph{external tensor product} of $\pi_1$ and $\pi_2$ is the
 representation $\, \pi_1 \boxtimes \pi_2 \,$ of the Lie algebra
 $\, \mathfrak{g}_1 \times \mathfrak{g}_2 \,$ on the tensor product
 $V_1 \otimes V_2$ of the vector spaces $V_1$ and $V_2$ given by
 \begin{equation} \label{eq:EXTTP}
  \begin{array}{cccc}
   \pi_1 \boxtimes \pi_2 :
   & \mathfrak{g}_1 \times \mathfrak{g}_2 & \longrightarrow &
     \mathfrak{gl}(V_1 \otimes V_2) \\[2mm]
   &              (X_1,X_2)               &   \longmapsto   &
     \pi_1(X_1) \otimes 1 \, + \, 1 \otimes \pi_2(X_2)
  \end{array}~.
 \end{equation}
\end{Definition}
\begin{remark}
 Observe that the resulting Lie algebra is the same in both cases,
 namely the direct sum of $\mathfrak{g}_1$ and $\mathfrak{g}_2$.
 Still, in the case of the tensor product, we prefer to use the symbol
 $\times$, rather than $\oplus$, since experience shows that this helps
 to avoid confusion.
 (On the other hand, we avoid the symbol $\otimes$ in this context,
 since the concept of tensor product of Lie algebras does not exist.)
 In particular, using both notations in parallel is helpful to
 characterize inclusions of subalgebras by representations: the
 symbol $\oplus$ indicates that the inclusion into $\mathfrak{g}$
 is by the direct sum construction, whereas the symbol $\times$
 indicates that the inclusion into $\mathfrak{g}$ is by the tensor
 product construction.
 This will greatly simplify the tables.
\end{remark}
The main difference between the two constructions lies in the fact that,
even when both $\pi_1$ and $\pi_2$ are irreducible, the external direct
sum $\, \pi_1 \boxplus \pi_2 \,$ is always a reducible representation
of $\, \mathfrak{g}_1 \oplus \mathfrak{g}_2$, whereas the external
tensor product $\, \pi_1 \boxtimes \pi_2$ \linebreak is always an
irreducible representation of $\, \mathfrak{g}_1 \times \mathfrak{g}_2$;
moreover, it can be shown that every irreducible representation of
$\, \mathfrak{g}_1 \times \mathfrak{g}_2 \,$ is given by this construction.
(See~\cite[Proposition~3.1.8, p.~123]{GW} for the proof of an entirely
analogous statement for groups, which can be adapted to Lie algebras
replacing the group algebra by the universal enveloping algebra.)
Finally, by iteration of both constructions, it is clear that we can
define the direct sum and the tensor product of an arbitrary finite
number of representations.

When working with external direct sums and tensor products, it is useful
to consider the following generalization of the concept of equivalence
between representations, adapted from~\cite[p.~55]{KL}:
\begin{Definition} \label{def:EQMAUT}
 Two representations $\, \pi_1: \mathfrak{g}_1 \longrightarrow \mathfrak{gl}%
 (V_1) \,$ and $\, \pi_2: \mathfrak{g}_2 \longrightarrow \mathfrak{gl}(V_2)$
 \linebreak
 of Lie algebras $\mathfrak{g}_1$ and $\mathfrak{g}_2$ on vector spaces $V_1$
 and $V_2$, respectively, are said to be \emph{quasiequivalent} if there
 is a pair $(\phi_{21},g_{21})$ consisting of a Lie algebra isomorphism
 $\, \phi_{21}: \mathfrak{g}_1 \longrightarrow \mathfrak{g}_2 \,$ and a
 linear isomorphism $\, g_{21}: V_1 \longrightarrow V_2 \,$ such that
 \begin{equation} \label{eq:EQMAUT}
  \pi_2^{}(\phi_{21}^{}(X_1^{}))~=~g_{21}^{} \, \pi_1^{}(X_1^{}) \, g_{21}^{-1}
 \quad \mbox{for $\, X_1^{} \in \mathfrak{g}_1^{}$}~.
 \end{equation}
\end{Definition}
\noindent
In this case, we may of course assume, without loss of generality,
that $\, \mathfrak{g}_1 = \mathfrak{g} = \mathfrak{g}_2$ \linebreak
and $\, \phi \in \mathrm{Aut}(\mathfrak{g})$.
Obviously, the usual notion of equivalence is recovered when $\, \phi
\in \mathrm{Inn}(\mathfrak{g})$, which motivates the terminology.

The concept of quasiequivalence is needed to formulate an important case
distinction in the calculation of normalizers for direct sums and tensor
products of irreducible representations.%
\footnote{In what follows, $\mathbb{C}^\times$ denotes the multiplicative
group of non-zero complex numbers.}
\begin{lemma} \label{lem:CNDS}
 Let $\, \pi_1: \mathfrak{g}_1 \longrightarrow \mathfrak{gl}(V_1) \,$
 and $\, \pi_2: \mathfrak{g}_2 \longrightarrow \mathfrak{gl}(V_2) \,$
 be irreducible representations of Lie algebras\/ $\mathfrak{g}_1$ and\/
 $\mathfrak{g}_2$ on finite-dimensional complex vector spaces\/ $V_1$
 and\/ $V_2$, respectively, and let $\, \pi: \mathfrak{g} \longrightarrow
 \mathfrak{gl}(V) \,$ be their external direct sum: $\, \pi = \pi_1 \boxplus
 \pi_2 \,$, $\mathfrak{g} = \mathfrak{g}_1 \oplus \mathfrak{g}_2 \,$,
 $V = V_1 \oplus V_2 \,$.
 For simplicity, assume\/ $\pi_1$ and\/~$\pi_2$ to be faithful (if not,
 replace\/ $\mathfrak{g}_1$ by\/ $\mathfrak{g}_1/\ker \pi_1$ and\/
 $\mathfrak{g}_2$ by\/ $\mathfrak{g}_2/\ker \pi_2$).
 Then considering\/ $\mathfrak{g}_1$, $\mathfrak{g}_2$ and\/~$\mathfrak{g}$
 as matrix Lie algebras (more precisely, $\mathfrak{g}_i \subset \mathfrak{g}
 \subset \mathfrak{gl}(V)$ \linebreak via\/~$\pi$ and also $\, \mathfrak{g}_i
 \subset \mathfrak{gl}(V_i) \,$ via\/~$\pi_i$), their centralizers and
 normalizers in the respective general linear groups are related as follows.
 \begin{enumerate}
  \item For the centralizers,
        \begin{equation} \label{eq:CENTDS1}
         \begin{array}{c}
          Z_{GL(V)}^{}(\mathfrak{g}_1)~
          =~\mathbb{C}^\times \, \mathrm{id}_{V_1}^{} \oplus\, GL(V_2)~, \\[2mm]
          Z_{GL(V)}^{}(\mathfrak{g}_2)~
          =~GL(V_1) \oplus\, \mathbb{C}^\times \, \mathrm{id}_{V_2}^{}~,
         \end{array}
        \end{equation}
        and taking the intersection,
        \begin{equation} \label{eq:CENTDS2}
          Z_{GL(V)}^{}(\mathfrak{g})~
          =~\mathbb{C}^\times \, \mathrm{id}_{V_1}^{} \oplus\,
            \mathbb{C}^\times \, \mathrm{id}_{V_2}^{}~.
        \end{equation}
  \item For the normalizers,
        \begin{equation} \label{eq:NORMDS1}
         \begin{array}{c}
          N_{GL(V)}^{}(\mathfrak{g}_1)~
          =~N_{GL(V_1)}^{}(\mathfrak{g}_1) \oplus\, GL(V_2)~, \\[2mm]
          N_{GL(V)}^{}(\mathfrak{g}_2)~
          =~GL(V_1) \oplus N_{GL(V_2)}^{}(\mathfrak{g}_2)~,
         \end{array}
        \end{equation}
        whereas
        \begin{equation} \label{eq:NORMDS2}
          N_{GL(V)}^{}(\mathfrak{g})~
          =~N_{GL(V_1)}^{}(\mathfrak{g}_1) \oplus
            N_{GL(V_2)}^{}(\mathfrak{g}_2)~,
        \end{equation}
        except when\/ $\pi_1$ and\/ $\pi_2$ are quasiequivalent, in which case
        \begin{equation} \label{eq:NORMDS3}
          N_{GL(V)}^{}(\mathfrak{g})~
          =~\bigl( N_{GL(V_1)}^{}(\mathfrak{g}_1) \oplus
                   N_{GL(V_2)}^{}(\mathfrak{g}_2) \bigr) : S_2~.
        \end{equation}
 \end{enumerate}
 The generalization to more than two direct summands is the obvious one.
\end{lemma}
\noindent
Here, for the sake of brevity, we use the following notation: given any two
subgroups $H_1$ of $GL(V_1)$ and $H_2$ of $GL(V_2)$, $H_1 \oplus H_2 \,$
will denote the subgroup of $GL(V_1 \oplus V_2)$ given by
\[
 H_1 \oplus H_2~=~\{ \, h_1 \oplus h_2~|~h_1 \in H_1 \,,\, h_2 \in H_2 \, \}~.
\]
\begin{proof} \hspace*{-2em}
 \footnote{The explicit proof given here serves as a kind of ``warm-up
 exercise'' for the proof of the next lemma.}
 The proof of the inclusions ``$\supset$'' in equations~(\ref{eq:CENTDS1})-%
 (\ref{eq:NORMDS3}) is straightforward provided we properly specify the action
 of~$S_2$ involved in the definition of the semi\-direct product in equation~%
 (\ref{eq:NORMDS3}) which, roughly speaking, consists in switching the
 direct summands.
 More precisely, choosing the pair $(\phi_{21},g_{21})$ as in Definition~%
\ref{def:EQMAUT} and setting $\, \phi_{12}^{} = \phi_{21}^{-1}$, $g_{12}^{}
 = g_{21}^{-1}$, we define an involutive automorphism $\phi$ on
 $\, \mathfrak{g} = \mathfrak{g}_1 \oplus\, \mathfrak{g}_2 \,$ and
 a linear involution $g$ on~$\, V = V_1 \oplus V_2 \,$ by setting
 \begin{equation} \label{eq:SWITCHDS}
  \begin{array}{cccc}
   \phi: & \mathfrak{g}_1 \oplus\, \mathfrak{g}_2 & \longrightarrow &
           \mathfrak{g}_1 \oplus\, \mathfrak{g}_2 \\[1mm]
         &                (X_1,X_2)               &   \longmapsto   &
             ( \phi_{12}(X_2) , \phi_{21}(X_1) )   \\[3mm]
   g: & V_1 \oplus V_2 & \longrightarrow &      V_1 \oplus V_2      \\[1mm]
      &    (v_1,v_2)   &   \longmapsto   & (g_{12} v_2 , g_{21} v_1)
  \end{array}
 \end{equation}
 and see that equation~(\ref{eq:EQMAUT}) becomes equivalent to the
 condition that $g$ normalizes $\mathfrak{g}$, since
 \begin{eqnarray*}
  (\phi(X) \, g)(v_1,v_2) \!\!
  &=&\!\! \bigl( (\phi_{12}(X_2) \, g_{12})(v_2) \,,
                 (\phi_{21}(X_1) \, g_{21})(v_1) \bigr)~, \\
  (g \, X)(v_1,v_2) \!\!
  &=&\!\! \bigl( (g_{12} \, X_2)(v_2) \,, (g_{21} \, X_1)(v_1) \bigr)~.
 \end{eqnarray*}
 To prove the converse inclusions ``$\subset$'', suppose first that
 $\, g \in GL(V) \,$ normalizes $\mathfrak{g}$ and write $\phi$ for
 the automorphism of~$\mathfrak{g}$ induced by conjugation with~$g$:
 \[
  g X~=~\phi(X) \, g \quad \mbox{for $\, X \in \mathfrak{g}$}~.
 \]
 Explicitly, this means that if we write $\phi$ in the form
 \begin{equation} \label{eq:BLOCK1}
  \phi \bigl( X_1 , X_2 \bigr)~
   =~\bigl( \phi_{11}(X_1) + \phi_{12}(X_2) \,,\,
            \phi_{21}(X_1) + \phi_{22}(X_2) \bigr)
 \end{equation}
 with $\, \phi_{11} \in \mathrm{End}(\mathfrak{g}_1)$, $\phi_{12} \in
 L(\mathfrak{g}_2,\mathfrak{g}_1)$, $\phi_{21} \in L(\mathfrak{g}_1,
 \mathfrak{g}_2)$, $\phi_{22} \in \mathrm{End}(\mathfrak{g}_2)$, and
 similarly $g$ in the form
 \begin{equation} \label{eq:BLOCK2}
   g \bigl( v_1 , v_2 \bigr)~
   =~\bigl( g_{11} v_1 + g_{12} v_2 \,,\, g_{21} v_1 + g_{22} v_2 \bigr)
 \end{equation}
 with $\, g_{11} \in \mathrm{End}(V_1)$, $g_{12} \in L(V_2,V_1)$,
 $g_{21} \in L(V_1,V_2)$, $g_{22} \in \mathrm{End}(V_2)$, then for
 all $\, X = (X_1,X_2) \,$ in $\, \mathfrak{g} = \mathfrak{g}_1
 \oplus \mathfrak{g}_2 \,$ and $\, (v_1,v_2) \in V_1 \oplus V_2$,
 the expression
 \[
  \bigl( g X \bigr) \bigl( v_1 , v_2 \bigr)~
  =~g \bigl( X_1 v_1 , X_2 v_2 \bigr)~
  =~\bigl( g_{11} X_1 v_1 + g_{12} X_2 v_2 \,,\,
           g_{21} X_1 v_1 + g_{22} X_2 v_2 \bigr)
 \]
 must be equal to the expression
 \begin{eqnarray*}
  \bigl( \phi(X) \, g \bigr) \bigl( v_1 , v_2 \bigr) \!\!
  &=&\!\! \phi(X) \bigl( g_{11} v_1 + g_{12} v_2 \,,\,
                         g_{21} v_1 + g_{22} v_2 \bigr) \\[1mm]
  &=&\!\! \Bigl( \bigl( \phi_{11}(X_1) + \phi_{12}(X_2) \bigr)
                 \bigl( g_{11} v_1 + g_{12} v_2 \bigr) \,, \\[-1mm]
  & &~           \bigl( \phi_{21}(X_1) + \phi_{22}(X_2) \bigr)
                 \bigl( g_{21} v_1 + g_{22} v_2 \bigr) \Bigr)
 \end{eqnarray*}
 leading us to the following system of equations, in which the ones in
 the first two lines are obtained by leaving $v_1$ arbitrary and setting
 $\, v_2 = 0$, with $X_1$ arbitrary and $\, X_2 = 0$ (left column) and
 with $\, X_1 = 0 \,$ and $X_2$ arbitrary (right column), while the ones
 in the last two lines are obtained by setting $\, v_1 = 0 \,$ and leaving
 $v_2$ arbitrary, with $\, X_1 = 0 \,$ and $X_2$ arbitrary (left column)
 and with $X_1$ arbitrary and $\, X_2 = 0$ (right column):
 \begin{equation} \label{eq:BLOCK3}
  \begin{array}{ccc}
   g_{11} X_1~=~\phi_{11}(X_1) \, g_{11} &,& \phi_{12}(X_2) \, g_{11}~=~0~,
   \\[2mm]
   g_{21} X_1~=~\phi_{21}(X_1) \, g_{21} &,& \phi_{22}(X_2) \, g_{21}~=~0~,
   \\[2mm]
   g_{12} X_2~=~\phi_{12}(X_2) \, g_{12} &,& \phi_{11}(X_1) \, g_{12}~=~0~,
   \\[2mm]
   g_{22} X_2~=~\phi_{22}(X_2) \, g_{22} &,& \phi_{21}(X_1) \, g_{22}~=~0~.
  \end{array}
 \end{equation}
 Now assume that $\, g \in GL(V) \,$ normalizes $\mathfrak{g}_1$:
 then these calculations can be applied if we take $\, X_2 = 0 \,$ with
 $\, \phi_{11} \in \mathrm{Aut}(\mathfrak{g}_1)$, $\phi_{21} \equiv 0 \,$
 and $\phi_{12}$, $\phi_{22}$ undetermined.
 Obviously, in this situation,
 the first equation in the first line states that $\, g_{11} \in GL(V_1) \,$
 normalizes~$\mathfrak{g}_1$, the first equation in the second line implies
 $\, g_{21} = 0$, the second equation in the third line implies $\, g_{12}
 = 0 \,$ since $\pi_1$ is irreducible, and all other equations are trivial:
 this proves the inclusion ``$\subset$'' in the first line of equation~%
 (\ref{eq:NORMDS1}) and, as an immediate corollary, that in the first
 line of equation~(\ref{eq:CENTDS1}) (the second one being completely
 analogous).
 \\[2mm]
 Finally, if $\, g \in GL(V) \,$ normalizes $\mathfrak{g}$, as supposed
 above, we can use the above system of equations to argue as follows.
 If $g_{11}$ ($g_{22}$) is invertible, then according to the second
 equation in the first (fourth) line, $\phi_{12} \equiv 0$ ($\phi_{21}
 \equiv 0$), and according to the first equation in the third (second)
 line, $g_{12} = 0$ ($g_{21} = 0$), since $\pi_2$ ($\pi_1$) is irreducible.
 But $g$ being invertible, this forces $g_{22}$ ($g_{11}$) to be
 invertible as well, so repeating the other part of the argument,
 we conclude that $g$ and hence $\phi$ are block diagonal.
 If on the other hand both $g_{11}$ and $g_{22}$ are not invertible,
 then $\ker g_{11}$ is a non-trivial $\mathfrak{g}_1$-invariant sub%
 space of~$V_1$ and $\ker g_{22}$ is a non-trivial $\mathfrak{g}_2$-%
 invariant subspace of~$V_2$, implying $\, g_{11} = 0 \,$ and
 $\, g_{22} = 0$, since $\pi_1$ and $\pi_2$ are irreducible.
 But $g$ being invertible, this forces $g_{12}$ and $g_{21}$
 to be invertible and, up to a constant multiple, to be each
 other's inverse, due to Schur's lemma, so we conclude that $g$
 and hence $\phi$ are block off-diagonal and $\pi_1$ and $\pi_2$
 are quasiequivalent.
\end{proof}

\noindent
Similarly, we have
\begin{lemma} \label{lem:CNTP}
 Let $\, \pi_1: \mathfrak{g}_1 \longrightarrow \mathfrak{gl}(V_1) \,$
 and $\, \pi_2: \mathfrak{g}_2 \longrightarrow \mathfrak{gl}(V_2) \,$
 be irreducible representations of Lie algebras\/ $\mathfrak{g}_1$ and\/
 $\mathfrak{g}_2$ on finite-dimensional complex vector spaces\/ $V_1$
 and\/~$V_2$, respectively, and let $\, \pi: \mathfrak{g} \longrightarrow
 \mathfrak{gl}(V) \,$ be their external tensor product: $\, \pi = \pi_1
 \boxtimes \pi_2 \,$, $\mathfrak{g} = \mathfrak{g}_1 \!\times
 \mathfrak{g}_2 \,$, $V = V_1 \otimes V_2 \,$.
 For simplicity, assume\/ $\pi_1$ and\/~$\pi_2$ to be faithful (if not,
 replace\/ $\mathfrak{g}_1$ by\/ $\mathfrak{g}_1/\ker \pi_1$ and\/
 $\mathfrak{g}_2$ by\/ $\mathfrak{g}_2/\ker \pi_2$) and traceless
 (this is automatic if\/ $\mathfrak{g}_1$ and\/ $\mathfrak{g}_2$
 are semisimple or, more generally, perfect).
 Then considering\/ $\mathfrak{g}_1$, $\mathfrak{g}_2$ and\/~%
 $\mathfrak{g}$ as Lie algebras of traceless matrices (more precisely,
 $\mathfrak{g}_i \subset \mathfrak{g} \subset \mathfrak{sl}(V) \,$
 via\/~$\pi$ and also $\, \mathfrak{g}_i \subset \mathfrak{sl}(V_i) \,$
 via\/~$\pi_i$), their centralizers and normalizers in the respective
 general linear groups are related as follows.
 \begin{enumerate}
  \item For the centralizers,
        \begin{equation} \label{eq:CENTTP1}
         \begin{array}{c}
          Z_{GL(V)}^{}(\mathfrak{g}_1^{})~
          =~\{ \,\mathrm{id}_{V_1}^{} \} \otimes\, GL(V_2^{})~, \\[2mm]
          Z_{GL(V)}^{}(\mathfrak{g}_2^{})~
          =~GL(V_1^{}) \otimes \{ \,\mathrm{id}_{V_2}^{} \}~,
         \end{array}
        \end{equation}
        and taking the intersection,
        \begin{equation} \label{eq:CENTTP2}
          Z_{GL(V)}^{}(\mathfrak{g})~=~\mathbb{C}^\times \, \mathrm{id}_V^{}~.
        \end{equation}
  \item For the normalizers,
        \begin{equation} \label{eq:NORMTP1}
         \begin{array}{c}
          N_{GL(V)}^{}(\mathfrak{g}_1^{})~
          =~N_{GL(V_1)}^{}(\mathfrak{g}_1^{}) \otimes\, GL(V_2^{})~, \\[2mm]
          N_{GL(V)}^{}(\mathfrak{g}_2^{})~
          =~GL(V_1^{}) \otimes N_{GL(V_2)}^{}(\mathfrak{g}_2^{})~,
         \end{array}
        \end{equation}
        and when both $\mathfrak{g}_1$ and $\mathfrak{g}_2$ are simple,
        \begin{equation} \label{eq:NORMTP2}
          N_{GL(V)}^{}(\mathfrak{g})~
          =~N_{GL(V_1)}^{}(\mathfrak{g}_1^{}) \otimes
            N_{GL(V_2)}^{}(\mathfrak{g}_2^{})~,
        \end{equation}
        except when\/ $\pi_1$ and\/ $\pi_2$ are quasiequivalent, in which case
        \begin{equation} \label{eq:NORMTP3}
          N_{GL(V)}^{}(\mathfrak{g})~
          =~\bigl( N_{GL(V_1)}^{}(\mathfrak{g}_1^{}) \otimes
                   N_{GL(V_2)}^{}(\mathfrak{g}_2^{}) \bigr) : S_2^{}~.
        \end{equation}
 \end{enumerate}
 The generalization to more than two tensor factors is the obvious one.
\end{lemma}
\noindent
Here, for the sake of brevity, we use the following notation: given any two
subgroups $H_1$ of $GL(V_1)$ and $H_2$ of $GL(V_2)$, $H_1 \otimes H_2 \,$
will denote the subgroup of $GL(V_1 \otimes V_2)$ given by
\[
 H_1 \!\otimes H_2~=~\{ \, h_1 \!\otimes h_2~|~
                           h_1 \in H_1 \,,\, h_2 \in H_2 \, \}~.
\]
\begin{remark} \label{rmk:PARTR}
 In the formulation of Lemma~\ref{lem:CNTP}, we have tacitly incorporated
 the statement that if $\pi_1$ and $\pi_2$ are faithful and traceless,
 then so is $\pi$.
 To show this, we make use of the notion of partial traces of endomorphisms
 in tensor products.
 Given two finite-dimensional complex vector spaces $V_1$ and~$V_2$, we can
 use the canonical isomorphism
 \[
  \mathrm{End}(V_1 \otimes V_2)~
  \cong~\mathrm{End}(V_1) \otimes \mathrm{End}(V_2)
 \]
 to argue that the bilinear maps
 \[
  \begin{array}{ccc}
   \mathrm{End}(V_1) \times \mathrm{End}(V_2)
   & \longrightarrow & \mathrm{End}(V_1) \\[1mm]
   (A_1,A_2)
   &   \longmapsto   & \tr(A_2) \, A_1
  \end{array}
 \]
 and
 \[
  \begin{array}{ccc}
   \mathrm{End}(V_1) \times \mathrm{End}(V_2)
   & \longrightarrow & \mathrm{End}(V_2) \\[1mm]
   (A_1,A_2)
   &   \longmapsto   & \tr(A_1) \, A_2
  \end{array}
 \]
 extend uniquely to linear maps
 \[
  \tr_2 : \mathrm{End}(V_1 \otimes V_2)~\longrightarrow~\mathrm{End}(V_1)
 \]
 and
 \[
  \tr_1 : \mathrm{End}(V_1 \otimes V_2)~\longrightarrow~\mathrm{End}(V_2)
 \]
 known as \emph{partial traces} and characterized by the property that
 \begin{equation}
  \tr_2^{}(A_1^{} \otimes A_2^{})~=~\tr_{V_2}(A_2^{}) \, A_1^{}~,
 \end{equation}
 and
 \begin{equation}
  \tr_1^{}(A_1^{} \otimes A_2^{})~=~\tr_{V_1}(A_1^{}) \, A_2^{}~,
 \end{equation}
 respectively.
 Obviously, the total trace is obtained by composing any one of these
 partial traces with the ordinary trace on the remaining tensor factor:
 \begin{equation}
  \tr_{V_1}(\tr_2^{}(A))~=~\tr_{V_1 \otimes V_2}(A)~=~\tr_{V_2}(\tr_1^{}(A))~.
 \end{equation}
 In particular,
 \begin{equation}
  \tr_{V_1 \otimes V_2}(A_1^{} \otimes A_2^{})~
  =~\tr_{V_1}(A_1^{}) \, \tr_{V_2}(A_2^{})~.
 \end{equation}
 As a result, we see that $A_1$ and $A_2$ can be recovered from
 $\, A_1^{} \otimes\, \mathrm{id}_{V_2} + \, \mathrm{id}_{V_1} \otimes A_2^{}$,
 provided both $A_1$ and $A_2$ are traceless, namely by the formula
 \begin{equation} \label{eq:PARTR1}
  \begin{array}{ccc}
   A_1~=~{\displaystyle \frac{1}{\dim V_2}} \; \tr_2^{}
          \bigl( A_1^{} \otimes\, \mathrm{id}_{V_2} + \,
                 \mathrm{id}_{V_1} \otimes A_2^{} \bigr) \quad
   & \mathrm{if} & \tr_{V_2}(A_2^{}) = 0~, \\[3mm]
   A_2~=~{\displaystyle \frac{1}{\dim V_1}} \; \tr_1^{}
          \bigl( A_1^{} \otimes\, \mathrm{id}_{V_2} + \,
                 \mathrm{id}_{V_1} \otimes A_2^{} \bigr) \quad
   & \mathrm{if} & \tr_{V_1}(A_1^{}) = 0~,
  \end{array}
 \end{equation}
 which implies that if $\pi_1$ and $\pi_2$ are faithful and traceless,
 then so is $\pi$.
 For later use, we also note the following formulas, valid for any
 $\, B \in \mathrm{End}(V_1 \otimes V_2) \,$ (they are obvious when $B$
 is tensor decomposable, $B = B_1 \otimes B_2$, and hence hold in general
 since both sides of each equation are linear in~$B$):
 \begin{equation} \label{eq:PARTR2}
  \begin{array}{ccc}
   \tr_2^{} \bigl( ( A_1^{} \otimes\, \mathrm{id}_{V_2} ) \, B \bigr)~
   =~A_1^{} \, \tr_2^{}(B) &,&
   \tr_2^{} \bigl( B \, ( A_1^{} \otimes\, \mathrm{id}_{V_2} ) \bigr)~
   =~\tr_2^{}(B) \, A_1^{}~. \\[3mm]
   \tr_1^{} \bigl( ( \,\mathrm{id}_{V_1} \otimes A_2^{} ) \, B \bigr)~
   =~A_2^{} \, \tr_1^{}(B) &,&
   \tr_1^{} \bigl( B \, ( \,\mathrm{id}_{V_1} \otimes A_2^{} ) \bigr)~
   =~\tr_1^{}(B) \, A_2^{}~.
  \end{array}
 \end{equation}
\end{remark}
\begin{remark} \label{rmk:INPRTP}
 Another useful tool that plays an important role in the proof of Lemma~%
 \ref{lem:CNTP} are the injections and projections that relate the tensor
 product space $\, V = V_1 \otimes V_2 \,$ to its factors $V_1$ and~$V_2$.
 Namely, fixed vectors $v_2$ in~$V_2$ and $v_1$ in~$V_1$ induce inclusions
 \[
  \begin{array}{cccc}
   \mathrm{i}_{v_2}: & V_1^{} & \longrightarrow &  V_1^{} \otimes V_2^{} \\[1mm]
                     & v_1^{} &   \longmapsto   & v_1^{} \!\otimes v_2^{}
  \end{array}
  \qquad \mbox{and} \qquad
  \begin{array}{cccc}
   \mathrm{i}_{v_1}: & V_2^{} & \longrightarrow &  V_1^{} \otimes V_2^{} \\[1mm]
                     & v_2^{} &   \longmapsto   & v_1^{} \!\otimes v_2^{}
  \end{array}
 \]
 whereas fixed linear forms $\, v_2^* \in V_2^* \,$ on~$V_2^{}$ and
 $\, v_1^* \in V_1^* \,$ on~$V_1^{}$ induce projections
 \[
  \begin{array}{cccc}
   \mathrm{pr}_{v_2^*}: & V_1^{} \otimes V_2^{}
   & \longrightarrow & V_1^{} \\[1mm]
                        & v_1^{} \!\otimes v_2^{}
   &   \longmapsto   & \langle v_2^*,v_2^{} \rangle \, v_1^{} 
  \end{array}
  \qquad \mbox{and} \qquad
  \begin{array}{cccc}
   \mathrm{pr}_{v_1^*}: & V_1^{} \otimes V_2^{}
   & \longrightarrow & V_2^{} \\[1mm]
                        & v_1^{} \!\otimes v_2^{}
   &   \longmapsto   & \langle v_1^*,v_1^{} \rangle \, v_2^{} 
  \end{array}
 \]
 which are equivariant under the respective actions of~$\mathfrak{g}_1$
 on~$V_1$ and~$V$ and of~$\mathfrak{g}_2$ on~$V_2$ and~$V$.
 This implies%
 \footnote{Note that the space $\, V = V_1 \otimes V_2 \,$ is
 irreducible under $\, \mathfrak{g} = \mathfrak{g}_1 \!\times
 \mathfrak{g}_2 \,$ but of course not under $\mathfrak{g}_1$ or~%
 $\mathfrak{g}_2$; rather, it is the direct sum of $\dim V_2$ copies
 of irreducible representations of~$\mathfrak{g}_1$, all equivalent to
 $\pi_1$, and also the direct sum of $\dim V_1$ copies of irreducible
 representations of~$\mathfrak{g}_2$, all equivalent to $\pi_2$, each
 of the corresponding decompositions into irreducible subspaces being
 highly non-unique.}
 that the $\mathfrak{g}_1$-irreducible subspaces of~$V$ are precisely
 the subspaces of the form $\, \mathrm{im} \; \mathrm{i}_{v_2} = V_1
 \otimes v_2 \,$ with $\, v_2 \in V_2 \setminus \{0\}$, where
 \[
  \begin{array}{ccc}
   V_1 \otimes v_2~=~V_1 \otimes v_2' &
   & \mbox{if $v_2$ and $v_2'$ are linearly dependent} \\[1mm]
   (V_1 \otimes v_2) \cap (V_1 \otimes v_2')~=~\{0\} &
   & \mbox{if $v_2$ and $v_2'$ are linearly independent} \\[1mm]
  \end{array}
 \]
 and similarly that the $\mathfrak{g}_2$-irreducible subspaces of~$V$ are
 precisely the subspaces of the form $\, \mathrm{im} \; \mathrm{i}_{v_1}
 = v_1 \otimes V_2 \,$ with $\, v_1 \in V_1 \setminus \{0\}$, where
 \[
  \begin{array}{ccc}
   v_1 \otimes V_2~=~v_1' \otimes V_2 &
   & \mbox{if $v_1$ and $v_1'$ are linearly dependent} \\[1mm]
   (v_1 \otimes V_2) \cap (v_1' \otimes V_2)~=~\{0\} &
   & \mbox{if $v_1$ and $v_1'$ are linearly independent} \\[1mm]
  \end{array}
 \]
 Of course, the only non-trivial statement here is that \emph{all} irreducible
 subspaces are of this form; here is a proof for the first case (the second one
 being completely analogous):
 Let $W$ be a $\mathfrak{g}_1$-irreducible subspace of~$V$, $W \neq \{0\}$.
 For any $\, v_2^* \in V_2^*$, the projection $\mathrm{pr}_{v_2^*}$, being
 $\mathfrak{g}_1$-equivariant, maps the $\mathfrak{g}_1$-invariant subspace
 $W$ of~$V$ to a $\mathfrak{g}_1$-invariant subspace of~$V_1$, and $V_1$
 being $\mathfrak{g}_1$-irreducible, this subspace can only be $\{0\}$
 or all of~$V_1$.
 Similarly, the kernel of~$\mathrm{pr}_{v_2^*}$ intersects~$W$
 in a $\mathfrak{g}_1$-invariant subspace of~$W$, and $W$ being
 $\mathfrak{g}_1$-irreducible, this subspace can only be $\{0\}$
 or all of~$W$.
 This means that the restriction $\mathrm{pr}_{v_2^*}\vert_W^{}$
 of~$\mathrm{pr}_{v_2^*}$ to~$W$ is either $0$ or an isomorphism.
 It cannot be zero for all $\, v_2^* \in V_2^* \,$ because there
 is \emph{no} non-zero vector in~$V$ that is annihilated by all the
 projections $\mathrm{pr}_{v_2^*}$, $v_2^* \in V_2^*$ (this is easy
 to see by expanding in a basis).
 Thus we have shown existence of a $\mathfrak{g}_1$-equivariant iso%
 morphism from~$W$ to~$V_1$; let us, for the time being, call~it~$\phi$.
 But then, for any $\, v_2^* \in V_2^*$, $\, \mathrm{pr}_{v_2^*} \circ\,
 \phi^{-1} \in \mathrm{End}(V_1) \,$ centralizes $\mathfrak{g}_1$ and,
 $V_1$ being $\mathfrak{g}_1$-irreducible, must, according to Schur's
 lemma, be a multiple of the identity: $\, \mathrm{pr}_{v_2^*} \circ\,
 \phi^{-1} = \lambda(v_2^*) \, \mathrm{id}_{V_1}$.
 Clearly, $\lambda(v_2^*)$ is linear in $v_2^*$ because so is
 $\, \mathrm{pr}_{v_2^*} \circ\, \phi^{-1}$, so there is a unique
 vector $\, v_2^{} \in V_2^{} \,$ such that $\, \lambda(v_2^*)
 = \langle v_2^*,v_2^{} \rangle \,$ and hence $\, \mathrm{pr}_{v_2^*}
 \circ\, \phi^{-1} = \langle v_2^*,v_2^{} \rangle \, \mathrm{id}_{V_1}
 = \, \mathrm{pr}_{v_2^*} \circ\, \mathrm{i}_{v_2}$.
 Since this equation holds for all $\, v_2^* \in V_2^*$, it follows
 that $\, \phi^{-1} = \, \mathrm{i}_{v_2}$, once again because there
 is \emph{no} non-zero vector in~$V$ that is annihilated by all the
 projections $\mathrm{pr}_{v_2^*}$, $v_2^* \in V_2^*$, q.e.d..
\end{remark}
\begin{proof} \hspace*{-2em}
 \footnote{A rigorous proof of this lemma does not seem to be readily
 available in the literature and will therefore be given here.}
 As in the previous lemma, the proof of the inclusions ``$\supset$'' in
 equations (\ref{eq:CENTTP1})-(\ref{eq:NORMTP3}) is straightforward provided
 we properly specify the action of~$S_2$ involved in the definition of the
 semidirect product in equation~(\ref{eq:NORMTP3}) which, roughly speaking,
 consists in switching the tensor factors.
 More precisely, choosing the pair $(\phi_{21},g_{21})$ as in Definition~%
\ref{def:EQMAUT} and setting $\, \phi_{12}^{} = \phi_{21}^{-1}$, $g_{12}^{}
 = g_{21}^{-1}$, we define an involutive automorphism $\, \phi \,$
 on~$\, \mathfrak{g} = \mathfrak{g}_1 \!\times \mathfrak{g}_2 \,$ and
 a linear involution $\, g \,$ on~$\, V = V_1 \otimes V_2$ \linebreak
 by setting
 \begin{equation} \label{eq:SWITCHTP}
  \begin{array}{cccc}
   \phi: & \mathfrak{g}_1 \!\times \mathfrak{g}_2 & \longrightarrow &
           \mathfrak{g}_1 \!\times \mathfrak{g}_2 \\[1mm]
         &                (X_1,X_2)               &   \longmapsto   &
             ( \phi_{12}(X_2) , \phi_{21}(X_1) )   \\[3mm]
   g: &  V_1 \otimes V_2  & \longrightarrow &        V_1 \otimes V_2
   \\[1mm]
      & v_1 \otimes\, v_2 &   \longmapsto   & g_{12} v_2 \otimes\, g_{21} v_1
  \end{array}
 \end{equation}
 and see that equation~(\ref{eq:EQMAUT}) becomes equivalent to the
 condition that $g$ normalizes $\mathfrak{g}$, since
 \begin{eqnarray*}
  (\phi(X) \, g)(v_1 \otimes\, v_2) \!\!
  &=&\!\! (\phi_{12}(X_2) \, g_{12})(v_2) \otimes\, g_{21}(v_1) \, + \,
          g_{12}(v_2) \otimes\, (\phi_{21}(X_1) \, g_{21})(v_1)~, \\[1mm]
  (g \, X)(v_1 \otimes\, v_2) \!\!
  &=&\!\! g_{12}(v_2) \otimes\, (g_{21} \, X_1)(v_1) \, + \,
          (g_{12} \, X_2)(v_2) \otimes\, g_{21}(v_1)~.
 \end{eqnarray*}
 To prove the converse inclusions ``$\subset$'', suppose first that
 $\, g \in GL(V) \,$ normalizes $\mathfrak{g}$ and write $\phi$ for
 the automorphism of~$\mathfrak{g}$ induced by conjugation with~$g$:
 \[
  g X g^{-1}~=~\phi(X) \quad \mbox{for $\, X \in \mathfrak{g}$}~.
 \]
 Explicitly, this means that if we write $\phi$ in the form
 \begin{equation} \label{eq:BLOCK4}
  \phi \bigl( X_1 , X_2 \bigr)~
   =~\bigl( \phi_{11}(X_1) + \phi_{12}(X_2) \,,\,
            \phi_{21}(X_1) + \phi_{22}(X_2) \bigr)
 \end{equation}
 with $\, \phi_{11} \in \mathrm{End}(\mathfrak{g}_1)$, $\phi_{12} \in
 L(\mathfrak{g}_2,\mathfrak{g}_1)$, $\phi_{21} \in L(\mathfrak{g}_1,
 \mathfrak{g}_2)$, $\phi_{22} \in \mathrm{End}(\mathfrak{g}_2)$, as before,
 then for all $\, X = (X_1,X_2) \,$ in $\, \mathfrak{g} = \mathfrak{g}_1
 \!\times \mathfrak{g}_2$, we have
 \begin{eqnarray*}
 \lefteqn{g \, \bigl( X_1^{} \otimes\, \mathrm{id}_{V_2} + \,
               \mathrm{id}_{V_1} \otimes X_2^{} \bigr) \, g^{-1}}
 \hspace{5mm} \\[1mm]
 &=&\!\! \bigl( \phi_{11}(X_1) + \phi_{12}(X_2) \bigr)
         \otimes\, \mathrm{id}_{V_2} + \,
         \mathrm{id}_{V_1} \otimes 
         \bigl( \phi_{21}(X_1) + \phi_{22}(X_2) \bigr)~,
 \end{eqnarray*}
 leading us, according to equation~(\ref{eq:PARTR1}), to the following system
 of equations:
 \begin{equation} \label{eq:BLOCK5}
  \begin{array}{c}
   \phi_{11}(X_1)~=~{\displaystyle \frac{1}{\dim V_2}} \; \tr_2^{}
   \bigl( \, g ( X_1^{} \otimes\, \mathrm{id}_{V_2} ) \, g^{-1} \bigr)~, \\[3mm]
   \phi_{21}(X_1)~=~{\displaystyle \frac{1}{\dim V_1}} \; \tr_1^{}
   \bigl( \, g ( X_1^{} \otimes\, \mathrm{id}_{V_2} ) \, g^{-1} \bigr)~, \\[3mm]
   \phi_{12}(X_2)~=~{\displaystyle \frac{1}{\dim V_2}} \; \tr_2^{}
   \bigl( \, g ( \,\mathrm{id}_{V_1} \otimes X_2^{} ) \, g^{-1} \bigr)~, \\[3mm]
   \phi_{22}(X_2)~=~{\displaystyle \frac{1}{\dim V_1}} \; \tr_1^{}
   \bigl( \, g ( \,\mathrm{id}_{V_1} \otimes X_2^{} ) \, g^{-1} \bigr)~.
  \end{array}
 \end{equation}
 Moreover, it is a standard fact (which we have already used extensively)
 that auto\-morphisms of semisimple Lie algebras permute their simple ideals,
 without mixing them, so under the additional hypotheses of equations~%
 (\ref{eq:NORMTP2}) and~(\ref{eq:NORMTP3}) ($\mathfrak{g} = \mathfrak{g}_1
 \times \mathfrak{g}_2$ with $\mathfrak{g}_1$ and $\mathfrak{g}_2$ simple),
 there are precisely two distinct possibilities: either (a)~$\phi$
 preserves both $\mathfrak{g}_1$ and~$\mathfrak{g}_2$, i.e., $\phi_{11}
 \in \mathrm{Aut}(\mathfrak{g}_1) \,$ and $\, \phi_{22} \in \mathrm{Aut}%
 (\mathfrak{g}_2)$, whereas $\, \phi_{21} \equiv 0$ \linebreak and
 $\, \phi_{12} \equiv 0$, or (b)~$\phi$ switches $\mathfrak{g}_1$
 and~$\mathfrak{g}_2$, i.e., $\phi_{11} \equiv 0 \,$ and $\, \phi_{22}
 \equiv 0$, whereas $\, \phi_{21}: \mathfrak{g}_1 \longrightarrow
 \mathfrak{g}_2 \,$ and $\, \phi_{12}: \mathfrak{g}_2 \longrightarrow
 \mathfrak{g}_1 \,$ are Lie algebra isomorphisms.
 \\[2mm]
 Now assume that $\, g \in GL(V) \,$ normalizes $\mathfrak{g}_1$ (but without
 imposing the additional hypothesis that $\mathfrak{g}_1$ and $\mathfrak{g}_2$
 should be simple): then the above calculations can be applied if we take
 $\, X_2 = 0 \,$ with $\, \phi_{11} \in \mathrm{Aut}(\mathfrak{g}_1)$,
 $\phi_{21} \equiv 0 \,$ and $\phi_{12}$, $\phi_{22}$ undetermined.
 Here, the main issue is to show that $g$ must be tensor decomposable.
 To~this~end, observe first of all that for each $\, X_1^{} \in
 \mathfrak{g}_1^{}$, there exists $\, X_1^\prime \in \mathfrak{g}_1^{} \,$
 such that $\; g ( X_1^{} \otimes\, \mathrm{id}_{V_2} ) \, g^{-1}
 =  X_1^\prime \otimes\, \mathrm{id}_{V_2} \,$; then applying $\tr_2^{}$
 to this equation and using equation~(\ref{eq:BLOCK5}), we get
 $\, X_1^\prime = \phi_{11}^{}(X_1^{})$, i.e.,
 \[
  g ( X_1^{} \otimes\, \mathrm{id}_{V_2} ) \, g^{-1}~
  =~\phi_{11}^{}(X_1^{}) \otimes\, \mathrm{id}_{V_2}~.
 \]
 Next, for arbitrary vectors $\, v_2^{} \in V_2^{} \,$ and covectors
 (linear forms) $\, v_2^* \in V_2^*$, we use the inclusions and
 projections introduced in Remark~\ref{rmk:INPRTP} to define
 \[
   \mathrm{pr}_{(v_2^*,v_2^{})}(g)~
   =~\mathrm{pr}_{v_2^*} \circ g \circ\, \mathrm{i}_{v_2^{}}^{}~
   \in~\mathrm{End}(V_1^{})~.
 \]
 Then
 \begin{eqnarray*}
  \mathrm{pr}_{(v_2^*,v_2^{})}(g) \, X_1^{} \!\!
  &=&\!\! \mathrm{pr}_{v_2^*} \circ g \circ\, \mathrm{i}_{v_2}^{} \circ X_1^{}~
   =~     \mathrm{pr}_{v_2^*} \circ g ( X_1^{} \otimes\, \mathrm{id}_{V_2})
          \circ \mathrm{i}_{v_2}^{}~, \\
  &=&\!\! \mathrm{pr}_{v_2^*} \circ
          ( \phi_{11}^{}(X_1^{}) \otimes\, \mathrm{id}_{V_2} ) \,
          g \circ \mathrm{i}_{v_2}^{}~
   =~     \phi_{11}^{}(X_1^{}) \circ \mathrm{pr}_{v_2^*} \circ
          g \circ\, \mathrm{i}_{v_2}^{} \\
  &=&\!\! \phi_{11}^{}(X_1^{}) \; \mathrm{pr}_{(v_2^*,v_2^{})}(g)~.
 \end{eqnarray*}
 This implies that both the kernel and the image of each~$\mathrm{pr}_%
 {(v_2^*,v_2^{})}(g)$ are $\mathfrak{g}_1$-invariant subspaces of~$V_1^{}$,
 so irreducibility forces all of these subspaces to be either $\{0\}$ or the
 entire space~$V_1$, which means that each $\mathrm{pr}_{(v_2^*,v_2^{})}(g)$
 is either $0$ or invertible.
 They cannot all be zero because there is \emph{no} non-zero linear operator
 in~$V$ that is annihilated by all these projections (this is easy to see by
 expanding in a basis).
 Thus we have shown existence of an invertible linear transformation on~$V_1$,
 which we shall denote by $g_{11}^{}$, which is equivariant in the sense that
 \[
  g_{11}^{} X_1^{} \, g_{11}^{-1}~=~\phi_{11}^{}(X_1^{})
  \quad \mbox{for $\, X_1^{} \in \mathfrak{g}_1^{}$}~,
 \]
 i.e., $\, g_{11}^{} \in GL(V_1^{}) \,$ normalizes~$\mathfrak{g}_1^{}$
 so that conjugation with $g_{11}^{}$ induces $\phi_{11}^{}$.
 But then, for any $\, (v_2^*,v_2^{}) \in V_2^* \times V_2^{}$,
 $g_{11}^{-1} \circ\, \mathrm{pr}_{(v_2^*,v_2^{})}(g) \in
 \mathrm{End}(V_1^{}) \,$ centralizes $\mathfrak{g}_1^{}$ and,
 $V_1^{}$ being $\mathfrak{g}_1^{}$-irreducible, must, according
 to Schur's lemma, be a multiple of the identity:
 \[
  g_{11}^{-1} \circ\, \mathrm{pr}_{(v_2^*,v_2^{})}(g)~
  =~\lambda(v_2^*,v_2^{}) \, \mathrm{id}_{V_1}~.
 \]
 Clearly, $\lambda$ is bilinear because $g_{11}^{-1} \circ\, \mathrm{pr}_%
 {(v_2^*,v_2^{})}(g)$ is bilinear in $(v_2^*,v_2^{})$, so there is a unique
 linear transformation $\, g_{22}^{} \in \mathrm{End}(V_2^{}) \,$ such that
 $\, \lambda(v_2^*,v_2^{}) = \langle v_2^*,g_{22}^{} v_2^{} \rangle$,
 implying
 \[
   \mathrm{pr}_{v_2^*} \circ g \circ\, \mathrm{i}_{v_2^{}}^{}~
   =~\mathrm{pr}_{(v_2^*,v_2^{})}(g)~
   =~\langle v_2^*,g_{22}^{} v_2^{} \rangle \, g_{11}^{}~
   =~\mathrm{pr}_{v_2^*} \circ (g_{11}^{} \otimes g_{22}^{})
     \circ\, \mathrm{i}_{v_2^{}}^{}~.
 \]
 Since $v_2^*$ and $v_2^{}$ are arbitrary, this implies $\, g = g_{11}^{} \!
 \otimes g_{22}^{}$, with $\, g_{11}^{} \in N_{GL(V_1)}^{}(\mathfrak{g}_1^{})$
 \linebreak
 and $\, g_{22}^{} \in GL(V_2^{}) \,$ arbitrary: this proves the inclusion
 ``$\subset$'' in the first line of equation~(\ref{eq:NORMTP1}) and, as an
 immediate corollary, that in the first line of equation~(\ref{eq:CENTTP1})
 (the second one being completely analogous).
 \\[2mm]
 Finally, we return to the case where $\, g \in GL(V) \,$ is supposed to
 normalize $\mathfrak{g}$ and both $\mathfrak{g}_1$  and~$\mathfrak{g}_2$
 are supposed to be simple.
 If (a)~$\phi$ preserves both $\mathfrak{g}_1$ and~$\mathfrak{g}_2$,
 we can immediately apply the previous part of the proof to conclude that
 \[
  g~\in~N_{GL(V)}^{}(\mathfrak{g}_1^{}) \cap N_{GL(V)}^{}(\mathfrak{g}_2^{})~
  =~N_{GL(V_1)}^{}(\mathfrak{g}_1^{}) \otimes N_{GL(V_2)}^{}(\mathfrak{g}_2^{})~.
 \]
 Thus assume that (b)~$\phi$ switches $\mathfrak{g}_1$ and~$\mathfrak{g}_2$,
 i.e., $\phi_{11} \equiv 0 \,$ and $\, \phi_{22} \equiv 0$, whereas
 $\, \phi_{21}: \mathfrak{g}_1 \longrightarrow \mathfrak{g}_2 \,$ and
 $\, \phi_{12}: \mathfrak{g}_2 \longrightarrow \mathfrak{g}_1 \,$ are
 Lie algebra isomorphisms.
 Then we can argue in a similar manner as before, observing first
 of all that for each $\, X_1^{} \in \mathfrak{g}_1^{} \,$ and each
 $\, X_2^{} \in \mathfrak{g}_2^{}$, there exist $\, X_1^\prime \in
 \mathfrak{g}_2^{} \,$ and $\, X_2^\prime \in \mathfrak{g}_1^{} \,$
 such that $\; g ( X_1^{} \otimes\, \mathrm{id}_{V_2} ) \, g^{-1}
 = \,\mathrm{id}_{V_1} \otimes X_1^\prime \;$ and
 $\; g ( \,\mathrm{id}_{V_1} \otimes X_2^{} ) \, g^{-1}
 = X_2^\prime \otimes\, \mathrm{id}_{V_2} \,$; then applying
 $\tr_1^{}$ to the first equation and $\tr_2^{}$ to the second
 equation and using equation~(\ref{eq:BLOCK5}),  we get
 $\, X_1^\prime = \phi_{21}^{}(X_1^{}) \,$ and $\, X_2^\prime
 = \phi_{12}^{}(X_2^{})$, i.e.,
 \[
  \begin{array}{c}
   g ( X_1^{} \otimes\, \mathrm{id}_{V_2} ) \, g^{-1}~
   =~\,\mathrm{id}_{V_1} \otimes \phi_{21}^{}(X_1^{})~, \\[2mm]
   g ( \,\mathrm{id}_{V_1} \otimes X_2^{} ) \, g^{-1}~
   =~\phi_{12}^{}(X_2^{}) \otimes\, \mathrm{id}_{V_2}~.
  \end{array}
 \]
 Next, for arbitrary vectors $\, v_1^{} \in V_1^{}$, $\, v_2^{} \in V_2^{} \,$
 and covectors (linear forms) $\, v_1^* \in V_1^*$, $\, v_2^* \in V_2^*$, we
 use the inclusions and projections introduced in Remark~\ref{rmk:INPRTP}
 to define
 \[
  \begin{array}{c}
   \mathrm{pr}_{(v_1^*,v_2^{})}(g)~
   =~\mathrm{pr}_{v_1^*} \circ g \circ\, \mathrm{i}_{v_2^{}}^{}~
   \in~L(V_1^{},V_2^{})~, \\[2mm]
   \mathrm{pr}_{(v_2^*,v_1^{})}(g)~
   =~\mathrm{pr}_{v_2^*} \circ g \circ\, \mathrm{i}_{v_1^{}}^{}~
   \in~L(V_2^{},V_1^{})~.
  \end{array}
 \]
 Then
 \begin{eqnarray*}
  \mathrm{pr}_{(v_1^*,v_2^{})}(g) \circ X_1^{} \!\!
  &=&\!\! \mathrm{pr}_{v_1^*} \circ g \circ\, \mathrm{i}_{v_2}^{} \circ X_1^{}~
   =~     \mathrm{pr}_{v_1^*} \circ g ( X_1^{} \otimes\, \mathrm{id}_{V_2})
          \circ \mathrm{i}_{v_2}^{}~, \\
  &=&\!\! \mathrm{pr}_{v_1^*} \circ
          ( \,\mathrm{id}_{V_1} \otimes \phi_{21}^{}(X_1^{}) ) \,
          g \circ \mathrm{i}_{v_2}^{}~
   =~     \phi_{21}^{}(X_1^{}) \circ \mathrm{pr}_{v_1^*} \circ
          g \circ\, \mathrm{i}_{v_2}^{} \\
  &=&\!\! \phi_{21}^{}(X_1^{}) \circ\, \mathrm{pr}_{(v_1^*,v_2^{})}(g)~, \\[3mm]
  \mathrm{pr}_{(v_2^*,v_1^{})}(g) \circ X_2^{} \!\!
  &=&\!\! \mathrm{pr}_{v_2^*} \circ g \circ\, \mathrm{i}_{v_1}^{} \circ X_2^{}~
   =~     \mathrm{pr}_{v_2^*} \circ g ( \,\mathrm{id}_{V_1} \otimes X_2^{} )
          \circ \mathrm{i}_{v_1}^{}~, \\
  &=&\!\! \mathrm{pr}_{v_2^*} \circ
          ( \phi_{12}^{}(X_2^{}) \otimes\, \mathrm{id}_{V_2} ) \,
          g \circ \mathrm{i}_{v_1}^{}~
   =~     \phi_{12}^{}(X_2^{}) \circ \mathrm{pr}_{v_2^*} \circ
          g \circ\, \mathrm{i}_{v_1}^{} \\
  &=&\!\! \phi_{12}^{}(X_2^{}) \circ\, \mathrm{pr}_{(v_2^*,v_1^{})}(g)~.
 \end{eqnarray*}
 This implies that the kernel of each~$\mathrm{pr}_{(v_1^*,v_2^{})}(g)$ and
 the image of each~$\mathrm{pr}_{(v_2^*,v_1^{})}(g)$ is a $\mathfrak{g}_1$-%
 invariant subspace of~$V_1^{}$ and similarly the kernel of each~%
 $\mathrm{pr}_{(v_2^*,v_1^{})}(g)$ and the image of each~%
 $\mathrm{pr}_{(v_1^*,v_2^{})}(g)$ is a $\mathfrak{g}_2$-%
 invariant subspace of~$V_2^{}$, so irreducibility forces all
 of these subspaces to be either $\{0\}$ or the entire space,
 which means that each $\mathrm{pr}_{(v_1^*,v_2^{})}(g)$ and each~%
 $\mathrm{pr}_{(v_2^*,v_1^{})}(g)$ is either $0$ or a linear isomorphism.
 They cannot all be zero because there is \emph{no} non-zero linear
 operator in~$V$ that is annihilated by all these projections (this
 is easy to see by expanding in a basis).
 Thus we have shown existence of linear isomorphisms from~$V_1^{}$
 to~$V_2^{}$ and from $V_2^{}$ to~$V_1^{}$ which we shall denote by
 $g_{21}$ and by~$g_{12}$, respectively, and which are equivariant
 in the sense that
 \[
  \begin{array}{c}
   \phi_{21}^{}(X_1^{})~=~g_{21}^{} \circ X_1^{} \circ\, g_{21}^{-1}~, \\[2mm]
   \phi_{12}^{}(X_2^{})~=~g_{12}^{} \circ X_2^{} \circ\, g_{12}^{-1}~.
  \end{array}
 \]
 Moreover, using Schur's lemma, we see that the product of $g_{21}$ and~%
 $g_{12}$ must be a multiple of the identity, so multiplying one of them
 with an appropriate scalar factor, we may assume without loss of generality
 that they are each other's inverse.
 Thus we conclude not only that the representations $\pi_1$ of~%
 $\mathfrak{g}_1$ and $\pi_2$ of~$\mathfrak{g}_2$ are quasiequivalent
 but also that any invertible linear transformation $g$ on~$V$ that
 normalizes~$\mathfrak{g}$ and switches $\mathfrak{g}_1$ and~%
 $\mathfrak{g}_2$ provides a realization of this quasiequivalence,
 in the sense of Definition~\ref{def:EQMAUT}.
 Finally, the product of any such invertible linear transformation
 on~$V$ with the inverse of any other one is of course an invertible
 linear transformation on~$V$ that normalizes both $\mathfrak{g}_1$
 and~$\mathfrak{g}_2$ and hence, by what has been proved before,
 belongs to \mbox{$N_{GL(V_1)}^{}(\mathfrak{g}_1^{}) \otimes
 N_{GL(V_2)}^{}(\mathfrak{g}_2^{})$.}
\end{proof}
\begin{remark} \label{rmk:COUNEX}
 To see that equation~(\ref{eq:NORMTP2}) may fail to be true if
 $\mathfrak{g}_1$ and~$\mathfrak{g}_2$ are not simple (even when
 they are semisimple), consider the following situation: take
 $\mathfrak{g}_1 = \mathfrak{g}_0 \!\times \tilde{\mathfrak{g}}_1$,
 $\mathfrak{g}_2 = \mathfrak{g}_0 \!\times \tilde{\mathfrak{g}}_2$
 and $\pi_1 = \pi_0 \boxtimes \tilde{\pi}_1 \,$ on $\, V_1 = V_0
 \otimes \tilde{V}_1$, $\pi_2 = \pi_0 \boxtimes \tilde{\pi}_2$
 \linebreak
 on $\, V_2 = V_0 \otimes \tilde{V}_2$, where $\mathfrak{g}_0$,
 $\tilde{\mathfrak{g}}_1$ and $\tilde{\mathfrak{g}}_2$ are mutually
 non-isomorphic simple Lie algebras and $\pi_0$, $\tilde{\pi}_1$ and
 $\tilde{\pi}_2$ are irreducible representations of $\mathfrak{g}_0$,
 $\tilde{\mathfrak{g}}_1$ and $\tilde{\mathfrak{g}}_2$ on the vector
 spaces $V_0$, $\tilde{V}_1$ and $\tilde{V}_2$, respectively.
 Then obviously, $N_{GL(V)}(\mathfrak{g})$ contains a linear
 involution $g$ which switches the two copies of $V_0$ but
 preserves $\tilde{V}_1$ and $\tilde{V}_2$, so it neither
 preserves nor switches $V_1$ and~$V_2$.
\end{remark}
\noindent
The same methods allow us to prove the following fact about invariant
bilinear or sesquilinear forms in tensor products, which will be used
later on:
\begin{lemma} \label{lem:BFTP}
 Let $\, \pi_1: \mathfrak{g}_1 \longrightarrow \mathfrak{gl}(V_1) \,$
 and $\, \pi_2: \mathfrak{g}_2 \longrightarrow \mathfrak{gl}(V_2) \,$
 be irreducible representations of Lie algebras\/ $\mathfrak{g}_1$ and\/
 $\mathfrak{g}_2$ on finite-dimensional complex vector spaces\/ $V_1$
 and\/~$V_2$, respectively, and let $\, \pi: \mathfrak{g} \longrightarrow
 \mathfrak{gl}(V) \,$ be their external tensor product: $\, \pi = \pi_1
 \boxtimes \pi_2 \,$, $\mathfrak{g} = \mathfrak{g}_1 \!\times
 \mathfrak{g}_2 \,$, $V = V_1 \otimes V_2 \,$.
 Then every invariant bilinear or sesquilinear form $(.\,,.)$ on~$V$
 is the tensor product of invariant bilinear or sesquilinear forms
 $(.\,,.)_1$ on~$V_1$ and $(.\,,.)_2$ on~$V_2$:
 \begin{equation} \label{eq:BFTP}
  (v_1^{} \otimes v_2^{} \,, v_1^\prime \otimes v_2^\prime)~
  =~(v_1^{},v_1^\prime)_1^{} \, (v_2^{},v_2^\prime)_2^{}~.
 \end{equation}
 The generalization to more than two tensor factors is the obvious one.
\end{lemma}
\begin{proof}
 Given an invariant bilinear or sesquilinear form $(.\,,.)$ on~$V$,
 let us define, for any two vectors $\, v_2^{},v_2^\prime \in V_2^{}$,
 an invariant bilinear or sesquilinear form $(.\,,.)_{(v_2^{},v_2^\prime)}$
 on~$V_1^{}$ by
 \[
  (v_1^{},v_1^\prime)_{(v_2^{},v_2^\prime)}~
  =~(v_1^{} \otimes v_2^{} \,, v_1^\prime \otimes v_2^\prime)~.
 \]
 That this form is really invariant is easy to check, since for
 $\, X_1 \in \mathfrak{g}_1$,
 \[
  \begin{array}{l}
   \bigl( X_1^{} v_1^{},v_1^\prime \bigr)_{(v_2^{},v_2^\prime)} +
   \bigl( v_1^{},X_1^{} v_1^\prime \bigr)_{(v_2^{},v_2^\prime)} \\[3mm]
   \hspace*{5mm}
   =~\bigl( (X_1^{} \otimes \mathrm{id}_{V_2}) (v_1^{} \otimes v_2^{}) \,,
            v_1^\prime \otimes v_2^\prime \bigr) + 
     \bigl( v_1^{} \otimes v_2^{} \,,
            (X_1^{} \otimes \mathrm{id}_{V_2}) (v_1^\prime \otimes v_2^\prime)
            \bigr)~.
  \end{array}
 \]
 Therefore, its kernel (formed by the vectors orthogonal to all vectors
 in~$V_1$) is a $\mathfrak{g}_1$-invariant subspace of~$V_1^{}$, so
 irreducibility forces it to be either $\{0\}$ or the entire space~$V_1$,
 which means that each of the forms $(.\,,.)_{(v_2^{},v_2^\prime)}$ is
 either $0$ or non-degenerate.
 They cannot all be zero, except when the original form $(.\,,.)$ on~$V$
 is zero, in which case the claim is trivial.
 Otherwise, we obtain existence of a non-degenerate invariant bilinear
 or sesquilinear form $(.\,,.)_1$ on~$V_1$.
 However, due to irreducibility, any two such forms are proportional,
 so we conclude that, for any $\, v_2^{},v_2^\prime \in V_2^{}$,
 $(.\,,.)_{(v_2^{},v_2^\prime)} = \lambda(v_2^{},v_2^\prime) (.\,,.)_1^{}$.
 Clearly, $\lambda$ is bilinear or sequilinear because $(.\,,.)_{(v_2^{},
 v_2^\prime)}$ is bilinear or sesquilinear in $(v_2^{},v_2^\prime)$,
 so we are done.
\end{proof}

With these concepts at our disposal, assume now that $\mathfrak{s}\,$ is a
reducible subalgebra of one of the classical Lie algebras $\mathfrak{g}$;
note that we do not exclude the case that $\mathfrak{s}\,$ is simple.
Viewing the inclusion of~$\mathfrak{s}\,$ into~$\mathfrak{g}$ as a
(faithful) representation $\pi$ of~$\mathfrak{s}\,$ on $\, V \cong
\mathbb{C}^n \,$ and using the fact that representations of reductive
Lie algebras are completely reducible, we may decompose $\pi$ into its
irreducible constituents, which are certain irreducible representations
$\pi_i$ of~$\mathfrak{s}\,$ on $\, V_i \cong \mathbb{C}^{n_i}$
($1 \leqslant i \leqslant r$), as follows:
\begin{equation} \label{eq:ISOTDEC}
 V~=~\bigoplus_{i=1}^r \; V_i \otimes W_i \qquad,\qquad
 \pi(X)~=~\bigoplus_{i=1}^r \; \pi_i(X) \otimes \mathrm{id}_{W_i}
 \quad \mbox{for $\, X \in \mathfrak{s}$}~.
\end{equation}
Here, the $\, W_i \cong \mathbb{C}^{m_i}$ ($1 \leqslant i \leqslant r$)
are auxiliary spaces on which $\mathfrak{s}$ acts trivially: for each~$i$,
$m_i$ is the \emph{multiplicity} with which $\pi_i$ appears in~$\pi$, while
the subspace $\, V_i \otimes W_i \,$ of~$V$ is known as the corresponding
\emph{isotypic component}.
The main advantage of this decomposition into isotypic components is that
it is unique, whereas the finer one into irreducible subspaces is not.
Another useful fact is that every irreducible subspace of~$V$ is contained
in precisely one isotypic component.
Finally, the representation $\pi$ is said to be \emph{multiplicity free}
if all $m_i$ are equal to $1$: in this case, the isotypic components become
identical with the irreducible subspaces, and these become unique.
\begin{lemma}
 Let\/ $\mathfrak{s}$ be a\/ $\Sigma$-primitive reducible subalgebra of
 one of the classical Lie algebras\/ $\mathfrak{g}$, where\/ $\Sigma$ is
 any subgroup of\/~$\mathrm{Out(\mathfrak{g})}$. Then the repesentation
 of\/~$\mathfrak{s}$ defined by its inclusion in\/~$\mathfrak{g}$ is
 multiplicity free.
\end{lemma}
\begin{proof}
 This is a simple consequence of the fact that, in general, the isotypic
 decomposition (\ref{eq:ISOTDEC}) implies that linear transformations
 $\, A \in \mathrm{End}(V) \,$ of the form
 \[
  A~=~\bigoplus_{i=1}^r \; \mathrm{id}_{V_i} \otimes A_i
 \]
 with $\, A_i \in \mathrm{End}(W_i) \,$ commute with all $\pi(X)$, $X \in
 \mathfrak{s}$, and hence if $\, m_i>1 \,$ for some~$i$, the centralizer
 of~$\mathfrak{s}$ in~$\mathfrak{g}$ contains elements that do no belong
 to~$\mathfrak{s}$, so $\mathfrak{s}$ cannot be self-normalizing.
\end{proof}
\noindent
As a result, it becomes clear that for complex representations, that is,
when $\, \mathfrak{g} = \mathfrak{su}(n)$, and with an appropriate choice
of basis, $\mathfrak{s}$ is contained in the \emph{complex flag algebra}
\begin{equation} \label{eq:CFLAG}
 \mathfrak{f}~
 =~\mathfrak{s} \bigl( \mathfrak{u}(n_1) \oplus \ldots \oplus
                       \mathfrak{u}(n_r) \bigr) \quad
 (n = n_1 + \ldots + n_r)~.
\end{equation}
For self-conjugate representations, that is, when $\, \mathfrak{g}
= \mathfrak{so}(n) \,$ or $\, \mathfrak{g} = \mathfrak{sp}(n)$, we can
say somewhat more: since $\mathfrak{g}$ is pointwise fixed under conjugation
with~$\tau$ and hence so is $\mathfrak{s}$, $\tau$~permutes the afore%
mentioned irreducible subspaces and therefore these fall into two classes:
single subspaces that are $\tau$-invariant and pairs of subspaces that are
switched under the action of $\tau$.
As a result, again with an appropriate choice of basis, $\mathfrak{s}$
is contained in the \emph{generalized real flag algebra}
\begin{equation} \label{eq:RFLAG}
 \begin{array}{c}
  \mathfrak{f}~
  =~\mathfrak{so}(n_1) \oplus \ldots \oplus \mathfrak{so}(n_r) \oplus
    \mathfrak{u}(n_{r+1}) \oplus \ldots \oplus \mathfrak{u}(n_{r+s}) \\[1mm]
  (n = n_1 + \ldots + n_r + 2n_{r+1} + \ldots + 2n_{r+s})
 \end{array}~,
\end{equation}
when $\, \mathfrak{g} = \mathfrak{so}(n)$ (using the inclusion
$\, \mathfrak{u}(n_i) \subset \mathfrak{so}(2n_i)$), but is
contained in the \emph{generalized quaternionic flag algebra}
\begin{equation} \label{eq:QFLAG}
 \begin{array}{c}
  \mathfrak{f}~
  =~\mathfrak{sp}(n_1) \oplus \ldots \oplus \mathfrak{sp}(n_r) \oplus
    \mathfrak{u}(n_{r+1}) \oplus \ldots \oplus \mathfrak{u}(n_{r+s}) \\[1mm]
  (n = n_1 + \ldots + n_r + 2n_{r+1} + \ldots + 2n_{r+s})
 \end{array}~,
\end{equation}
when $\, \mathfrak{g} = \mathfrak{sp}(n)$ (using the inclusion
$\, \mathfrak{u}(n_i) \subset \mathfrak{sp}(2n_i)$).
Moreover, all these (generalized) flag algebras $\mathfrak{f}$
are $\mathrm{Inn}(\mathfrak{g},\mathfrak{s})$-invariant, and
when $\mathfrak{s}\,$ is $\mathbb{Z}_2$-invariant, they are also
$\mathrm{Aut}_{\mathbb{Z}_2}(\mathfrak{g},\mathfrak{s})$-invariant.%
\footnote{From now on, whenever we speak about $\mathbb{Z}_2$-invariant
or $\mathrm{Aut}_{\mathbb{Z}_2}(\mathfrak{g},\mathfrak{h})$-invariant
or $\mathbb{Z}_2$-primitive subalgebras of one of the classical Lie
algebras~$\mathfrak{g}$, we always assume that $\, \mathfrak{g}
= \mathfrak{su}(n) \,$ or $\, \mathfrak{g} = \mathfrak{so}(n) \,$
with $n$ even, even when this is not stated explicitly.}
[\,Indeed, given $\, \phi \in \mathrm{Inn}(\mathfrak{g},\mathfrak{s})$,
write it in the form $\, \phi = \mathrm{Ad}(g) \,$ with $\, g \in G_0$.
Then since $\phi\,$ preserves $\mathfrak{s}$, $g$ permutes the afore%
mentioned irreducible subspaces; moreover, when $\, \mathfrak{g} =
\mathfrak{so}(n) \,$ or $\, \mathfrak{g} = \mathfrak{sp}(n)$, $g$
commutes with $\tau$ and hence transforms irreducible subspaces
in the same pair to irreducible subspaces in the same pair.
But then $\, \phi = \mathrm{Ad}(g) \,$ preserves~$\mathfrak{f}$,
since $\mathfrak{f}$ consists precisely of \emph{all} elements
of~$\mathfrak{g}$ that map each of the aforementioned irreducible
subspaces to itself.
The same argument prevails when $\mathfrak{s}\,$ is $\mathbb{Z}_2$-%
invariant and  $\, \phi \in \mathrm{Aut}_{\mathbb{Z}_2}(\mathfrak{g},%
\mathfrak{s}) \setminus \mathrm{Inn}(\mathfrak{g},\mathfrak{s})$,
since we can still write $\, \phi = \mathrm{Ad}(g) \,$ where $g$ is
now an antiunitary transformation on~$\mathbb{C}^n$ in the case of~%
$\mathfrak{su}(n)$ and an orthogonal transformation on~$\mathbb{R}^n$
of determinant $-1$ in the case of $\mathfrak{so}(n)$: the hypothesis
that $\mathfrak{s}$ should be $\mathbb{Z}_2$-invariant merely states
that the set of such automorphisms $\phi$ is not empty.\,]
But this implies that if $\mathfrak{s}$ is assumed to be primitive,
it must in fact be equal to $\mathfrak{f}$, and the same conclusion
holds if $\mathfrak{s}$ is assumed to be $\mathbb{Z}_2$-primitive,
provided $\mathfrak{f}$ is $\mathbb{Z}_2$-invariant. Note also that
for $\, \mathfrak{g} = \mathfrak{su}(n)$, this additional condition
of $\mathbb{Z}_2$-invariance is automatically satisfied by the
complex flag algebras of equation~(\ref{eq:CFLAG}), but that
for $\, \mathfrak{g} = \mathfrak{so}(n)$ with $n$ even, it
imposes restrictions on the generalized real flag algebras
of equation~(\ref{eq:RFLAG}) that are allowed; these will
be derived below.

Having shown that the reducible primitive or $\mathbb{Z}_2$-primitive
subalgebras of classical Lie algebras must be sought among
the (generalized) flag algebras or $\mathbb{Z}_2$-invariant
(generalized) flag algebras $\mathfrak{f}$ introduced above, our next
step will be to figure out the additional constraints on the numbers
$n_i$ imposed by primitivity or $\mathbb{Z}_2$-primitivity.
Beginning with the primitive case, we first note that any automorphism
of~$\mathfrak{su}(n)$, $\mathfrak{so}(n)$ or~$\mathfrak{sp}(n)$ which
preserves a (generalized) flag algebra as defined in equations~%
(\ref{eq:CFLAG}), (\ref{eq:RFLAG}) or~(\ref{eq:QFLAG}), respectively,
can only permute blocks of equal type and equal size.
As a result, there are various situations in which a (generalized) flag
algebra $\mathfrak{f}$ cannot be primitive, simply because
one easily finds a strictly larger proper subalgebra $\tilde{\mathfrak{f}}$
of~$\mathfrak{g}$ which is $\mathrm{Inn}(\mathfrak{g},\mathfrak{f})$-%
invariant, namely:
\begin{itemize}
 \item If, in the case of $\mathfrak{so}(n)$ or~$\mathfrak{sp}(n)$,
       $\mathfrak{f}$ contains ``orthogonal blocks'' or ``symplectic
       blocks'', respectively, together with ``unitary blocks'',
       take $\tilde{\mathfrak{f}}$ to be $\, \mathfrak{so}(p)
       \oplus \mathfrak{so}(2q) \,$ or $\, \mathfrak{sp}(p) \oplus
       \mathfrak{sp}(2q)$, with $\, p = n_1 + \ldots + n_r \,$ and
       \linebreak $\, q = n_{r+1} + \ldots + n_{r+s} \,$, respectively,
       to conclude that $\mathfrak{f}$ cannot be primitive.
 \item If, in the case of $\mathfrak{so}(n)$ or~$\mathfrak{sp}(n)$,
       $\mathfrak{f}$ contains no ``orthogonal blocks'' or ``symplectic
       blocks'', respectively, so then $\, r=0$, but contains more
       than one ``unitary block'', take $\tilde{\mathfrak{f}}$ to be
       $\, \mathfrak{so}(2n_1) \oplus \ldots \oplus \mathfrak{so}(2n_s) \,$
       or \linebreak $\, \mathfrak{sp}(2n_1) \oplus \ldots \oplus
       \mathfrak{sp}(2n_s)$, respectively, to conclude that
       $\mathfrak{f}$ cannot be primitive.
 \item If, in the case of $\mathfrak{su}(n)$ or in the case
       of~$\mathfrak{so}(n)$ with ``orthogonal blocks'' only
       or of~$\mathfrak{sp}(n)$ with ``symplectic blocks'' only,
       $\mathfrak{f}$ contains more than two blocks whose sizes are
       not all equal, and arranging the blocks according to their
       size, say in decreasing order, with $i$ denoting the first
       index for which $\, n_i > n_{i+1}$, take $\tilde{\mathfrak{f}}$
       to be $\, \mathfrak{s}(\mathfrak{u}(p) \oplus \mathfrak{u}(q)) \,$
       or $\, \mathfrak{so}(p) \oplus \mathfrak{so}(q) \,$ or
       $\, \mathfrak{sp}(p) \oplus \mathfrak{sp}(q)$, respectively,
       with $\, p = n_1 + \ldots + n_i \,$ and $\, q = n_{i+1}
       + \ldots + n_r \,$, to conclude that  $\mathfrak{f}$
       cannot be primitive.
\end{itemize}
Thus we are left with the following candidates for reducible primitive
subalgebras:
\[
 \begin{array}{ccl}
  \mathfrak{g}~=~\mathfrak{su}(n): &\quad&
  \left\{ \begin{array}{cc}
           \mathfrak{s} \bigl( \mathfrak{u}(p) \oplus \mathfrak{u}(q) \bigr) &
           (n = p+q) \\
           \mathfrak{s} \bigl( \mathfrak{u}(p) \oplus \ldots \oplus
                               \mathfrak{u}(p) \bigr)
           & (\mbox{$l$ summands, } n = p\>\!l\,,\, l \geqslant 3)
          \end{array} \right. \\[3ex]
  \mathfrak{g}~=~\mathfrak{so}(n): &\quad&
  \left\{ \begin{array}{cc}
           \mathfrak{so}(p) \oplus \mathfrak{so}(q) & (n = p+q) \\
           \mathfrak{so}(p) \oplus \ldots \oplus \mathfrak{so}(p)
           & (\mbox{$l$ summands, } n = p\>\!l\,,\, l \geqslant 3) \\
           \mathfrak{u}(p) & (n = 2p)
          \end{array} \right. \\[4ex]
  \mathfrak{g}~=~\mathfrak{sp}(n): &\quad&
  \left\{ \begin{array}{cc}
           \mathfrak{sp}(2p) \oplus \mathfrak{sp}(2q) & (n = 2(p+q)) \\
           \mathfrak{sp}(2p) \oplus \ldots \oplus \mathfrak{sp}(2p)
           & (\mbox{$l$ summands, } n = 2p\>\!l\,,\, l \geqslant 3) \\
           \mathfrak{u}(p) & (n = 2p)
          \end{array} \right.
 \end{array}
\]
Finally, in the case of $\mathfrak{so}(n)$ with $n$ even, we must check which
of these are also $\mathbb{Z}_2$-invariant and hence $\mathbb{Z}_2$-primitive.
Now if $\mathfrak{f}$ contains at least one ``orthogonal block'' (even a
trivial one of the form $\, \mathfrak{so}(1) = \{0\} \,$ will do), then
$\mathfrak{f}$ is $\mathbb{Z}_2$-invariant, because we can define the
orthogonal transformation $g$ on~$\mathbb{R}^n$ of determinant $-1$ that
implements the outer automorphism of $\mathfrak{so}(n)$ to be orthogonal
of determinant $-1$ in that block and the identity in all other blocks.
But if $\mathfrak{f}$ consists just of a single ``unitary block''
($\mathfrak{f} = \mathfrak{u}(p) \,$ with $\, n=2p$), then $g$ is
complex conjugation in $\, \mathbb{C}^p = \mathbb{R}^n$, and this has
determinant $(-1)^p$, so $\mathfrak{f}$ will be $\mathbb{Z}_2$-invariant
if and only if $p$ is odd.

Looking at the result, we note that it includes the first statement of
\linebreak Theorem~\ref{theo:DYNEXC} above: the only reducible primitive
or $\mathbb{Z}_2$-primitive subalgebras of classical Lie algebras which
are simple are given by the inclusion~(\ref{eq:RSMPSA}) (the above case
of $\, \mathfrak{so}(p) \oplus \mathfrak{so}(q) \,$ with $\, p=n-1 \,$
and $\, q=1$, say, taking into account that
$\, \mathfrak{so}(1) = \{0\}$).

Passing to the case of irreducible subalgebras, we note first of all
that any irreducible subalgebra $\mathfrak{s}$ of a classical Lie
algebra $\mathfrak{g}$ is necessarily semisimple since its center,
consisting of multiples of the identity and at the same time belonging
to~$\mathfrak{g}$ (whose elements are traceless), reduces~to~$\{0\}$.
Thus if $\mathfrak{s}$ is not simple, we may decompose the irreducible
representation space $\, V \cong \mathbb{C}^n \,$ of~$\mathfrak{s}$
into the tensor product $\, V = V_1 \otimes \ldots \otimes V_r \,$ of
irreducible representation spaces $\, V_i \cong \mathbb{C}^{n_i} \,$
of the simple ideals $\mathfrak{s}_i$ of~$\mathfrak{s}$ ($1 \leqslant
i \leqslant r$), where $\, n = n_1 \ldots n_r \,$ (see, e.g., the first
part of Theorem~3.3 in~\cite[p.~272]{Dy2} or Proposition~3.1.8 in~%
\cite[p.~123]{GW}).
Given the fact that irreducible representations map simple Lie algebras
into Lie algebras of traceless matrices and applying Lemma~\ref{lem:BFTP}
(for the case of sesquilinear forms), it becomes clear that for complex
representations, i.e., when $\, \mathfrak{g} = \mathfrak{su}(n)$, and
with an appropriate choice of basis, $\mathfrak{s}$ is contained in
\begin{equation} \label{eq:CALGTP}
 \mathfrak{f}^\times~
 =~\mathfrak{su}(n_1) \times \ldots \times \mathfrak{su}(n_r) \quad
 (n = n_1 \ldots n_r)~.
\end{equation}
For self-conjugate representations, that is, when $\, \mathfrak{g}
= \mathfrak{so}(n) \,$ or $\, \mathfrak{g} = \mathfrak{sp}(n)$, we
can say somewhat more: since in these cases, the representation
of~$\mathfrak{s}$ on~$\mathbb{C}^n$ preserves a bilinear form
(which is symmetric when $\, \mathfrak{g} = \mathfrak{so}(n) \,$
and antisymmetric when $\, \mathfrak{g} = \mathfrak{sp}(n)$),
Lemma~\ref{lem:BFTP} (for the case of bilinear forms) implies that
the same is true for each of the representations of $\mathfrak{s}_i$
on~$\mathbb{C}^{n_i}$; moreover, since the tensor product of two
symmetric or two antisymmetric bilinear forms is symmetric whereas
that of a symmetric and an antisymmetric bilinear form is anti%
symmetric, we conclude that $\mathfrak{s}$ is contained in
\begin{equation} \label{eq:RALGTP}
 \begin{array}{c}
  \mathfrak{f}^\times~
  =~\mathfrak{so}(n_1) \times \ldots \times \mathfrak{so}(n_r) \times
    \mathfrak{sp}(n_{r+1}) \times \ldots \times \mathfrak{sp}(n_{r+s}) \\[1mm]
  (n = n_1 \ldots n_r \, n_{r+1} \ldots n_{r+s})
 \end{array}~,
\end{equation}
with $s$ even when $\, \mathfrak{g} = \mathfrak{so}(n)$ and $s$ odd
when $\, \mathfrak{g} = \mathfrak{sp}(n) \,$ (see, e.g., the second
part of Theorem~3.3 in~\cite[p.~272]{Dy2}).
Moreover, once again, all these subalgebras $\mathfrak{f}^\times$
are $\mathrm{Inn}(\mathfrak{g},\mathfrak{s})$-invariant, and when
$\mathfrak{s}\,$ is $\mathbb{Z}_2$-invariant, they are also
$\mathrm{Aut}_{\mathbb{Z}_2}(\mathfrak{g},\mathfrak{s})$-invariant.%
\addtocounter{footnote}{-1}\footnotemark\
[\,Indeed, given $\, \phi \in \mathrm{Inn}(\mathfrak{g},\mathfrak{s})$,
write it in the form $\, \phi = \mathrm{Ad}(g) \,$  with $\, g \in G_0$.
Then since $\phi\,$ preserves $\mathfrak{s}$, $\phi\,$ permutes the
simple ideals of~$\mathfrak{s}$ and, as can be shown by iterating the
argument used in the proof of Lemma~\ref{lem:CNTP}, $g$ permutes the
tensor factors of~$V$ correspondingly, i.e., there exists a permutation
$\pi$ of $\, \{1,\ldots,r\}$ \linebreak such that $\, \phi\,\mathfrak{s}_i
= \mathfrak{s}_{\pi(i)} \,$ and $\; g (v_1 \otimes \ldots \otimes v_r)
= g_{1,\pi^{-1}(1)}(v_{\pi^{-1}(1)}) \otimes \ldots \otimes g_{r,\pi^{-1}(r)}
(v_{\pi^{-1}(r)})$ \linebreak where $\phi_{\pi(i),i}$ is a Lie algebra
isomorphism from $\mathfrak{s}_i$ to $\mathfrak{s}_{\pi(i)}$ and
$g_{i,\pi^{-1}(i)}$ is a linear isomorphism from $V_{\pi^{-1}(i)}$
to~$V_i$ such that $\, \phi_{\pi(i),i}^{}(X_i^{})
= g_{i,\pi^{-1}(i)}^{} \, X_i^{} \, g_{i,\pi^{-1}(i)}^{-1}$.
But~then $\, \phi = \mathrm{Ad}(g) \,$ preserves~$\mathfrak{f}^\times$,
since $\mathfrak{f}^\times$ consists precisely of \emph{all} elements $X$
of~$\mathfrak{g}$ that are tensor decomposable.
The same argument prevails when $\mathfrak{s}\,$ is $\mathbb{Z}_2$-invariant
and  $\, \phi \in \mathrm{Aut}_{\mathbb{Z}_2}(\mathfrak{g},\mathfrak{s})
\setminus \mathrm{Inn}(\mathfrak{g},\mathfrak{s})$, since we can still
write $\, \phi = \mathrm{Ad}(g)$ \linebreak where $g$ is now an antiunitary
transformation on~$\mathbb{C}^n$ in the case of~$\mathfrak{su}(n)$ and an
orthogonal transformation on~$\mathbb{R}^n$ of determinant~$-1$ in the
case of $\mathfrak{so}(n)$: the hypothesis that $\mathfrak{s}$ should
be $\mathbb{Z}_2$-invariant merely states that the set of such
automorphisms $\phi$ is not empty.\,]
But this implies that if $\mathfrak{s}$ is assumed to be primitive, it
must in fact be equal to $\mathfrak{f}^\times$, and the same conclusion
holds if $\mathfrak{s}$ is assumed to be $\mathbb{Z}_2$-primitive,
provided $\mathfrak{f}^\times$ is $\mathbb{Z}_2$-invariant.
Note also that for $\, \mathfrak{g} = \mathfrak{su}(n)$, this
additional condition of $\mathbb{Z}_2$-invariance is automatically
satisfied by all the subalgebras of equation~(\ref{eq:CALGTP}),
provided we choose complex conjugation in $\, V = V_1 \otimes
\ldots \otimes V_r \,$ to be the tensor product of complex
conjugations in each tensor factor, but that for $\, \mathfrak{g}
= \mathfrak{so}(n)$ with $n$ even, it imposes restrictions on the
subalgebras of equation~(\ref{eq:RALGTP}) that are allowed;
these will be derived below.

Having shown that the irreducible primitive or $\mathbb{Z}_2$-primitive
subalgebras of classical Lie algebras must be sought among the subalgebras
or $\mathbb{Z}_2$-invariant subalgebras~$\mathfrak{f}^\times$ introduced
above, our next step will be to figure out the additional constraints
on the numbers $n_i$ imposed by primitivity or $\mathbb{Z}_2$-primitivity.
\linebreak
Beginning with the primitive case, we first note that any automorphism
of~$\mathfrak{su}(n)$, $\mathfrak{so}(n)$ or~$\mathfrak{sp}(n)$ which
preserves a subalgebra as defined in equations~(\ref{eq:CALGTP}) or~%
(\ref{eq:RALGTP}), respectively, can only permute tensor factors of
equal type and equal size. \linebreak
As a result, there are various situations in which such a subalgebra
$\mathfrak{f}^\times$ cannot be primitive, simply because one easily
finds a strictly larger proper subalgebra $\tilde{\mathfrak{f}}^\times$
of~$\mathfrak{g}$ which is $\mathrm{Inn}(\mathfrak{g},\mathfrak{f}^\times)$-%
invariant, namely:
\begin{itemize}
 \item If, in the case of~$\mathfrak{so}(n)$ or~$\mathfrak{sp}(n)$,
       $\mathfrak{f}^\times$ contains ``orthogonal blocks'' as well
       as ``symplectic blocks'', then it must contain
       precisely one ``orthogonal block'' and precisely a pair of
       ``symplectic blocks'' of type $\mathfrak{sp}(2)$ in the case
       of $\mathfrak{so}(n)$ ($s$ even) or a single ``symplectic
       block'' in the case of~$\mathfrak{sp}(n)$ ($s$ odd): otherwise,
       take  $\tilde{\mathfrak{f}}^\times$ to be $\, \mathfrak{so}(p)
       \times \mathfrak{so}(q) \,$ in the case of $\mathfrak{so}(n)$
       ($s$ even) or $\, \mathfrak{so}(p) \times \mathfrak{sp}(q) \,$
       in the case of $\mathfrak{sp}(n)$ ($s$ odd), with $\, p = n_1
       \ldots n_r \,$ and \linebreak $\, q = n_{r+1} \ldots n_{r+s} \,$,
       to conclude that $\mathfrak{f}^\times$ cannot be primitive.%
       \footnote{The possible occurrence of a pair of ``symplectic
       blocks'' of type $\mathfrak{sp}(2)$ is due to the exceptional
       isomorphism $\, \mathfrak{sp}(2) \times \mathfrak{sp}(2)
       \cong \mathfrak{so}(4)$, which entails that in this case,
       and only in this case, $\tilde{\mathfrak{f}}^\times$ is
       just $\mathfrak{f}^\times$ itself.}
 \item If, in the case of $\mathfrak{su}(n)$ or in the case of~%
       $\mathfrak{so}(n)$ or~$\mathfrak{sp}(n)$ with ``orthogonal
       blocks'' or ``symplectic blocks'' only, $\mathfrak{f}^\times$
       contains more than two blocks whose sizes are not all equal,
       and arranging the blocks according to their size, say in
       decreasing order, with $i$ denoting the first index for
       which $\, n_i > n_{i+1}$, take $\tilde{\mathfrak{f}}^\times$
       to be as follows:
 \begin{itemize}
  \item if $\, \mathfrak{g} = \mathfrak{su}(n)$, take
        $\, \tilde{\mathfrak{f}}^\times = \mathfrak{su}(p)
        \times \mathfrak{su}(q)$, with $\, p = n_1 \ldots n_i \,$
        and $\, q = n_{i+1} \ldots n_r \,$,
  \item if $\, \mathfrak{g} = \mathfrak{so}(n) \,$ and all blocks
        are orthogonal, take $\, \tilde{\mathfrak{f}}^\times
        = \mathfrak{so}(p) \times \mathfrak{so}(q)$, with
        $\, p = n_1 \ldots n_i \,$ and $\, q = n_{i+1} \ldots n_r \,$,
  \item if $\, \mathfrak{g} = \mathfrak{so}(n) \,$ and all blocks
        are symplectic (their number $s$ being even and $\, n_1
        = \ldots = n_i > 2$), take
       \[
        \begin{array}{cc}
         \tilde{\mathfrak{f}}^\times~
         =~\mathfrak{so}(p) \times \mathfrak{so}(q) & \mbox{if $i$ is even}
         \\[1mm]
         \tilde{\mathfrak{f}}^\times~
         =~\mathfrak{sp}(p) \times \mathfrak{sp}(q) & \mbox{if $i$ is odd}
        \end{array}
       \]
        with $\, p = n_1 \ldots n_i \,$ and $\, q = n_{i+1} \ldots n_s \,$,
  \item if $\, \mathfrak{g} = \mathfrak{sp}(n) \,$ and all blocks
        are symplectic (their number $s$ being odd and $\, n_1
        = \ldots = n_i > 2$), take
       \[
        \begin{array}{ccc}
         \tilde{\mathfrak{f}}^\times~
         =~\mathfrak{so}(p) \times \mathfrak{sp}(q) & \mbox{if $i$ is even}
         \\[1mm]
         \tilde{\mathfrak{f}}^\times~
         =~\mathfrak{sp}(p) \times \mathfrak{so}(q) & \mbox{if $i$ is odd}
        \end{array}
       \]
       with $\, p = n_1 \ldots n_i \,$ and $\, q = n_{i+1} \ldots n_s \,$,
 \end{itemize}
       to conclude that $\mathfrak{f}^\times$ cannot be primitive, except
       in the last two cases, provided that we choose $\, i=1 \,$ and\,%
       \addtocounter{footnote}{-1}\footnotemark
       \[
        \begin{array}{cc}
                 s=2         & \mbox{for $\mathfrak{so}(n)$ ($s$ even)}~,
         \\[1mm]
         s=3~,~n_2 = 2 = n_3 & \mbox{for $\mathfrak{sp}(n)$ ($s$ odd)}~.
        \end{array}
       \]
\end{itemize}
Thus we are left with the following candidates for irreducible primitive
subalgebras:
\[
 \begin{array}{ccl}
  \mathfrak{g}~=~\mathfrak{su}(n): &\quad&
  \left\{ \begin{array}{cc}
           \mathfrak{su}(p) \times \mathfrak{su}(q) &
           (n = p\>\!q) \\
           \mathfrak{su}(p) \times \ldots \times \mathfrak{su}(p)
           & (\mbox{$l$ factors, } n = p\>\!^l,\, l \geqslant 3)
          \end{array} \right. \\[3ex]
  \mathfrak{g}~=~\mathfrak{so}(n): &\quad&
  \left\{ \begin{array}{cc}
           \mathfrak{so}(p) \times \mathfrak{so}(q) & (n = p\>\!q) \\
           \mathfrak{sp}(2p) \times \mathfrak{sp}(2q) & (n = 4p\>\!q) \\
           \mathfrak{so}(p) \times \ldots \times \mathfrak{so}(p)
           & (\mbox{$l$ factors, } n = p\>\!^l, l \geqslant 3) \\
           \mathfrak{sp}(2p) \times \ldots \times \mathfrak{sp}(2p)
           & (\mbox{$l$ factors, } n = (2p)^l, l \geqslant 3 \mbox{ even})
          \end{array} \right. \\[6ex]
  \mathfrak{g}~=~\mathfrak{sp}(n): &\quad&
  \left\{ \begin{array}{cc}
           \mathfrak{sp}(2p) \times \mathfrak{so}(q) & (n = 2p\>\!q) \\
           \mathfrak{sp}(2p) \times \ldots \times \mathfrak{sp}(2p)
           & (\mbox{$l$ factors, } n = (2p)^l, l \geqslant 3 \mbox{ odd})
          \end{array} \right.
 \end{array}
\]
Finally, in the case of $\mathfrak{so}(n)$ with $n$ even, we must check which
of these are also $\mathbb{Z}_2$-invariant and hence $\mathbb{Z}_2$-primitive.
In order to see whether it is possible to construct an orthogonal trans%
formation $g$ on~$\mathbb{R}^n$ of determinant $-1$ that implements the
outer automorphism of $\mathfrak{so}(n)$ in such a way as to preserve
$\mathfrak{f}^\times$, we proceed case by case, making use of the fact
that, as before, $g$ can at most permute tensor factors of equal type
and equal size. Moreover, if $g$ preserves the tensor factors, then
it must be the tensor product of orthogonal transformations in each
tensor factor, each of which must have determinant $\pm 1$, and we
can compute its determinant from the following general formula for
determinants in tensor products:
\begin{equation} \label{eq:DETTP}
 {\det}_{V \otimes W}^{}(A \otimes B)~
 =~({\det}_V^{} A)^{\dim W} ({\det}_W^{} B)^{\dim V}~.
\end{equation}
\begin{itemize}
 \item $\mathfrak{f}^\times = \mathfrak{so}(p) \times \mathfrak{so}(q)$
       $(n = pq)$.
       Note that $p$ or $q$ (i.e., at least one of them) must be even
       (since $n$ is) and that if $\, p \neq q$, $g$ must preserve the
       tensor factors, while if $\, p = q$, $g$ may be chosen to switch
       them.
       In the first case, it follows from equation~(\ref{eq:DETTP})
       that we can achieve $\, \det g = -1 \,$ if and only if $p$
       or $q$ (i.e., the other one of them) is odd, namely by
       taking $g$ to be the tensor product of an orthogonal
       transformation of determinant~$-1$ in the even-dimensional
       tensor factor with the identity in the odd-dimensional
       tensor factor.
       In the second case, the switch operator has determinant
       $(-1)^{p(p-1)/2}$, so we can achieve $\, \det g = -1 \,$
       if and only if $p/2$ is odd.
 \item $\mathfrak{f}^\times = \mathfrak{sp}(2p) \times \mathfrak{sp}(2q)$
       $(n = 4pq)$.
       Note that if $\, p \neq q$, $g$ must preserve the tensor factors,
       while if $\, p = q$, $g$ may be chosen to switch them.
       In the first case, it follows from equation~(\ref{eq:DETTP})
       that we cannot achieve $\, \det g = -1$.
       In~the second case, the switch operator has determinant
       $(-1)^{p(2p-1)}$, so we can achieve $\, \det g = -1 \,$
       if and only if $p$ is odd.
 \item $\mathfrak{f}^\times = \mathfrak{so}(p) \times \ldots \times
       \mathfrak{so}(p)$ ($l$ factors, $n = p^l$, $l \geqslant 3$).
       Note that $p$ must be even (since $n$ is), so if $g$ is assumed
       to preserve the tensor factors, it follows from equation~%
       (\ref{eq:DETTP}) that we cannot achieve $\, \det g = -1$.
       Moreover, if $g$ switches two of the $l$ tensor factors but is
       the identity on the remaining $l-2$, then according to equation~%
       (\ref{eq:DETTP}) it has determinant $\, (-1)^{(p(p-1)/2) \cdot%
       p^{l-2}} = +1$, and since every permutation can be written as the
       product of transpositions, it follows that we cannot achieve
       $\, \det g = -1$.
 \item $\mathfrak{f}^\times = \mathfrak{sp}(2) \times \ldots \times
       \mathfrak{sp}(2p)$ ($l$ factors, $n = (2p)^l$, $l \geqslant 3$).
       Again, if $g$ is assumed to preserve the tensor factors, it
       follows from equation~(\ref{eq:DETTP}) that we cannot achieve
       $\, \det g = -1$.
       Similarly, if $g$ switches two of the $l$ tensor factors but is
       the identity on the remaining $l-2$, then according to equation~%
       (\ref{eq:DETTP}) it has determinant $\, (-1)^{(p(p-1)/2) \cdot
       p^{l-2}} = +1$, and since every permutation can be written as
       the product of transpositions, it follows that we cannot
       achieve $\, \det g = -1$.
\end{itemize}
To summarize, we see that among the irreducible primitive subalgebras
$\mathfrak{f}$ of $\mathfrak{so}(n)$ with $n$ even, as listed above,
the only ones which are $\mathbb{Z}_2$-invariant and hence also
$\mathbb{Z}_2$-primitive are $\, \mathfrak{so}(p) \times
\mathfrak{so}(q) \,$ ($n = pq$) with $\, p \neq q \,$
and either $p$ or $q$ odd, $\, \mathfrak{so}(p) \times
\mathfrak{so}(p) \,$ ($n = p^2$) with $p$ even and $p/2$
odd, and $\, \mathfrak{sp}(2p) \times \mathfrak{sp}(2p) \,$
($n = 4 p^2$) with $p$ odd.

This proves the classification of the primitive and $\mathbb{Z}_2$-%
primitive subalgebras of classical Lie algebras: the simple ones are
determined according to Theorem~\ref{theo:DYNEXC} above and the results
for the remaining ones are summarized in Tables~\ref{tab:SUA},
\ref{tab:SOA} and~\ref{tab:SPA}, providing the complete list of
(conjugacy classes of) non-simple primitive and $\mathbb{Z}_2$-%
primitive subalgebras of the classical Lie algebras $\mathfrak{su}(n)$,
$\mathfrak{so}(n)$ and $\mathfrak{sp}(n)$, respectively, with the
following conventions.

The first column indicates the isomorphism type of the subalgebra.
The second column indicates the conditions which are necessary for
the inclusion to exist and also for avoiding repetitions, most of
which are due to well-known canonical isomorphisms between classical
Lie algebras of low rank:
\begin{eqnarray*}
 &\mathfrak{so}(3)~=~\mathfrak{su}(2)~~~,~~~
  \mathfrak{sp}(2)~=~\mathfrak{su}(2)~~~,~~~
  \mathfrak{so}(2)~=~\mathbb{R}& \\[1mm]
 &\mathfrak{sp}(4)~=~\mathfrak{so}(5)~~~,~~~
  \mathfrak{so}(4)~=~\mathfrak{su}(2) \oplus \mathfrak{su}(2)& \\[1mm]
 &\mathfrak{so}(6)~=~\mathfrak{su}(4)&
\end{eqnarray*}
The third column indicates the intrinsic nature of the subalgebra,
classifying it into three types as above, namely abelian~($a$),
truly semisimple~($s$) and truly reductive~($r$), together with
the type of inclusion: reducible (red.) or irreducible (irred.).
Finally, Table~\ref{tab:SOA} has a fourth column, only relevant
for $\mathfrak{so}(n)$ with $n$ even, where we indicate under
what additional conditions the corresponding subalgebra is also
$\mathbb{Z}_2$-primitive.
Note that the subalgebras of~$\mathfrak{so}(2n)$ that are not
$\mathbb{Z}_2$-primitive fall into two distinct conjugacy classes
that are transformed into each other by an outer automorphism.
(In~Table~\ref{tab:SUA}, this additional column is superfluous
 since all primitive subalgebras are also $\mathbb{Z}_2$-primitive.)
All tables are divided into two parts: in the upper part, we list
the maximal subalgebras and in the lower part, we list the non-maximal
subalgebras.

The main reason for organizing these data in this form is for
better reference, for instance in the next section.


\begin{table}[!htb]
\begin{center}
{\small
\begin{tabular}{|c|c|c|c|} \hline
 \rule[-2mm]{0mm}{7mm} subalgebra \rule[-2mm]{0mm}{7mm} & conditions & type \\
 \hline\hline
 \rule[-1ex]{0mm}{5ex} $\R \oplus \mathfrak{su}(p) \oplus \mathfrak{su}(q)$ &
 $n = p+q\,$, $p \geqslant q \geqslant 1\,$, $p \geqslant 2$ &
 $r$ $\!-\!$ red. \\
 \rule[-1ex]{0mm}{4ex} $\R$ & $n = 2\,$ & $a$ $\!-\!$ red. \\
 \rule[-2ex]{0mm}{5ex} $\mathfrak{su}(p) \times \mathfrak{su}(q)$ &
 $n = p\>\!q\,$, $p \geqslant q \geqslant 2$ & $s$ $\!-\!$ irred. \\ \hline
 \rule[-1ex]{0mm}{6ex} $\R^{l-1} \oplus
                        \bigoplus\limits_{k=1}^l \mathfrak{su}(p)$ &
 $n = p\>\!l\,$, $l \geqslant 3\,$, $p \geqslant 2$ & $r$ $\!-\!$ red. \\
 \rule[-1ex]{0mm}{4ex} $\R^{n-1}$ & $n \geqslant 3$ & $a$ $\!-\!$ red. \\
 \rule[-3ex]{0mm}{7.5ex} $\prod\limits_{k=1}^l \mathfrak{su}(p)$ &
 $n = p\>\!^{l}\,$, $l \geqslant 3\,$, $p \geqslant 2$ & $s$ $\!-\!$ irred. \\
 \hline
\end{tabular}
}
\end{center}
\begin{center}
\caption{\label{tab:SUA}
          Non-simple primitive and $\mathbb{Z}_2$-primitive
          subalgebras of $\mathfrak{su}(n)$ ($n \geqslant 2$)}
\end{center}
\end{table}

\begin{table}[!htb]
\begin{center}
{\small
\begin{tabular}{|c|c|c|c|c|} \hline
 \rule[-5mm]{0mm}{13mm} subalgebra \rule[-5mm]{0mm}{13mm} & conditions & type &
 \begin{minipage}{7.5em}
  \begin{center}
   $\mathbb{Z}_2$-primitive \\ ($n$ even)
  \end{center}
 \end{minipage} \\
 \hline\hline
 \rule[-1ex]{0mm}{5ex} $\mathfrak{so}(p) \oplus \mathfrak{so}(q)$ &
 $n = p+q\,$, $p \geqslant q \geqslant 3$ & $s$ $\!-\!$ red. & yes \\
 \rule[-1ex]{0mm}{5ex} $\R \oplus \mathfrak{so}(p)$ &
 $n = p+2\,$, $p \geqslant 3$ & $r$ $\!-\!$ red. & yes \\
 \rule[-1ex]{0mm}{5ex} $\mathfrak{su}(2)\oplus\mathfrak{su}(2)$ &
 $n = 5$ & $s$ $\!-\!$ red. & $-$ \\
 \rule[-1ex]{0mm}{5ex} $\R \oplus \mathfrak{su}(p) = \mathfrak{u}(p)$ &
 $n = 2p\,$, $p \geqslant 3\,$ & $r$ $\!-\!$ red. & $p$ odd \\
 \rule[-2ex]{0mm}{6ex} $\mathfrak{so}(p) \times \mathfrak{so}(q)$ &
 $n = p\>\!q\,$, $p \geqslant q \geqslant 3\,$, $p,q \neq 4$ &
 $s$ $\!-\!$ irred. &
 \begin{minipage}{7.5em}
  \begin{center}
    $p \neq q$: $p$ or $\!q$ odd \\ $p = q$: $p/2$ odd
  \end{center}
 \end{minipage} \\
 \rule[-2ex]{0mm}{5ex} $\mathfrak{sp}(2p) \times \mathfrak{sp}(2q)$ &
 $n = 4p\>\!q\,$, $p \geqslant q \geqslant 1$ &
 $s$ $\!-\!$ irred. & $p = q$: $p$ odd \\ \hline
 \rule[-1ex]{0mm}{6ex} $ \bigoplus\limits_{k=1}^l \mathfrak{so}(p)$ &
 $n = p\>\!l\,$, $l \geqslant 3\,$, $p \geqslant 3$ & $s$ $\!-\!$ red. & yes \\
 \rule[-2ex]{0mm}{5ex} $\R^l$ &
 $n = 2l\,$, $l \geqslant 3$ & $a$ $\!-\!$ red. & yes \\ 
 \rule[-1ex]{0mm}{5ex} $ \prod\limits_{k=1}^l \mathfrak{so}(p)$ &
 $n = p\>\!^l\,$, $l \geqslant 3\,$, $p \geqslant 3\,$, $p \neq 4\,$ &
 $s$ $\!-\!$ irred. & no \\
 \rule[-1ex]{0mm}{5ex} $\prod\limits_{k=1}^l \mathfrak{sp}(2p)$ &
 $n=(2p)^l\,$, $l \geqslant 4\,$, $l$ even, $p \geqslant 1$ &
 $s$ $\!-\!$ irred. & no \\  
 \rule[-2ex]{0mm}{6ex}
 $\mathfrak{so}(p) \!\times\! \mathfrak{sp}(2) \!\times\! \mathfrak{sp}(2)$ &
 $n = 4p\,$, $p \geqslant 3\,$, $p \neq 4$ & $s$ $\!-\!$ irred. & $p$ odd \\
 \hline
\end{tabular}
}
\end{center}
\begin{center}
 \caption{\label{tab:SOA}
          Non-simple primitive and $\mathbb{Z}_2$-primitive
          subalgebras of $\mathfrak{so}(n)$ ($n \geqslant 5$)}
\end{center}
\vspace{-0.5ex}
\end{table}

\begin{table}[!htb]
\begin{center}
{\small
\begin{tabular}{|c|c|c|} \hline
 \rule[-2mm]{0mm}{7mm} subalgebra \rule[-2mm]{0mm}{7mm} & conditions & type \\
 \hline\hline
 \rule[-1ex]{0mm}{5ex} $\mathfrak{sp}(2p) \oplus \mathfrak{sp}(2q)$ &
 $n = 2(p+q)\,$, $p \geqslant q \geqslant 1$ & $s$ $\!-\!$ red. \\
 \rule[-1ex]{0mm}{4ex} $\R \oplus \mathfrak{su}(p) = \mathfrak{u}(p)$ &
 $n = 2p\,$, $p \geqslant 2$ & $r$ $\!-\!$ red. \\
 \rule[-2ex]{0mm}{5ex} $\mathfrak{sp}(2p) \times \mathfrak{so}(q)$ &
 $n = 2p\>\!q\,$, $p \geqslant 1\,$, $q \geqslant 3\,$, $q \neq 4$ &
 $s$ $\!-\!$ irred. \\ \hline
 \rule[-1ex]{0mm}{6ex} $\bigoplus\limits_{k=1}^l \mathfrak{sp}(2p)$ &
 $n = 2p\>\!l\,$, $l \geqslant 3\,$, $p \geqslant 1$ & $s$ $\!-\!$ red. \\
 \rule[-1ex]{0mm}{5ex} $\prod\limits_{k=1}^l \mathfrak{sp}(2p)$ &
 $n = (2p)^l\,$, $l \geqslant 3\,$, $l$ odd, $p \geqslant 1$ &
 $s$ $\!-\!$ irred. \\
 \rule[-2ex]{0mm}{6ex}
 $\mathfrak{sp}(2p) \!\times\! \mathfrak{sp}(2) \!\times\! \mathfrak{sp}(2)$ &
 $n = 8p\,$, $p \geqslant 2$ & $s$ $\!-\!$ irred. \\ \hline
\end{tabular}
}
\end{center}
\begin{center}
 \caption{\label{tab:SPA}
          Non-simple primitive subalgebras of $\mathfrak{sp}(n)$
          ($n$ even, $n \geqslant 4$)}
\end{center}
\vspace{-2ex}
\end{table}

\section{Maximal subgroups of the classical groups}

Our goal in this last section is to compute the normalizers, within
the corresponding classical groups $G$, of the primitive subalgebras
of the classical Lie algebras and thus obtain a complete list of all
maximal subgroups of the classical groups.

Beginning with the case of simple primitive subalgebras, we distinguish two
types: (i)~the ``classical'' simple inclusions $\, \mathfrak{so}(n) \subset
\mathfrak{su}(n)$, $\mathfrak{sp}(n) \subset \mathfrak{su}(n) \,$ and
$\, \mathfrak{so}(n-1) \subset \mathfrak{so}(n) \,$ and (ii)~the
remaining ``non-trivial'' simple inclusions given by the irreducible
representations that do not belong to Dynkin's list of exceptions.
For the ``classical'' simple inclusions it is easy to find the
corresponding maximal subgroups: $O(n) \subset SU(n)$,
$\SPG(n) \subset SU(n) \,$ and $\, O(n-1) \subset SO(n)$.

Passing to the non-simple primitive subalgebras, which have been
determined in the previous section, what remains to be done is to
compute the corresponding list of maximal subgroups of the
classical groups.
The final results are summarized in Table~\ref{tab:SUG} for the
group $SU(n)$, in Tables~\ref{tab:SOG1} and~\ref{tab:SOG2} for the
group $SO(n)$ and in Table~\ref{tab:SPG} for the group $\SPG(n)$.
For each maximal subgroup $H$ of~$G$, we list in the first column
the isomorphism type of the connected one-component $H_0$ of~$H$
and in the second column the component group $H/H_0$, taking into
account that $H$ must be the normalizer $N_G(H_0)$ of~$H_0$ in~$G$.
This information completely characterizes $H$.

To carry out the concrete calculations, we shall employ two tools.
The first tool consists in the introduction of an ``intermediate
subgroup'' between $H_0$ and $H$, which will be denoted by $H_i$
and will allow us to control the size and the structure of $H/H_0$.

Suppose $G$ is a connected semisimple Lie group with Lie algebra~%
$\mathfrak{g}$ and $H$ is a maximal subgroup of $G$ with connected
one-component~$H_0$ and Lie algebra~$\mathfrak{h}$.
Assume that $\mathfrak{h}$ is not an ideal of~$\mathfrak{g}$.
In the particular case where $\mathfrak{g}$ is simple~-- which is
the case of interest here~-- this is equivalent to assuming that
$\mathfrak{h} \neq \{0\}$.
According to Theorem~\ref{theo:maxsg1}, $H$ is then equal to
the normalizer $N_G(H_0)$ of~$H_0$ in~$G$.

Now consider the centralizer $Z_G(H_0)$ of~$H_0$ in~$G$ and define
$H_i$ to be the closed subgroup of~$G$ generated by the two closed
subgroups $H_0$ and $Z_G(H_0)$:
\begin{equation} \label{eq:NORMCENT0}
 H_i = H_0 \, Z_G(H_0) = Z_G(H_0) \, H_0~.
\end{equation}
Then we have\footnote{As usual, the symbol ``$\lhd$'' is to be read
as ``is normal subgroup of''.}
\begin{equation} \label{eq:NORMCENT1}
 H_0 \lhd H_i \lhd H~.
\end{equation}
Considering the action of~$H$ on~$H_0$ by conjugation in~$G$ (recall
that $N_G(H_0)$ is defined as the set of elements $g$ of~$G$ such that
conjugation by~$g$ leaves $H_0$ invariant), $H_i$ consists of those
elements of $H$ that act as inner automorphisms on~$H_0$, while
$Z_G(H_0)$ consists of those elements of $H$ that act trivially
on~$H_0$.
In other words, under the adjoint representation $\, \mathrm{Ad}:
G \longrightarrow \mathrm{Aut}(\mathfrak{g}) \,$ (which maps into
$\mathrm{Inn}(\mathfrak{g})$ because $G$ is connected and maps onto
$\mathrm{Inn}(\mathfrak{g})$ because $G$ is semisimple), $H$ is the
inverse image of $\mathrm{Inn}(\mathfrak{g},\mathfrak{h})$ while
$H_i$ is the inverse image of $\mathrm{Inn}(\mathfrak{h})$, so
\begin{equation} \label{eq:NORMCENT2}
 H/H_i~\cong~\mathrm{Inn}(\mathfrak{g},\mathfrak{h})\,/\,
             \mathrm{Inn}(\mathfrak{h})~.
\end{equation}
Also note that $H$, $H_i$ and $H_0$ all have the same Lie algebra
$\mathfrak{h}$ (since this is true for $H_0$ and $H$) and thus
the quotient groups $H/H_0$, $H/H_i$ and $H_i/H_0$ are discrete.
Moreover, by the second isomorphism theorem of group theory, we have
\begin{equation} \label{eq:NORMCENT3}
 H/H_i~=~\frac{H/H_0}{H_i/H_0}~.
\end{equation}
If $G$ is compact, the above mentioned discrete groups are in fact
finite and their orders satisfy the relation
\begin{equation} \label{eq:NORMCENT4}
 |H/H_0|~=~|H_i/H_0|\,|H/H_i|~.
\end{equation}
In almost all cases of interest, we may use additional properties to
extract more information about the structure of $H_i/H_0$ and of $H/H_i$,
and hence of $H/H_0$ as well.
\begin{itemize}
 \item $H_i/H_0\,$: Since $\, H_0 \cap Z_G(H_0) = Z(H_0)$, the first 
       isomorphism theorem implies
       \begin{equation} \label{eq:NORMCENT5}
        H_i/H_0~\cong~Z_G(H_0) \, / \, Z(H_0)~.
       \end{equation}
       Moreover, the inclusion $\, Z(G) \subset Z_G(H_0) \,$ is valid
       in general. Conversely, if $G$ is a classical group and the
       inclusion of $H_0$ in $G$ is an irreducible representation,
       then by Schur's lemma we have $\, Z_G(H_0) \subset (\mathbb{K}
       \cdot 1) \cap G \subset Z(G)$ \linebreak and thus
       \begin{equation} \label{eq:NORMCENT6}
        Z_G(H_0)~=~Z(G) \quad , \quad
        H_i/H_0~\cong~Z(G) \, / \, Z(H_0)~.
       \end{equation}
       More generally, by Schur's lemma, $Z_G(H_0)$ acts by scalar
       multiplication on each irreducible component.
       Finally, we note that if $\mathfrak{h}$ has maximal rank,
       that is, if $H_0$ contains a maximal torus $T$ of~$G$, then
       $\, Z_G(H_0) \subset T \subset H_0 \,$ and hence $\, H_i = H_0$.
 \item $H/H_i\,$: This quotient group may identified with the group
       $\, \mathrm{Out}(\mathfrak{h}) \cap \mathrm{Inn}(\mathfrak{g})$
       \linebreak
       of those outer automorphisms of~$\mathfrak{h}$ that can be extended
       to inner auto\-morphisms of~$\mathfrak{g}$ and thus can be implemented
       by conjugation by elements of the ambient group~$G$.
       If $\mathfrak{h}$ is semisimple, it may also be identified with the
       group $\, \mathrm{Aut}(\Gamma) \cap \mathrm{Inn}(\mathfrak{g}) \,$ of
       those automorphisms of the Dynkin diagram $\Gamma$ of $\mathfrak{h}$
       that can be extended to inner automorphisms of $\mathfrak{g}$ and
       thus can be implemented by conjugation by elements of the ambient
       group $G$:
       \begin{equation} \label{eq:NORMCENT7}
        H/H_i~\cong~\mathrm{Aut}(\Gamma) \cap \mathrm{Inn}(\mathfrak{g})~.
       \end{equation}
       In what follows, we shall calculate $H/H_i$ using the results of
       the previous section, in particular Lemma~\ref{lem:CNDS} and
       Lemma~\ref{lem:CNTP}.
\end{itemize}
\begin{remark}
 According to the relation~(\ref{eq:NORMCENT3}), the group $H/H_0$ is
 an upwards extension of the group $H_i/H_0$ by the group $H/H_i$, but
 simple examples show that this extension is not necessarily split,
 that is, a semi-direct product.
\end{remark}

The second tool is the following lemma which, in combination with
Lemma~\ref{lem:CNDS}, allows to handle the reducible inclusions.
\begin{lemma} \label{lem:BLDIAG}
 Assume that $G$ is one of the classical groups $SU(n)$, $SO(n)$ or $\SPG(n)$
 and consider the standard embeddings
 \[
  \mathfrak{g}(p) \oplus \mathfrak{g}(q) \subset \mathfrak{g}(n) \quad
  (n = p + q)
 \]
 by block diagonal matrices, where $\mathfrak{g}(r)$ is a generic symbol
 standing for the corresponding type of classical Lie algebra (i.e.,
 $\mathfrak{su}(r)$, $\mathfrak{so}(r)$ or $\mathfrak{sp}(r)$, respectively).
 Then if $\, g \in G \,$ normalizes both $\mathfrak{g}(p)$ and
 $\mathfrak{g}(q)$, $g$ must itself be block diagonal.
 The same statement holds for more than two direct summands.
\end{lemma}
\begin{proof}
 Writing $\, g = \bigl( \begin{smallmatrix} A & B \\ C & D \end{smallmatrix}
 \bigr) \,$ and assuming that $g$ normalizes matrices of the form $\bigl(
 \begin{smallmatrix} X & 0 \\ 0 & 0 \end{smallmatrix} \bigr)$ as well
 as matrices of the form $\bigl( \begin{smallmatrix} 0 & 0 \\ 0 & Y
 \end{smallmatrix} \bigr)$, it follows immediately that $B$ and $C$
 must vanish, since the $X$'s and $Y$'s form irreducible sets of
 matrices.
\end{proof}

With these generalities out of the way, we proceed to the case by case
analysis of each of the entries of Tables~\ref{tab:SUG}--\ref{tab:SPG}.
\vspace{1ex}

\subsection{Maximal subgroups of $SU(n)$}

\begin{enumerate}
 \item The inclusion
       \[
        S \bigl( U(p) \times U(q) \bigr)~\subset~SU(n) \qquad (n = p+q)
       \]
       is given by the direct sum of the defining representations of
       $U(p)$ and $U(q)$. Explicitly, it is obtained as the restriction
       of the inclusion of $\, U(p) \times U(q) \,$ into $U(n)$ by block
       diagonal matrices. In fact, the inverse image of $\, SU(n) \subset
       U(n)$ under this inclusion is the subgroup
       \[
        H_0~=~S \bigl( U(p) \times U(q) \bigr)~
        =~\{ \, (A,B) \in U(p) \times U(q)~|~\det(A) \det(B) = 1 \, \}
       \]
       consisting of the block diagonal matrices belonging to $SU(n)$,
       which is connected.
       It follows from Lemma~\ref{lem:CNDS} that the centralizer $Z_G(H_0)$
       of~$H_0$ in~$G$ is contained in $H_0\,$; hence $\, H_i = H_0$.
       Similarly, we compute $\, H = N_G(H_0)$ \linebreak by combining
       Lemma~\ref{lem:CNDS} and Lemma~\ref{lem:BLDIAG} to conclude
       that the elements $h$ of~$H$ must be of the form $\, h =
       \bigl( \begin{smallmatrix} A & 0 \\ 0 & D \end{smallmatrix}
       \bigr) \,$ or $\, h = \bigl( \begin{smallmatrix} 0 & B \\ C & 0
       \end{smallmatrix} \bigr)$, with $\, p = q \,$ in the second case.
       In the first~case, it follows that $\, h \in H_0$, and thus
       $\, H = H_0 \,$ if $\, p \neq q$.
       In the second case, we may write $h$ as the product of an
       element of~$H_0$ with the matrix $\, \bigl( \begin{smallmatrix}
       0 & 1 \\ -1 & 0 \end{smallmatrix} \bigr) \in SU(n) \,$ ($n = 2p$),
       whose square is $-1_n$ and belongs to~$H_0$, so we may conclude
       that $\, H/H_0 = \mathbb{Z}_2 \,$ if $\, p = q$.
 \item The inclusion
       \[
        SU(p) \times_{\mathbb{Z}_d} SU(q)~\subset~SU(n) \qquad (n = p\>\!q)
       \]
       is given by the tensor product of the defining representations
       of $U(p)$ and $U(q)$, and $\, d = \mathrm{gcd}(p,q) \,$ is the
       greatest common divisor of $p$ and $q$. \linebreak
       Explicitly, it is induced by the homomorphism
       \[
        \begin{array}{ccc}
         U(p) \times U(q) & \longrightarrow &    U(n)     \\[1mm]
              (A,B)       &   \longmapsto   & A \otimes B
        \end{array}
       \]
       which, in contrast to the situation encountered in the first item,
       is not injective but rather has a non-trivial kernel given by
       \[
        \{ \, (\exp(i\alpha) 1_p\,,\exp(-i\alpha) 1_q)~|~
              \alpha \in \mathbb{R} \, \}~\cong~U(1)~.
       \]
       Using the formula
       \[
        \det(A \otimes B)~=~(\det(A))^q \, (\det(B))^p~,
       \]
       we see that the inverse image of $\, SU(n) \subset U(n) \,$ under this
       homomorphism is the subgroup 
       \[
        S' \bigl( U(p) \times U(q) \bigr)~
        =~\{ \, (A,B) \in U(p) \times U(q)~|~
             \det(A)^q \det(B)^p = 1 \, \}~.
       \]
       and the restricted homomorphism
       \[
        \begin{array}{ccc}
         S' \bigl( U(p) \times U(q) \bigr) & \longrightarrow &    SU(n)   
         \\[1mm]
                       (A,B)               &   \longmapsto   & A \otimes B
        \end{array}
       \]
       still has the same kernel
       \[
        \{ \, (\exp(i\alpha) 1_p\,,\exp(-i\alpha) 1_q)~|~
              \alpha \in \mathbb{R} \, \}~\cong~U(1)~.
       \]
       Its intersection with the connected subgroup
       $\, SU(p) \times SU(q) \,$ is
       \[
        \{ \, (\exp(2 \pi i k/d) 1_p\,,\exp(- 2 \pi i k/d) 1_q)~|~
              0 \leqslant k < d \, \}
       \]
       and is isomorphic to $\mathbb{Z}_d$.
       Factoring this out, we obtain the desired inclusion.
       Note also that the center $Z(G)$ of~$G$ is contained in the
       quotient group $\, S' \bigl( U(p) \times U(q) \bigr) / U(1) \,$
       which is generated by $Z(G)$ and by the quotient group
       $\, \bigl( SU(p) \times SU(q) \bigr) / \mathbb{Z}_d \,$;
       hence
       \[
        H_0~\cong~\bigl( SU(p) \times SU(q) \bigr) / \mathbb{Z}_d~,
       \]
       and
       \[
        H_i~\cong~S' \bigl( U(p) \times U(q) \bigr) / U(1)~,
       \]
       so according to equation (\ref{eq:NORMCENT6}),
       \[
        H_i/H_0~\cong~Z(G)/Z(H_0)~\cong~\mathbb{Z}_n / \mathbb{Z}_{n/d}~
                \cong~\mathbb{Z}_d~.
       \]
       Explicitly, a representative of the connected component
       of~$H_i$ corresponding to $k \mod d$ is given by the matrix
       $\exp(2 \pi i k/n) \, 1_n$, since when $k$ is a multiple of
       $d$, $k/n$ will be a multiple of $d/n$ which can be written
       in the form 
       \[
        \frac{d}{n}~=~\frac{r}{p} + \frac{s}{q}
       \]
       (where $r$ and $s$ are chosen such that $rq'+sp^{\,\prime} = 1$
       with $p^{\,\prime} = p/d$ and $q' = q/d$, which is possible since
       $p^{\,\prime}$ and $q'$ are relatively prime) and therefore
       $\, \exp(2 \pi i k/n) \, 1_n \in H_0$. \\[2mm]
       In order to compute the normalizer $H$ of the connected subgroup
       $H_0$ so defined, we first note that (for $p \geqslant 3$) the
       Dynkin diagram of $\mathfrak{su}(p)$ admits only one automorphism
       which can be implemented by an antilinear involution $\sigma_p$
       in~$\mathbb{C}^p$.
       Since the tensor product of linear\,/\,antilinear maps is a
       linear\,/\,antilinear map and the tensor product between a
       linear map and an antilinear map is not well defined, we see
       that there is only one non-trivial automorphism of the Dynkin
       diagram $\Gamma$ of $\mathfrak{su}(p) \times \mathfrak{su}(q)$
       induced by automorphisms of the Dynkin diagrams of the factors
       which can be extended to an automorphism of $\mathfrak{su}(n)$,
       namely the one implemented by an anti-linear involution of the
       form $\, \sigma_p \otimes \sigma_q$, but this will always be
       an outer automorphism of $\mathfrak{su}(n)$.
       Thus it becomes clear that for $\, p \neq q$, $\mathrm{Aut}%
       (\Gamma) \cap \mathrm{Inn}(\mathfrak{g})=\{1\}$, while for
       $\, p = q$, $\mathrm{Aut}(\Gamma) \cap \mathrm{Inn}(\mathfrak{g})$
       can at most be equal to the group $\Z_2$, corresponding to the
       possibility of switching the factors in the tensor product
       that exists in this case.
       Explicitly, a representative of the corresponding connected
       component of~$H$ is given by the transformation
       \[
        \begin{array}{ccc}
          \mathbb{C}^n   & \longrightarrow &
                 \mathbb{C}^n              \\[1mm]
         z = x \otimes y &   \longmapsto   &
         z^\tau = \exp(i\phi_p) \, y \otimes x
        \end{array}~,
       \]
       which belongs to  $SU(n)$ ($n = p^2$) if we choose the phase
       $\exp(i\phi_p)$ according to
       \[
        \exp(i\phi_p)~
        =~\left\{ \begin{array}{ccc}
                   \mbox{$+1$ or $-1$} & \mbox{if} & p=0 \mod 4 \\[1mm]
                           + 1         & \mbox{if} & p=1 \mod 4 \\[1mm]
                      \exp(i\pi/4)     & \mbox{if} & p=2 \mod 4 \\[1mm]
                           - 1         & \mbox{if} & p=3 \mod 4 \\[1mm]
                  \end{array} \right.~,
       \]
       since the permutation which maps $x \otimes y$ to $y \otimes x$,
       when written in a basis of $\mathbb{C}^n$ provided by a basis of
       $\mathbb{C}^p$ by taking tensor products between vectors of the
       latter, is the product of $\binom{p}{2} = p\>\!(p-1)/2$ transpositions
       \[
        e_i \otimes e_j \longleftrightarrow e_j \otimes e_i
        \qquad \mbox{$(1 \leqslant i < j \leqslant p)$}
       \]
       with $p$ fixed points $\, e_i \otimes e_i$ $(1 \leqslant i
       \leqslant p)$ and hence is an involution with determinant
       $(-1)^{p(p-1)/2}$.
       Therefore, we conclude that for $p = q$, \linebreak
       $\mathrm{Aut}(\Gamma) \cap \mathrm{Inn}(\mathfrak{g}) = \mathbb{Z}_2$.
       Finally, the structure of the group $H/H_0$ in this case can
       be deduced by observing that for $p \neq 2 \mod 4$, the map
       $.{}^\tau$ has order $2$ and hence $H/H_0$ is the direct product
       $\mathbb{Z}_p \times \mathbb{Z}_2$, whereas for $p = 2 \mod 4$,
       i.e., $p = 2r$ with $r$ odd, it has order $8$, but its square
       coincides with the matrix $\, \exp(2 \pi i r^2/p^{\,2}) \, 1_n
       = i \, 1_n \,$ which represents one of the non-trivial connected
       components of $H_i\,$, and the square of the latter belongs to
       $H_0$, so it follows that $H/H_0$ is a non-split upward extension
       of the group $\mathbb{Z}_p$ by the group $\mathbb{Z}_2$
       ($H/H_0 = \mathbb{Z}_p \,.\, \mathbb{Z}_2$) which can be
       explicitly constructed as the quotient group
       $\mathbb{Z}_p \times_{\mathbb{Z}_2} \mathbb{Z}_4$.
 \item Generalizing the procedure of the first item, we define the
       inclusion
       \[
        S \bigl( U(p) \times \ldots \times U(p) \bigr)~\subset~SU(n)
        \qquad (n = p\>\!l)
       \]
       by the direct sum of $l$ copies of the defining representation of
       $U(p)$, which realizes $\, U(p) \times \ldots \times U(p) \,$ by
       block diagonal matrices in $U(n)$, with $l$ blocks of size $p$
       along the diagonal. The inverse image of $\, SU(n) \subset U(n) \,$
       under this inclusion is the subgroup
       \[
        H_0~=~S \bigl( U(p) \times \ldots \times U(p) \bigr)
       \]
       consisting of the block diagonal matrices belonging to~$SU(n)$,
       which is connected.
       As before, it follows from Lemma~\ref{lem:CNDS} that the
       centralizer $Z_G(H_0)$ of~$H_0$ in~$G$ is contained in $H_0\,$;
       hence $\, H_i = H_0$.
       Similarly, we compute $\, H = N_G(H_0) \,$ by combining
       Lemma~\ref{lem:CNDS} and Lemma~\ref{lem:BLDIAG} to conclude
       that the elements $h$ of~$H$ must be such that, when represented
       as block $(l \times l)$-matrices (with entries that are themselves
       $(p \times p)$-matrices), they contain precisely one nonvanishing
       entry in each line and each column: obviously, any such matrix $h$
       defines a permutation $\sigma(h)$ of $\{1,\ldots,l\}$. Thus we
       obtain a group homomorphism from~$H$ to the permutation group~%
       $S_l$ which has kernel~$H_0$ and is surjective (we can use the
       same argument as in item~1 above to show that its image contains
       all transpositions), so we may conclude that $\, H/H_0 = S_l$.
 \item Generalizing the procedure of the second item, we define the
       inclusion 
       \[
        \bigl( SU(p) \times \ldots \times SU(p) \bigr) / \mathbb{Z}_p^{l-1}~
        \subset~SU(n) \qquad (n = p\>\!^l)
       \]
       by the tensor product of $l$ copies of the defining representation
       of $U(p)$. Explicitly, it is induced by the homomorphism
       \[
        \begin{array}{ccc}
         U(p) \times \ldots \times U(p)
         & \longrightarrow &             U(n)               \\[1mm]
               (A_1,\ldots,A_l)        
         &   \longmapsto   & A_1 \otimes \ldots \otimes A_l
        \end{array}
       \]
       which, once again, is not injective but rather has a non-trivial
       kernel given by
       \[
        \{ \, (\exp(i \alpha_1) \, 1_p\,,\ldots,\exp(i \alpha_l) \, 1_p)~|~
              \exp(i (\alpha_1 + \ldots + \alpha_l)) = 1 \, \}~
        \cong~U(1)^{l-1}~.
       \]
       Using the formula
       \[
        \det \bigl( A_1 \otimes \ldots \otimes A_l \bigr)~
        =~\bigl( \det(A_1) \,\ldots\, \det(A_l) \bigr)^{p^{\,l-1}}~,
       \]
       we see that the inverse image of $\, SU(n) \subset U(n) \,$ under
       this homomorphism is the subgroup
       \[
        S' \bigl( U(p) \times \ldots \times U(p) \bigr)~
        =~\begin{array}{l}
           \{ \, (A_1,\ldots,A_l) \in U(p) \times \ldots \times U(p)~|~ \\
           \hspace*{3em}
           \bigl( \det(A_1) \,\ldots\, \det(A_l) \bigr)^{p^{\,l-1}} = 1 \, \}~,
        \end{array}
       \]
       and the restricted homomorphism 
       \[
        \begin{array}{ccc}
         S' \bigl( U(p) \times \ldots \times U(p) \bigr)
         & \longrightarrow &             SU(n)              \\[1mm]
                         (A_1,\ldots,A_l)                     
         &   \longmapsto   & A_1 \otimes \ldots \otimes A_l
        \end{array}
       \]
       still has the same kernel
       \[
        \{ \, (\exp(i \alpha_1) \, 1_p\,,\ldots,\exp(i \alpha_l) \, 1_p)~|~
              \exp(i (\alpha_1 + \ldots + \alpha_l)) = 1 \, \}~
        \cong~U(1)^{l-1}~.
       \]
       Its intersection with the connected subgroup
       $\, SU(p) \times \ldots \times SU(p) \,$ is
       \[
        \{(\exp(2\pi i k_1/p) \, 1_p\,,\ldots,\exp(2\pi i k_l/p) \, 1_p)~|~
        \exp(2\pi i (k_1 + \ldots + k_l)/p) = 1\}
       \]
       and is isomorphic to $\mathbb{Z}_p^{l-1}$.
       Factoring this out, we obtain the desired inclusion.
       Note also that the center $Z(G)$ of $G$ is contained in the quotient
       group $\, S' \bigl( U(p) \times\ldots\times U(p) \bigr)/U(1)^{l-1} \,$
       which is generated by $Z(G)$ and by the quotient group
       $\, \bigl( SU(p) \times\ldots\times SU(p) \bigr)/\mathbb{Z}_p^{l-1}$;
       hence
       \[
        H_0~\cong~\bigl( SU(p) \times \ldots \times SU(p) \bigr) /
                  \mathbb{Z}_p^{l-1}~,
       \]
       and
       \[
        H_i~\cong~S' \bigl( U(p) \times \ldots \times U(p) \bigr) /
                  U(1)^{l-1}~,
       \]
       so according to equation~(\ref{eq:NORMCENT6}),
       \[
        H_i/H_0~\cong~Z(G)/Z(H_0)~\cong~\mathbb{Z}_n / \mathbb{Z}_p~
                \cong~\mathbb{Z}_{p\>\!^{l-1}}~.
       \]
       Explicitly, a representative of the connected component of
       $H_i$ corresponding to $k \mod p\,^{l-1}$ is given by the
       matrix $\exp(2 \pi i k/n) \, 1_n$, since when $k$ is a multiple
       of~$p\,^{l-1}$, $\exp (2 \pi i k/n) \, 1_n \in H_0$. \\[2mm]
       In order to compute the normalizer $H$ of the connected subgroup
       $H_0$ so defined, we note, as before, that there is no automorphism
       of the Dynkin diagram of $\, \mathfrak{su}(p) \times \ldots \times
       \mathfrak{su}(p) \,$ induced by automorphisms of the Dynkin diagrams
       of the factors that can be extended to an inner automorphism of~%
       $\mathfrak{su}(n)$.
       Thus it becomes clear that $\mathrm{Aut}(\Gamma) \cap \mathrm{Inn}%
       (\mathfrak{g})$ can at most be equal to the symmetric group $S_l$,
       corresponding to the possibility of permuting the $l$ factors in
       the tensor product.
       But for $\, l \geqslant 3$, we can implement any transposition
       through an involution of determinant $+1$ since, for example,
       the determinant of the transformation
       \[
        \begin{array}{ccc}
         \mathbb{C}^n   
         & \longrightarrow &
         \mathbb{C}^n    \\[1mm]
         z = x_1 \otimes x_2 \otimes x_3 \otimes \ldots \otimes x_l
         &   \longmapsto   &
         z^{\tau_{12}} = \pm \, x_2 \otimes x_1 \otimes x_3
                                    \otimes \ldots \otimes x_l
        \end{array}~,
       \]
       which represents the transposition $\tau_{12}$ ($\tau_{12}(1) = 2$,
       $\tau_{12}(2) = 1$, $\tau_{12}(i) = i$ for $i \geqslant 3$) is
       $(-1)^{p^{\,l-1}(p \mp 1)/2}$, and this is equal to $+1$ when $p$ is
       even and, with an appropriate choice of sign, also when $p$ is odd.
       Therefore, it follows that $\, \mathrm{Aut}(\Gamma) \cap 
       \mathrm{Inn}(\mathfrak{g}) = S_l \,$ and that $H/H_0$
       is the direct product $\, \mathbb{Z}_{p^{\,l-1}} \times S_l$.
\end{enumerate}

\begin{table}[!htb]
\begin{center}
\begin{tabular}{|c|c@{\hspace{0.5em}}c@{\hspace{0.5em}}c|} \hline
 \rule[-2mm]{0mm}{7mm} {\small connected component $H_0$} 
 \rule[-2mm]{0mm}{7mm} &
 \multicolumn{3}{c|}{{\small component group $H/H_0$}}  \\ \hline\hline
 \rule[-5ex]{0mm}{13ex} $\begin{array}{c}
                         S \bigl( U(p) \times U(q) \bigr) \\[1mm]
                         \mbox{(reducible)} \\ \mbox{(direct sum)}
                        \end{array}$ \rule[-5ex]{0mm}{13ex} &
 $\begin{array}{c} \{1\} \\[1mm] \Z_2 \end{array}$ &
 $\begin{array}{c} \mbox{for} \\[1mm] \mbox{for} \end{array}$ &
 $\begin{array}{c}
   n = p+q\,,\, p > q \geqslant 1 \\[1mm]
   n = p+q\,,\, p = q \geqslant 1
  \end{array}$ \\
 \rule[-6ex]{0mm}{15ex} $\begin{array}{c}
                         SU(p) \times_{\Z_d} SU(q) \\[1mm]
                         \mbox{(irreducible)} \\ \mbox{(tensor product)} \\
                         (d = \mathrm{gcd}(p,q))
                        \end{array}$ \rule[-6ex]{0mm}{15ex} &
 $\Z_d$ & $\mbox{for}$ & $n = p\>\!q\,,\, p > q \geqslant 2$ \\
 \rule[-6ex]{0mm}{13ex} $\begin{array}{c}
                         SU(p) \times_{\Z_p} SU(p) \\[1mm]
                         \mbox{(irreducible)} \\ \mbox{(tensor product)}
                        \end{array}$ \rule[-6ex]{0mm}{14ex} &
 $\begin{array}{c}
   \Z_p \times \Z_2 \\[1mm] \Z_p \times_{\Z_2} \Z_4
  \end{array}$ &
 $\begin{array}{c} \mbox{for} \\[1mm] \mbox{for} \end{array}$ &
 $\begin{array}{c}
   n = p^2\,,\, p \geqslant 2\,,\, p \neq 2 \mod 4 \\[1mm]
   n = p^2\,,\, p \geqslant 2 \,,\, p = 2 \mod 4
  \end{array}$ \\ \hline
 \rule[-6ex]{0mm}{15ex} $\begin{array}{c}
                          S \left( \prod\limits_{k=1}^l U(p) \right) \\[4mm]
                          \mbox{(reducible)} \\ \mbox{(direct sum)}
                         \end{array}$ \rule[-6ex]{0mm}{15ex} &
 $S_l$ & $\mbox{for}$ & $n = p\>\!l\,$, $l \geqslant 3\,$, $p \geqslant 1$ \\
 \rule[-8ex]{0mm}{17ex} $\begin{array}{c}
                          \prod\limits_{k=1}^l SU(p) \, / \, \Z_p^{l-1}
                          \\[4mm]
                          \mbox{(irreducible)} \\ \mbox{(tensor product)}
                        \end{array}$ \rule[-8ex]{0mm}{17ex} & 
 $\Z_{p\>\!^{l-1}} \times S_l$ & \mbox{for} &
 $n = p\>\!^{l}\,$, $l \geqslant 3\,$, $p \geqslant 2$
 \\ \hline
\end{tabular}
\end{center}
\begin{center}
\caption{\label{tab:SUG}
  Non-simple maximal subgroups $H$ of $SU(n)$ ($n \geqslant 2$)}
\end{center}
\end{table}

\subsection{Maximal subgroups of $SO(n)$}

\begin{enumerate}
 \item The inclusion
       \[
        SO(p) \times SO(q)~\subset~SO(n) \qquad (n = p+q)
       \]
       is given by the direct sum of the defining representations of
       $O(p)$ and $O(q)$. Explicitly, it is obtained as the restriction
       of the inclusion of $\, O(p) \times O(q) \,$ into $O(n)$ by block
       diagonal matrices. In fact, the inverse image of $\, SO(n) \subset
       O(n)$ under this inclusion is the subgroup
       \[
        H_+~=~S \bigl( O(p) \times O(q) \bigr)~
        =~\{ \, (A,B) \in O(p) \times O(q)~|~\det(A) \det(B) = 1 \, \}
       \]
       consisting of the block diagonal matrices belonging to $SO(n)$,
       which has two connected components: its one-component is the subgroup
       \[
        H_0~=~SO(p) \times SO(q)~
        =~\{ \, (A,B) \in O(p) \times O(q)~|~\det(A) = 1 = \det(B) \, \}~,
       \]
       while the other component is the set
       \[
        \{ \, (A,B) \in O(p) \times O(q)~|~\det(A) = -1 = \det(B) \, \}~.
       \]
       It follows from Lemma~\ref{lem:CNDS} that the centralizer $Z_G(H_0)$
       of~$H_0$ in $G$ is contained in $H_+\,$; more precisely, we have
       $\, H_i = H_+ \,$ if $p$ and $q$ are both odd and $\, H_i = H_0 \,$
       otherwise.
       Similarly, we compute $\, H = N_G(H_0) \,$ by combining Lemma~%
       \ref{lem:CNDS} and Lemma~\ref{lem:BLDIAG} to conclude that, as
       in the $SU(n)$ case, the elements $h$ of~$H$ must be of the form
       $\, h = \bigl( \begin{smallmatrix} A & 0 \\ 0 & D \end{smallmatrix}
       \bigr) \,$ or $\, h = \bigl( \begin{smallmatrix} 0 & B \\ C & 0
       \end{smallmatrix} \bigr)$, with $\, p = q \,$ in the second case.
       In~the first case, it follows that $\, h \in H_+$, and thus
       $\, H = H_+$, $H/H_0 = \mathbb{Z}_2 \,$ if $\, p \neq q$.
       In the second case, we may write $h$ as the product of an
       element of~$H_+$ with the matrix $\, \bigl( \begin{smallmatrix}
       0 & 1 \\ -1 & 0 \end{smallmatrix} \bigr) \in SO(n)$ \linebreak
       ($n = 2p$), whose square is $-1_n$ and belongs to~$H_+$
       but belongs to~$H_0$ only when $p$ is even, so we may con%
       clude that $\, H/H_0 = \mathbb{Z}_2 \times \mathbb{Z}_2 \,$
       if $\, p = q \,$ is even and $\, H/H_0 = \mathbb{Z}_4 \,$ if
       $\, p = q \,$ is odd.
 \item The inclusion
       \[
        SO(4)~\subset~SO(5)
       \]
       is a special case of the inclusion $\, SO(n-1) \subset SO(n) \,$
       and is the only one to appear here because $SO(n-1)$ is simple if
       $\, n \geqslant 6$.
       But whether simple or not, $\mathfrak{so}(n-1)$ is a maximal
       subalgebra of $\mathfrak{so}(n)$, and the normalizer of the
       corresponding connected subgroup $\; H_0 = \{ \bigl(
       \begin{smallmatrix} A & 0 \\ 0 & 1 \end{smallmatrix}
       \bigr) \, | \, A \in SO(n-1) \} \cong SO(n-1)$ \linebreak
       in~$SO(n)$ is the subgroup $\; H_+ = \{ \bigl(
       \begin{smallmatrix} A & 0 \\ 0 & \det A \end{smallmatrix}
       \bigr) \, | \, A \in O(n-1) \} \cong O(n-1)$.
       Therefore, $H/H_0 = \mathbb{Z}_2$.
 \item The inclusion
       \[
        U(p)~\subset~SO(n) \qquad (n = 2p)
       \]
       is given by
       \[
        A+iB~\longmapsto \left( \begin{array}{cc}
                                 A & B \\ -B & A
                                \end{array} \right).
       \]
       Note that $\mathfrak{u}(p)$ has maximal rank in $\mathfrak{so}(n)$
       and hence $\, H_i = H_0$.
       Moreover, $\mathrm{Aut}(\Gamma) = Z_2 \,$ where the non-trivial
       diagram automorphism is given by complex conjugation $\, A+iB
       \mapsto A-iB$, and this can be implemented by conjugating
       $\bigl( \begin{smallmatrix} A & B \\ -B & A \end{smallmatrix}
       \bigr)$ with the matrix $\bigl( \begin{smallmatrix} 1 & 0 \\ 0 & - 1
       \end{smallmatrix} \bigr)$, which belongs to $SO(n)$ if and only if
       $\, n = 2p \,$ with $p$ even. Therefore, $H/H_0 = \mathbb{Z}_2 \,$
       if $\, n = 2p \,$ with $p$ even and $\, H/H_0 = \{1\} \,$ if
       $\, n = 2p \,$ with $p$ odd.
 \item The inclusions
       \[
        \begin{array}{cc}
         SO(p) \times SO(q)~\subset~SO(n) \qquad &
         (n = p\>\!q \, , \mbox{$p$ or $q$ odd}) \\[2mm]
         SO(p) \times_{\mathbb{Z}_2} SO(q)~\subset~SO(n) \qquad &
         (n = p\>\!q \, , \mbox{$p$ and $q$ even})
        \end{array}
       \]
       are given by the tensor product of the defining representations of
       $O(p)$ and $O(q)$. Explicitly, they are induced by the homomorphism
       \[
        \begin{array}{ccc}
         O(p) \times O(q) & \longrightarrow &    O(n)     \\[1mm]
              (A,B)       &   \longmapsto   & A \otimes B
        \end{array}
       \]
       which, in contrast to the situation encountered in the first item,
       is not injective but rather has a non-trivial kernel given by
       \[
        \{ (1_p\,,1_q),(-1_p\,,-1_q) \}~\cong~\mathbb{Z}_2~.
       \]
       Using the formula
       \[
        \det(A \otimes B)~=~(\det(A))^q \, (\det(B))^p~,
       \]
       we see that the inverse image of $\, SO(n) \subset O(n) \,$ under this
       homomorphism is the subgroup
       \[
        S' \bigl( O(p) \times O(q) \bigr)~
        =~\left\{ \begin{array}{ccccc}
                   S \bigl( O(p) \times O(q) \bigr) &
                  \mbox{~if~} & \mbox{$p$ odd} &,& \mbox{$q$ odd} \\[1mm]
                   O(p) \times SO(q) &
                  \mbox{~if~} & \mbox{$p$ odd} &,& \mbox{$q$ even} \\[1mm]
                   SO(p) \times O(q) &
                  \mbox{~if~} & \mbox{$p$ even} &,& \mbox{$q$ odd} \\[1mm]
                   O(p) \times O(q) &
                  \mbox{~if~} & \mbox{$p$ even} &,& \mbox{$q$ even}
                 \end{array} \right.~.
       \]
       In the first three cases, the group $\, S' \bigl( O(p) \times O(q)
       \bigr) \,$ has two connected components such that the kernel of the
       homomorphism introduced above meets both components in exactly one
       element (since $(1_p\,,1_q)$ and $(-1_p\,,-1_q)$ belong to different
       connected components), so the restriction of this homo\-morphism
       to the connected one-component provides the desired inclusion.
       In~the~last case, the group $\, S' \bigl( O(p) \times O(q) \bigr) \,$
       has four connected components, but the kernel of the homomorphism
       introduced above is contained in the connected one-component
       (since $(1_p\,,1_q)$ and $(-1_p\,,-1_q)$ belong to the same
       connected component), so it is necessary to factor this out
       in order to obtain the following sequence of inclusions:
       \[
         SO(p) \times_{\mathbb{Z}_2} SO(q)~
         \subset~O(p) \times_{\mathbb{Z}_2} O(q)~\subset~SO(n)~.
       \]
       Note also that the center $Z(G)$ of $G$ is contained in~$H_0$;
       hence $H_i = H_0$. \\[2mm]
       In order to compute the normalizer $H$ of the connected subgroup
       $H_0$ so defined, we first note that (for $p$ even, $p \geqslant 4$)
       the Dynkin diagram of $\mathfrak{so}(p)$ admits only one non-trivial
       automorphism which can be implemented by a reflection $\sigma_p$
       in $\mathbb{R}^p$.%
       \footnote{The argument can also be applied in the case of~%
       $\mathfrak{so}(4)$, even though its Dynkin diagram is not connected.
       In the case of $\mathfrak{so}(8)$, there are other non-trivial
       automorphisms of the Dynkin diagram, but these cannot be
       implemented by linear maps on $\mathbb{R}^8$.}
       Computing determinants, we see that the auto\-morphisms of the Dynkin
       diagram $\Gamma$ of $\, \mathfrak{so}(p) \times \mathfrak{so}(q) \,$
       induced by auto\-morphisms of the Dynkin diagrams of the factors can
       always be extended to automorphisms of $\mathfrak{so}(n)$ and that
       these will be outer automorphisms (implemented by a matrix in $O(n)$
       that does not belong to $SO(n)$) if $p$ or $q$ is odd but will be
       inner automorphisms (implemented by a matrix in $SO(n)$) if $p$
       and $q$ are even.
       Thus it becomes clear that for $p \neq q$,%
       \[
        \mathrm{Aut}(\Gamma) \cap \mathrm{Inn}(\mathfrak{g})~
        =~\left\{ \begin{array}{ccc}
                        \{1\}       & \mbox{if} &
                   \mbox{$p \neq q$, $p$ or $q$ odd} \\[1mm]
                   \Z_2 \times \Z_2 & \mbox{if} &
                   \mbox{$p \neq q$, $p$ and $q$ even}
                  \end{array} \right.~,
       \]
       while for $p = q$, $\mathrm{Aut}(\Gamma) \cap \mathrm{Inn}%
       (\mathfrak{g}) \,$ will contain an additional factor, corres%
       ponding to the possibility of switching factors in the tensor
       product that exists in this case.
       Explicitly, a representative of the corresponding connected
       component of $H$ is given by the transformation
       \[
        \begin{array}{ccc}
          \mathbb{R}^n   & \longrightarrow & \mathbb{R}^n \\[1mm]
         z = x \otimes y &   \longmapsto   & z^\tau = \pm \, y \otimes x
        \end{array}~,
       \]
       which, by the same argument employed in the $SU(n)$ case, is
       an involution with determinant $(-1)^{p(p \mp 1)/2}$ and thus
       belongs to $SO(n)$, provided we make~an adequate choice of
       sign, except when $\, p = 2 \!\mod 4$.
       Therefore, we conclude that $\, \mathrm{Aut}(\Gamma) \cap
       \mathrm{Inn}(\mathfrak{g}) \,$ is equal to $\mathbb{Z}_2$
       if $p$ is odd, to $\, \mathbb{Z}_2 \times \mathbb{Z}_2$
       \linebreak
       if $\, p = 2 \!\mod 4$, to $\, (\mathbb{Z}_2 \times \mathbb{Z}_2)
       : \mathbb{Z}_2$ (or $\mathbb{Z}_2 \ltimes (\mathbb{Z}_2 \times
       \mathbb{Z}_2)$) if $\, p = 0 \!\mod 4$ \linebreak with $\, p > 4 \,$
       and to the full permutation group $S_4$ if $\, p = 4 \,$
       (this case is also covered in item 8 below).
       Note that, in the penultimate case, the structure of this
       group as a semi-direct product is given by conjugation with
       the transformation $.{}^\tau$ which maps $\, A \otimes B \,$
       to $\, B \otimes A \,$ and thus acts on $\, \mathbb{Z}_2
       \times \mathbb{Z}_2$ \linebreak by switching the factors.
 \item The inclusion
       \[
        \SPG(2p) \times_{\mathbb{Z}_2} \SPG(2q)~\subset~SO(n)
        \qquad (n=4p\>\!q)
       \]
       is given by the tensor product of the defining representations
       of $\SPG(2p)$ and $\SPG(2q)$.
       Indeed, noting that the tensor product of two pseudo-real/
       quaternionic representations is real, we can proceed in the
       same way as in the previous item, though with some
       simplifications: the homomorphism
       \[
        \begin{array}{ccc}
         \SPG(2p) \times \SPG(2q) & \longrightarrow &    O(n)     \\[1mm]
                  (A,B)           &   \longmapsto   & A \otimes B
        \end{array}
       \]
       always has image contained in $\, G = SO(n) \,$ and the same kernel
       as before:
       \[
        \{ (1_{2p}\,,1_{2q}),(-1_{2p}\,,-1_{2q}) \}~\cong~\mathbb{Z}_2~.
       \]
       Therefore, it provides the desired inclusion of the corresponding
       quotient group.
       Note also that the center $Z(G)$ of $G$ is contained in $H_0\,$;
       hence $H_i = H_0$. \\[2mm]
       In order to compute the normalizer $H$ of the connected
       subgroup $H_0$ so defined, we note that for $\, p \neq q$,
       $\, \mathrm{Aut}(\Gamma) \cap \mathrm{Inn}(\mathfrak{g}) = \{1\}$,
       while for $\, p = q$, $\mathrm{Aut}(\Gamma) \cap \mathrm{Inn}%
       (\mathfrak{g}) \,$ can at most be equal to the group $\mathbb{Z}_2$,
       corresponding to the possibility of switching the factors in the
       tensor product that exists in this case.
       The argument to decide whether this switch operator has determinant
       $+1$ or $-1$ is the same as before and leads to the conclusion that
       $\, \mathrm{Aut}(\Gamma) \cap \mathrm{Inn}(\mathfrak{g}) \,$ is
       trivial if $\, p = q \,$ is odd and is equal to $\mathbb{Z}_2$
       if $\, p = q \,$ is even.
 \item Generalizing the procedure of item 1, we define the inclusion
       \[
        SO(p) \times \ldots \times SO(p)~\subset~SO(n) \qquad (n=p\>\!l)
       \]
       by the direct sum of $l$ copies of the defining representation of
       $O(p)$, which realizes $\, O(p) \times \ldots \times O(p) \,$ by
       block diagonal matrices in $O(n)$, with $l$ blocks of size $p$
       along the diagonal. The inverse image of $\, SO(n) \subset O(n) \,$
       under this inclusion is the subgroup
       \[
        H_+~=~S \bigl( O(p) \times \ldots \times O(p) \bigr)~
       \]
       consisting of the block diagonal matrices belonging to $SO(n)$,
       which has $2^{l-1}$ connected components, its one-component being
       the subgroup
       \[
        H_0~=~SO(p) \times \ldots \times SO(p)~.
       \]
       As before, it follows from Lemma~\ref{lem:CNDS} that the centralizer
       $Z_G(H_0)$ of~$H_0$ in $G$ is~contained in $H_+\,$; more precisely,
       we have $\, H_i = H_+ \,$ if $p$ is odd and $\, H_i = H_0 \,$
       if $p$ is even.
       Similarly, we compute $\, H = N_G(H_0) \,$ by combining Lemma~%
       \ref{lem:CNDS} and Lemma~\ref{lem:BLDIAG} to conclude that, as
       in the $SU(n)$ case, the elements $h$ of~$H$ must be such that,
       when represented as block $(l \times l)$-matrices (with entries
       that are themselves $(p \times p)$-matrices), they contain
       precisely one nonvanishing entry in each line and each column:
       obviously, any such matrix $h$ defines a permutation $\sigma(h)$
       of $\{1,\ldots,l\}$. Thus we obtain a group homomorphism from~$H$
       to the permutation group~$S_l$ which has kernel~$H_+$ and is
       surjective (we can use the same argument as in item~1 above
       to show that its image contains all transpositions, their
       preimages being matrices whose square is $-1_p$ in two of
       the $l$ diagonal blocks and $1_p$ in the remaining $l-2$
       ones, so these belong to~$H_+$ but belong to~$H_0$ only
       when $p$ is even), so we may conclude that $H/H_0$ is an
       upwards extension of the group $\mathbb{Z}_2^{l-1}$ by the
       symmetric group $S_l$, which is a split extension, that is,
       a semi-direct product $\, H/H_0 = \mathbb{Z}_2^{l-1} : S_l$
       (or $S_l \ltimes \mathbb{Z}_2^{l-1}$), if $p$ is even and is
       a non-split extension $\, H/H_0 = \mathbb{Z}_2^{l-1} .\, S_l \,$
       if $p$ is odd.
 \item Generalizing the procedure of item 4, we define the inclusions
       \[
        \begin{array}{ccc}
         SO(p) \times \ldots \times SO(p)~\subset~SO(n) \qquad &
         (n = p^l \, , \, \mbox{$p$ odd}) \\[2mm]
         \bigl( SO(p) \times \ldots \times SO(p) \bigr) /\mathbb{Z}_2^{l-1}~
         \subset~SO(n) \qquad &
         (n = p^l \, , \, \mbox{$p$ even})
        \end{array}
       \]
       by the tensor product of $l$ copies of the defining representation
       of $O(p)$. Explicitly, they are induced by the homomorphism
       \[
        \begin{array}{ccc}
         O(p) \times \ldots \times O(p) & \longrightarrow &
         O(n)     \\[1mm]
               (A_1,\ldots,A_l)         &   \longmapsto   &
         A_1 \otimes \ldots \otimes A_l
        \end{array}
       \]
       which, once again, is not injective but rather has a non-trivial kernel
       given by
       \[
        \{ \, (\epsilon_1 1_p\,,\ldots,\epsilon_l 1_p)~|~
              \epsilon_1,\ldots,\epsilon_l = \pm 1 \,,\,
              \epsilon_1 \ldots\, \epsilon_l = 1 \, \}~
        \cong~\mathbb{Z}_2^{l-1}~.
       \]
       Using the formula
       \[
        \det \bigl( A_1 \otimes \ldots \otimes A_l \bigr)~
        =~\bigl( \det(A_1) \,\ldots\, \det(A_l) \bigr)^{p^{\,l-1}}~,
       \]
       we see that the inverse image of $\, SO(n) \subset O(n) \,$ under
       this homomorphism is the subgroup
       \[
        S' \bigl( O(p) \times \ldots \times O(p) \bigr)
        =~\left\{ \begin{array}{ccc}
                   S \bigl( O(p) \times \ldots \times O(p) \bigr) &
                   \mbox{if} & \mbox{$p$ odd} \\[2mm]
                   O(p) \times \ldots \times O(p) &
                   \mbox{if} & \mbox{$p$ even}
                 \end{array} \right.~.
       \]
       In the first case, the group $\, S' \bigl( O(p) \times \ldots \times
       O(p) \bigr) \,$ has $2^{l-1}$ connected components such that the
       kernel of the homomorphism introduced above meets each of them
       in exactly one element, so the restriction of this homomorphism
       to the connected one-component provides the desired inclusion.
       In the second case, the group $\, S' \bigl( O(p) \times \ldots
       \times O(p) \bigr) \,$ has $2^l$ connected components, but the
       kernel of the homomorphism introduced above is contained in the
       connected one-component, so it is necessary to factor this out
       in order to obtain the following sequence of inclusions:
       \[
         \bigl( SO(p) \times \ldots \times SO(p) \bigr) / \mathbb{Z}_2^{l-1}~
         \subset~
         \bigl( O(p) \times \ldots \times O(p) \bigr) / \mathbb{Z}_2^{l-1}~
         \subset~SO(n)~.
       \]
       Note also that the center $Z(G)$ of $G$ is contained in $H_0\,$;
       hence $H_i = H_0$.

       In order to compute the normalizer $H$ of the connected subgroup
       $H_0$ so defined, we note, as before, that every
       non-trivial automorphism of the Dynkin diagram $\Gamma$ 
       of $\, \mathfrak{so}(p) \times \ldots \times \mathfrak{so}(p) \,$
       induced by non-trivial automorphisms of the Dynkin 
       diagrams of the factors (of course these exist only when $p$ is
       even) can be extended to an inner automorphism of $\mathfrak{so}(n)$.
       Moreover, $\, \mathrm{Aut}(\Gamma) \cap \mathrm{Inn}(\mathfrak{g}) \,$
       will contain a permutation group as an additional factor, corresponding
       to the possibility of switching the factors in the tensor product.
       \linebreak
       But for $\, l \geqslant 3$, we can proceed as in the $SU(n)$ case
       to implement any transposition through an involution of 
       determinant $+1$.
       Therefore we conclude that $\, \mathrm{Aut}(\Gamma) \cap
       \mathrm{Inn}(\mathfrak{g}) \,$ is equal to the permutation
       group $S_l$ if $p$ is odd, to $\, \Z_2^l : S_l \,$ (or
       $\, S_l \ltimes \Z_2^l$) if $p$ is even with $\, p > 4 \,$
       and to the permutation group $S_{2l}$ if $\, p = 4 \,$ (this
       case is also covered in item 8 below).
       In the penultimate case, the structure of this group is that
       of a semi-direct product since $S_l$ acts on $\mathbb{Z}_2^l$
       by switching the factors.
 \item Generalizing the procedure of item 5, we define the inclusion
       \[
         \bigl( \SPG(2p) \times \ldots \times \SPG(2p) \bigr) /
         \mathbb{Z}_2^{l-1}~\subset~SO(n)
         \qquad (n=(2p)^l, \, \mbox{$l$ even})
       \]
       by the tensor product of $l$ copies of the defining representation
       of $\SPG(2p)$. \linebreak
       Indeed, noting that the tensor product of $l$ copies of a
       pseudo-real/ \linebreak quaternionic representation is real
       if $l$ is even, we can proceed as in the previous item,
       though with some simplifications: the homomorphism
       \[
        \begin{array}{ccc}
         \SPG(2p) \times \ldots \times \SPG(2p) & \longrightarrow &
         O(n)     \\[1mm]
                   (A_1,\ldots,A_l)             &   \longmapsto   &
         A_1 \otimes \ldots \otimes A_l
        \end{array}
       \]
       always has image contained in $\, G = SO(n) \,$ and has the same kernel
       as before:
       \[
        \{ \, (\epsilon_1 1_{2p}\,,\ldots,\epsilon_l 1_{2p})~|~
              \epsilon_1,\ldots,\epsilon_l = \pm 1 \,,\,
              \epsilon_1 \ldots\, \epsilon_l = 1 \, \}~
        \cong~\mathbb{Z}_2^{l-1}~.
       \]
       Therefore, it provides the required inclusion of the corresponding
       quotient group.
       Note also that the center $Z(G)$ of $G$ is contained in $H_0\,$;
       hence $H_i = H_0$.

       In order to compute the normalizer $H$ of the connected subgroup
       $H_0$ so defined, we~note, as before, that $\, \mathrm{Aut}(\Gamma)
       \cap \mathrm{Inn}(\mathfrak{g}) \,$ can at most be equal to the
       symmetric group $S_l$, corresponding to the possibility of
       permuting the factors in the tensor product.
       But for $\, l \geqslant 3$, we can proceed as in the $SU(n)$ case to
       implement any transposition through an involution of determinant $+1$.
       Therefore, we conclude that $\, \mathrm{Aut}(\Gamma) \cap \mathrm{Inn}%
       (\mathfrak{g}) = S_l$.
\end{enumerate}

\begin{table}[!htb]
\begin{center}
\begin{tabular}{|c|c@{\hspace{0.5em}}c@{\hspace{0.5em}}c|} \hline
 \rule[-2mm]{0mm}{7mm} {\small connected component $H_0$}
 \rule[-2mm]{0mm}{7mm} &
 \multicolumn{3}{c|}{{\small component group $H/H_0$}} \\ \hline\hline
 \rule[-5ex]{0mm}{13ex} $\begin{array}{c}
                          SO(p) \times SO(q) \\[1mm]
                          \mbox{(reducible)} \\ \mbox{(direct sum)}
                         \end{array}$ \rule[-5ex]{0mm}{13ex} &
 $\begin{array}{c} \Z_2 \\[1mm] \Z_4 \\[1mm] \Z_2 \times \Z_2 \end{array}$ &
 $\begin{array}{c} \mbox{for} \\[1mm] \mbox{for} \\[1mm] \mbox{for}
  \end{array}$ &
 $\begin{array}{c}
   n = p+q\,,\, p > q \geqslant 2 \\[1mm]
   n = p+q\,,\, p = q \geqslant 3~\mbox{odd} \\[1mm]
   n = p+q\,,\, p = q \geqslant 4~\mbox{even}
  \end{array}$ \\
 \rule[-2ex]{0mm}{6ex} $SO(4)$ \rule[-2ex]{0mm}{6ex} &
 $\Z_2$ & \mbox{for} & $n=5$ \\
 \rule[-4ex]{0mm}{10ex} $U(p)$ \rule[-4ex]{0mm}{10ex} &
 $\begin{array}{c} \{1\} \\[1mm] \Z_2 \end{array}$ &
 $\begin{array}{c} \mbox{for} \\[1mm] \mbox{for} \end{array}$ &
 $\begin{array}{c}
   n = 2p\,,\, p \geqslant 3\,,\, \mbox{$p$ odd} \\[1mm]
   n = 2p\,,\, p \geqslant 4\,,\, \mbox{$p$ even}
  \end{array}$ \\
 \rule[-7ex]{0mm}{16ex} $\begin{array}{c}
                         SO(p) \times SO(q) \\[1mm]
                         \mbox{(irreducible)} \\ \mbox{(tensor product)}
                        \end{array}$ \rule[-3ex]{0mm}{9ex} &
 $\begin{array}{c} \{1\} \\ \mbox{} \\[1mm] \Z_2 \end{array}$ &
 $\begin{array}{c} \mbox{for} \\ \mbox{} \\[1mm] \mbox{for} \end{array}$ &
 $\begin{array}{c}
   n = p\>\!q\,,\, p > q \geqslant 3 \\
   \mbox{$p$ or $q$ odd} \\[1mm]
   n = p\>\!q\,,\, p = q \geqslant 3 \\
   \mbox{$p=q$ odd}
  \end{array}$ \\
 \rule[-10ex]{0mm}{22ex} $\begin{array}{c}
                         SO(p) \times_{\Z_2} SO(q) \\[1mm]
                         \mbox{(irreducible)} \\ \mbox{(tensor product)}
                        \end{array}$ \rule[-3ex]{0mm}{9ex} &
 $\begin{array}{c}
   \Z_2 \times \Z_2 \\ \mbox{} \\[1mm]
   \Z_2 \times \Z_2 \\ \mbox{} \\[1mm]
   (\Z_2 \times \Z_2) : \Z_2
  \end{array}$ &
 $\begin{array}{c}
   \mbox{for} \\ \mbox{} \\[1mm] \mbox{for} \\ \mbox{} \\[1mm] \mbox{for}
  \end{array}$ &
 $\begin{array}{c}
   n = p\>\!q\,,\, p > q \geqslant 4 \\
   \mbox{$p$ and $q$ even} \\[1mm]
   n = p\>\!q\,,\, p = q \geqslant 6 \\
   \mbox{$p=q=2 \mod 4$} \\[1mm]
   n = p\>\!q\,,\, p = q \geqslant 8 \\
   \mbox{$p=q=0 \mod 4$}
  \end{array}$ \\
 \rule[-6ex]{0mm}{14ex} $\begin{array}{c}
                         \SPG(2p) \times_{\Z_2} \SPG(2q) \\[1mm]
                         \mbox{(irreducible)} \\ \mbox{(tensor product)}
                        \end{array}$ \rule[-4ex]{0mm}{10ex} &
 $\begin{array}{c} \{1\} \\[1mm] \{1\} \\[1mm] \Z_2 \end{array}$ &
 $\begin{array}{c}
   \mbox{for} \\[1mm] \mbox{for} \\[1mm] \mbox{for}
  \end{array}$ &
 $\begin{array}{c}
   n = 4p\>\!q\,,\, \mbox{$p > q \geqslant 1$} \\[1mm]
   n = 4p\>\!q\,,\, \mbox{$p = q$ odd} \\[1mm]
   n = 4p\>\!q\,,\, \mbox{$p = q$ even}
  \end{array}$ \\ \hline
\end{tabular}
\end{center}
\begin{center}
 \caption{\label{tab:SOG1}
          Non-simple maximal subgroups $H$ of $SO(n)$ ($n \geqslant 5$),
          part 1}
\end{center}
\end{table}

\begin{table}[!htb]
\begin{center}
\begin{tabular}{|c|c@{\hspace{0.5em}}c@{\hspace{0.5em}}c|} \hline
 \rule[-2mm]{0mm}{7mm} {\small connected component $H_0$}
 \rule[-2mm]{0mm}{7mm} &
 \multicolumn{3}{c|}{{\small component group $H/H_0$}} \\ \hline\hline
 \rule[-6ex]{0mm}{15ex} $\begin{array}{c}
                         \prod\limits_{k=1}^l SO(p) \\[4mm]
                         \mbox{(reducible)} \\ \mbox{(direct sum)}
                        \end{array}$ \rule[-3ex]{0mm}{9ex} &
 $\begin{array}{c} \Z_2^{l-1} .\, S_l \\[1mm]
                   \Z_2^{l-1} : S_l \end{array}$ &
 $\begin{array}{c} \mbox{for} \\[1mm] \mbox{for} \end{array}$ &
 $\begin{array}{c}
   n = p\>\!l\,,\, l \geqslant 3\,,\, p \geqslant 3\,,\, \mbox{$p$ odd}
   \\[1mm]
   n = p\>\!l\,,\, l \geqslant 3\,,\, p \geqslant 2\,,\, \mbox{$p$ even}
  \end{array}$ \\
 \rule[-6ex]{0mm}{15ex} $\begin{array}{c}
                         \prod\limits_{k=1}^l SO(p) \\[4mm]
                         \mbox{(irreducible)} \\ \mbox{(tensor product)}
                        \end{array}$ \rule[-3ex]{0mm}{9ex} &
 $S_l$ & $\mbox{for}$ &
 $n = p\>\!^l,\, l \geqslant 3\,,\, p \geqslant 3\,,\, \mbox{$p$ odd}$
 \\
 \rule[-6ex]{0mm}{15ex} $\begin{array}{c}
                         \prod\limits_{k=1}^l SO(p) \, / \, \Z_2^{l-1}
                         \\[4mm]
                         \mbox{(irreducible)} \\ \mbox{(tensor product)}
                        \end{array}$ \rule[-3ex]{0mm}{9ex} &
 $\Z_2^l : S_l$ & $\mbox{for}$ &
 $n = p\>\!^l,\, l \geqslant 3\,,\, p \geqslant 6\,,\, \mbox{$p$ even}$ \\
 \rule[-7ex]{0mm}{16ex} $\begin{array}{c}
                         \prod\limits_{k=1}^l \SPG(2p) \, / \, \Z_2^{l-1}
                         \\[4mm]
                         \mbox{(irreducible)} \\ \mbox{(tensor product)}
                        \end{array}$ \rule[-3ex]{0mm}{9ex} &
 $S_l$ & $\mbox{for}$ &
 $n=(2p)^l,\, l \geqslant 4\,,\, \mbox{$l$ even}, p \geqslant 1$
 \\ \hline
\end{tabular}
\end{center}
\begin{center}
 \caption{\label{tab:SOG2}
          Non-simple maximal subgroups $H$ of $SO(n)$ ($n \geqslant 5$),
          part 2}
\end{center}
\end{table}

\subsection{Maximal subgroups of $\SPG(n)$ ($n$ even)}

\begin{enumerate}
 \item The inclusion
       \[
        \SPG(2p) \times \SPG(2q)~\subset~\SPG(n) \qquad (n = 2(p+q))
       \]
       is given by the direct sum of the defining representations
       of $\SPG(2p)$ and $\SPG(2q)$. Again, it is obtained by taking
       block diagonal matrices, and the inverse image of $\SPG(n)$
       under this inclusion is the subgroup \linebreak $\, H_0 = \SPG(2p)
       \times \SPG(2q) \,$ consisting of the block diagonal
       matrices belonging to $\SPG(n)$, which is connected.
       The remainder of the argument is the same as in item~1
       of the $SU(n)$ case, leading to the conclusion that
       $\, H = H_0 \,$ if $\, p \neq q \,$ and $\, H/H_0 =
       \mathbb{Z}_2 \,$ if $\, p = q$.
 \item The inclusion
       \[
        U(p)~\subset~\SPG(n) \qquad (n = 2p)
       \]
       is given by
       \[
        A+iB~\longmapsto \left( \begin{array}{cc}
                                 A & B \\ -B & A
                                \end{array} \right).
       \]
       Note that $\mathfrak{u}(p)$ has maximal rank in $\mathfrak{sp}(n)$
       and hence $\, H_i = H_0$. Moreover, $\mathrm{Aut}(\Gamma) = Z_2 \,$
       where the non-trivial diagram automorphism is given by complex con%
       jugation $\, A+iB \mapsto A-iB$, and this can be implemented by con%
       jugating $\bigl( \begin{smallmatrix} A & B \\ -B & A \end{smallmatrix}
       \bigr)$ with the matrix $\bigl( \begin{smallmatrix} 1 & 0 \\ 0 & - 1
       \end{smallmatrix} \bigr)$, which belongs to $\SPG(n)$. Therefore,
       $H/H_0 = \mathbb{Z}_2$.
 \item The inclusions
       \[
        \begin{array}{cc}
         \SPG(2p) \times SO(q)~\subset~\SPG(n) \qquad &
         (n = 2p\>\!q \, , \, \mbox{$q$ odd}) \\[2mm]
         \SPG(2p) \times_{\mathbb{Z}_2} SO(q)~\subset~\SPG(n) \qquad &
         (n = 2p\>\!q \, , \, \mbox{$q$ even})
        \end{array}
       \]
       are given by the tensor product of the defining representations
       of $\SPG(2p)$ and $O(q)$.
       Indeed, noting that the tensor product of a pseudo-real/quaternionic
       representation and a real representation is pseudo-real/quaternionic,
       we can proceed in the same way as in items 4 and 5 of the $SO(n)$ case:
       the homomorphism
       \[
        \begin{array}{ccc}
         \SPG(2p) \times O(q) & \longrightarrow &   \SPG(n)   \\[1mm]
                (A,B)         &   \longmapsto   & A \otimes B
        \end{array}
       \]
       always has image contained in $\, G = \SPG(n) \,$ and has the same
       kernel as before:
       \[
        \{ (1_{2p}\,,1_q),(-1_{2p}\,,-1_q) \}~\cong~\mathbb{Z}_2~.
       \]
       On the other hand, the group $\, \SPG(2p) \times O(q) \,$ has two
       connected components.
       If~$q$~is odd, the kernel of the above homomorphism meets
       both components in exactly one element (since $(1_{2p}\,,1_q)$ and
       $(-1_{2p}\,,-1_q)$ belong to different connected components), so the
       restriction of this homomorphism to the connected one-component
       provides the desired inclusion.
       If $q$ is even, the kernel of the above homomorphism is contained
       in the connected one-component (since $(1_{2p}\,,1_q)$ and
       $(-1_{2p}\,,-1_q)$ belong to the same connected component),
       so it is necessary to factor this out in order to obtain
       the following sequence of inclusions:
       \[
         \SPG(2p) \times_{\mathbb{Z}_2} SO(q)~
         \subset~\SPG(2p) \times_{\mathbb{Z}_2} O(q)~\subset~\SPG(n)~.
       \]
       Note also that the center $Z(G)$ of $G$ is contained in~$H_0\,$;
       hence $\, H_i = H_0$.

       In order to compute the normalizer $H$ of the connected subgroup
       $H_0$ so defined, we can proceed in the same way as in items 4 and 5
       of the $SO(n)$ case to conclude that $\, \mathrm{Aut}(\Gamma)
       \cap \mathrm{Inn}(\mathfrak{g}) \,$ is trivial if $q$ is odd\,%
       \footnote{Here, we have to exclude the possibility of a switch
       operator between the tensor factors when $p=2$ and $q=5$, since
       $\, \mathfrak{sp}(4) \cong \mathfrak{so}(5)$, but this follows
       from Lemma~\ref{lem:CNTP}, taking into account that the two
       fundamental representations of this simple Lie algebra used
       to construct the embedding into $\mathfrak{sp}(20)$ are not
       quasiequivalent, since they have different dimensions ($4$
       and $5$, respectively).}
       and is equal to $\mathbb{Z}_2$ if $q$ is even, except when
       $\, p = 1 \,$ and $\, q = 4$, in which case it is equal to
       the full permutation group $S_3$ (this case is also covered
       in item 5 below).
 \item Generalizing the procedure of the first item, we define the inclusion
       \[
        \SPG(2p) \times \ldots \times \SPG(2p)~\subset~\SPG(n)
        \qquad (n = 2p\>\!l)
       \]
       by the direct sum of $l$ copies of the defining representation
       of $\SPG(2p)$.
       Again, it is obtained by taking block diagonal matrices, and the
       inverse image of $\SPG(n)$ under this inclusion is the subgroup
       $\, H_0 = \SPG(2p) \times \ldots \times \SPG(2p)$ \linebreak
       consisting of the block diagonal matrices belonging to $\SPG(n)$,
       which is connected.
       The remainder of the argument is the same as in item~3 of the
       $SU(n)$ case, leading to the conclusion that $\, H/H_0 = S_l$.
 \item Generalizing the procedure of previous items, we define the inclusion
       \[
         \bigl( \SPG(2p) \times \ldots \times \SPG(2p) \bigr) /
         \mathbb{Z}_2^{l-1}~\subset~\SPG(n)
         \qquad (n = (2p)^l, \, \mbox{$l$ odd})
       \]
       by the tensor product of $l$ copies of the defining representation
       of $\SPG(2p)$.
\end{enumerate}

\begin{table}[!htb]
\begin{center}
\begin{tabular}{|c|c@{\hspace{0.5em}}c@{\hspace{0.5em}}c|} \hline
 \rule[-2mm]{0mm}{7mm} {\small connected component $H_0$}
 \rule[-2mm]{0mm}{7mm} &
 \multicolumn{3}{c|}{{\small component group $H/H_0$}} \\ \hline\hline
 \rule[-5ex]{0mm}{13ex} $\begin{array}{c}
                          \SPG(2p) \times \SPG(2q) \\[1mm]
                          \mbox{(reducible)} \\ \mbox{(direct sum)}
                         \end{array}$ \rule[-3ex]{0mm}{9ex} &
 $\begin{array}{c} \{1\} \\[1mm] \Z_2 \end{array}$ &
 $\begin{array}{c} \mbox{for} \\[1mm] \mbox{for} \end{array}$ &
 $\begin{array}{c}
   n = 2(p+q)\,,\, p > q \geqslant 1 \\[1mm]
   n = 2(p+q)\,,\, p = q \geqslant 1
  \end{array}$ \\
 \rule[-2ex]{0mm}{6ex} $U(p)$ \rule[-1ex]{0mm}{5ex} &
 $\Z_2$ & $\mbox{for}$ & $n = 2p\,,\, p \geqslant 2$ \\
 \rule[-5ex]{0mm}{12ex} $\begin{array}{c}
                          \SPG(2p) \times SO(q) \\[1mm]
                          \mbox{(irreducible)} \\ \mbox{(tensor product)}
                         \end{array}$ \rule[-4ex]{0mm}{10ex} &
 $\{1\}$ & $\mbox{for}$ &
 $n = 2p\>\!q\,,\, p \geqslant 1\,,\, q \geqslant 3\,,\, \mbox{$q$ odd}$ \\
 \rule[-6ex]{0mm}{13ex} $\begin{array}{c}
                          \SPG(2p) \times_{\Z_2} SO(q) \\[1mm]
                          \mbox{(irreducible)} \\ \mbox{(tensor product)}
                         \end{array}$ \rule[-4ex]{0mm}{10ex} &
 $\Z_2$ & $\mbox{for}$ &
 $\begin{array}{c}
   n = 2p\>\!q\,,\, p \geqslant 1\,,\, q \geqslant 4\,,\, \mbox{$q$ even}
   \\[1mm]
   \mbox{except when $p=1$ and $q=4$}
  \end{array}$
 \\ \hline
 \rule[-6ex]{0mm}{15ex} $\begin{array}{c}
                          \prod\limits_{k=1}^l \SPG(2p) \\[4mm]
                          \mbox{(reducible)} \\ \mbox{(direct sum)}
                         \end{array}$ \rule[-4ex]{0mm}{11ex} &
 $S_l$ & $\mbox{for}$ & $n = 2p\>\!l\,,\, l \geqslant 3\,,\, p \geqslant 1$ \\
 \rule[-7ex]{0mm}{16ex} $\begin{array}{c}
                          \prod\limits_{k=1}^l \SPG(2p) \, / \, \Z_2^{l-1}
                          \\[4mm]
                          \mbox{(irreducible)} \\ \mbox{(tensor product)}
                         \end{array}$ \rule[-3ex]{0mm}{9ex} &
 $S_l$ & $\mbox{for}$ &
 $n = (2p)^l\,,\, l \geqslant 3\,,\, \mbox{$l$ odd}, p \geqslant 1$ \\ \hline
\end{tabular}
\end{center}
\begin{center}
 \caption{\label{tab:SPG}
   Non-simple maximal subgroups $H$ of $\SPG(n)$ ($n$ even, $n \geqslant 4$)}
\end{center}
\end{table}

\begin{enumerate} \item[]
       Indeed, noting that the tensor product of $l$ copies of a
       pseudo-real/ \linebreak quaternionic representation is
       pseudo-real/quaternionic if $l$ is odd, we can proceed
       as in item 8 of the $SO(n)$ case: the homomorphism
       \[
        \begin{array}{ccc}
         \SPG(2p) \times \ldots \times \SPG(2p) & \longrightarrow &
         \SPG(n)     \\[1mm]
                   (A_1,\ldots,A_l)             &   \longmapsto   &
         A_1 \otimes \ldots \otimes A_l
        \end{array}
       \]
       always has image contained in $\, G = \SPG(n) \,$ and the same kernel
       as before:
       \[
        \{ \, (\epsilon_1 1_{2p}\,,\ldots,\epsilon_l 1_{2p})~|~
              \epsilon_1,\ldots,\epsilon_l = \pm 1 \,,\,
              \epsilon_1 \ldots\, \epsilon_l = 1 \, \}~
        \cong~\mathbb{Z}_2^{l-1}~.
       \]
       Therefore, it provides the desired inclusion of the corresponding
       quotient group.
       Note also that the center $Z(G)$ of $G$ is contained in $H_0\,$;
       hence $\, H_i = H_0$.

 \pagebreak

       In order to compute the normalizer $H$ of the connected subgroup
       $H_0$ so defined, we note that $\, \mathrm{Aut}(\Gamma) \cap
       \mathrm{Inn}(\mathfrak{g}) \,$ can at most be equal to the
       symmetric group $S_l$, corresponding to the possibility of
       switching the factors in the tensor product.
       But for $\, l \geqslant 3$, we can, once more, proceed as in
       the $SU(n)$ case to implement any transposition through an
       involution of determinant $+1$. Therefore, we conclude that
       $\, \mathrm{Aut}(\Gamma) \cap \mathrm{Inn}(\mathfrak{g}) = S_l$.
\end{enumerate}

\section*{Acknowledgements}

This work has been partly motivated by the attempt to gain a thorough
understanding of the original work of Dynkin~\cite{Dy1,Dy2}, which led
to the master's thesis of the first author~\cite{Ant}, and grew out of
the master's thesis of the third author~\cite{Gav}, both under the
supervision of the second author.

The authors would like to thank one of the referees for his extensive
and careful report on a first version of this paper, containing a wealth
of constructive criticisms and suggestions that have helped to greatly
improve the manuscript. They would also like to thank the editor,
Prof.\ E.B.~Vinberg, for pointing out an important flaw in an
earlier version and providing the opportunity to correct it.


\end{document}